\documentclass{article}
\usepackage{geometry}    
\usepackage{amsfonts}    
\usepackage{amssymb}     
\usepackage[all]{xy}     
\usepackage{amsmath}     
\usepackage[dvips]{graphicx}       
\usepackage{color}
\usepackage{pst-tree}    
\usepackage{pstricks}

\newcommand{\go}[1]{\mathfrak {#1}}
\newcommand{\rrbr}{]\!]}
\newcommand{\llbr}{[\![}
\newcommand{\Symbole}{\begin{picture}(0,4)(0,0)
\qbezier(0,0)(0,0)(0,7) \qbezier(0,0)(0,0)(7,0)
 \qbezier(0,0)(0,0)(6,6)
 \qbezier(-0.5,7)(-0.5,7)(0.5,7) \qbezier(7,-0.5)(7,-0.5)(7,0.5)  \qbezier(5.5,6.5)(5.5,6.5)(6.5,5.5)
\end{picture}\hspace{9pt}}

\newtheorem{The}{Theorem}
\newtheorem{Def}[The]{Definition}

\newtheorem{Pro}[The]{Proposition}

\newtheorem{Rk}[The]{Remark}
\newtheorem{Ex}[The]{Example}

\begin{document}

\begin{center}
\Large{McKay correspondence and the branching law\\
for finite subgroups of $\mathbf{SL}_3\mathbb{C}$}\
\end{center}

\vspace{.4cm}
\begin{center}
Fr\'ed\'eric BUTIN\footnote{Universit\'e de Lyon, Universit\'e
Lyon 1, CNRS, UMR5208, Institut
Camille Jordan,
43 blvd du 11 novembre 1918, F-69622 Villeurbanne-Cedex, France,
email: butin@math.univ-lyon1.fr
}, Gadi S. PERETS\footnote{Universit\'e de Lyon, Universit\'e
Lyon 1, CNRS, UMR5208, Institut
Camille Jordan,
43 blvd du 11 novembre 1918, F-69622 Villeurbanne-Cedex, France,
email: gadi@math.univ-lyon1.fr
}\vspace{.5cm}\\
\end{center}

\begin{small}
\textbf{\textsc{Abstract}}\\
Given $\Gamma$ a finite subgroup of $\mathbf{SL}_3\mathbb{C}$, we determine how an arbitrary
finite dimensional irreducible representation of $\mathbf{SL}_3\mathbb{C}$ decomposes under the action of $\Gamma$.
To the subgroup $\Gamma$ we attach a generalized Cartan matrix $C_\Gamma$. Then, inspired by B. Kostant, we decompose the Coxeter element of the Kac-Moody algebra attached to $C_\Gamma$ as a product of reflections of a special form, thereby suggesting an algebraic form for the McKay correspondence in dimension $3$.\\
\end{small}

\section{\textsf{Introduction}}

\subsection{\textsf{Framework and results}}

\noindent Let $\Gamma$ be a finite subgroup of $\mathbf{SL}_3\mathbb{C}$. In this paper, we determine how the
finite dimensional irreducible representations of $\mathbf{SL}_3\mathbb{C}$ decompose under the action of the subgroup $\Gamma$.
These representations are indexed by $\mathbb{N}^2$. For $(m,\,n)\in \mathbb{N}^2$, let $V(m,\,n)$ denote the corresponding simple finite dimensional module. Let $\{\gamma_0,\dots,\,\gamma_l\}$ be the set of irreducible characters of $\Gamma$. We determine the numbers $m_i(m,\,n)$ --- the multiplicity of the character $\gamma_i$ in the representation $V(m,\,n)$. For that effect we introduce the  formal
power series: $$P_\Gamma(t,\,u)_i=\sum_{m=0}^\infty\sum_{n=0}^\infty m_i(m,\,n)t^mu^n.$$

\noindent We show that $m_i(t,\,u)$ is a rational function. We determine the rational functions which are obtained in that way for all the finite subgroups of $\mathbf{SL}_3\mathbb{C}$.\\

\noindent The proof uses an inversion of the recursion formula for the numbers $m_i(m,\,n)$.
The recursion formula is obtained through the decomposition of the tensor product of $V(m,\,n)$
with the natural representation of $\mathbf{SL}_3\mathbb{C}$. The key observation which leads
to this inversion is that a certain pair matrices are simultaneously diagonalizable. The eigenvalues of the matrices are values from the character table of the group $\Gamma$. This leads to the proof that the power series $$P_\Gamma(t,\,u)_i=\sum_{m=0}^\infty\sum_{n=0}^\infty m_i(m,\,n)t^mu^n$$ is rational. The actual calculation of this rational function then reduces to matrix multiplication.\\

\noindent This method applies indeed  to the $\mathbf{SL}_2\mathbb{C}$ case. It gives a complete (very short) proof of the results
obtained by B. Kostant in \cite{Kos85}, \cite{Kos06}, and by Gonzalez-Sprinberg and Verdier in \cite{GSV83}, and leads to an explicit determination of all the above multiplicities for the finite subgroups of $\mathbf{SL}_2\mathbb{C}$.\\
Although the results for $\mathbf{SL}_2\mathbb{C}$  are not new, the explicit relation of the rational functions with the eigenvalues of the Cartan matrix attached to the finite subgroup of $\mathbf{SL}_2\mathbb{C}$  doesn't seem to appear in the literature. In \cite{Kos85} this is established through the analysis of the orbit structure of the Coxeter element.\\

\noindent The construction of a minimal resolution of singularities of the orbifold $\mathbb{C}^3/\Gamma$ centralizes a lot of interest. It is related to the geometric McKay correspondence, cf. (for example)  \cite{BKR01}, \cite{GSV83}.
In this framework  Gonzalez-Spriberg-Verdier \cite{GSV83} use the Poincar\'e series determined above in their explicit construction of minimal resolution for singularities for $V=\mathbb{C}^2/\Gamma$ when $\Gamma$ is a finite subgroup of $\mathbf{SL}_2\mathbb{C}$.
  Following that approach the results of our calculation could be eventually used to construct explicit synthetic
  minimal resolution of singularities for orbifolds of the form $\mathbf{SL}_3\mathbb{C}/\Gamma$ where $\Gamma$ is a finite subgroup of $\mathbf{SL}_3\mathbb{C}$. This might clarify the description of the exceptional fiber of the minimal resolution
  of $\mathbf{SL}_3\mathbb{C}/\Gamma$ (see \cite{GNS04}).\\

\noindent An essential ingredient of the approach  of  B. Kostant  in \cite{Kos85} is the decomposition of a Coxeter element in the Weyl group attached to the Lie algebra corresponding to a subgroup $\Gamma$ of $\mathbf{SL}_2\mathbb{C}$ through the McKay correspondence as a product of simple reflections belonging to mutually orthogonal sets of roots.

\noindent Inspired by this approach, we attach to each finite subgroup $\Gamma$ of $\mathbf{SL}_3\mathbb{C}$ a generalized Cartan matrix $C_\Gamma$. We then factorize this matrix as a product of elements in the Weyl group of the Kac-Moody Lie algebra corresponding to $C_\Gamma$. These  elements are products of  simple reflections corresponding to roots in mutually orthogonal sets.\\

\subsection{\textsf{Organization of the paper}}

\noindent In \textbf{Section 2} we treat the  $\mathbf{SL}_2\mathbb{C}$ case. We show that the formal power series of the
multiplicities is a rational function by showing that it is an entry in a vector obtained as product of three matrices, two of which are scalar matrices the third one being a matrix with rational entries, by a scalar vector. We calculate the matrices for each finite subgroup of $\mathbf{SL}_3\mathbb{C}$. We give then the rational functions obtained.\\

\noindent In \textbf{Section 3} we apply the above method for the finite subgroups of $\mathbf{SL}_3\mathbb{C}$. Here we use the notations of [YY93] in which a classification of the finite subgroups of $\mathbf{SL}_3\mathbb{C}$ is presented.\\
Here again we prove the rationality of the formal power series of the multiplicities by showing that each such a series is an entry in the product of three matrices, two of them are scalar matrices
and the third being a matrix with rational entries, with a scalar vector.\\
For each finite subgroup of $\mathbf{SL}_3\mathbb{C}$ we give the the matrices involved in the product. To each subgroup $\Gamma$ we attach a generalized Cartan matrix $C_\Gamma$ (McKay correspondence in dimension 3)  we show its graph and its decomposition as a product of elements in the Weyl group of the Kac-Moody Lie algebra $\go{g}(C_\Gamma)$.\\

\noindent Then, for the series $A,\ B,\ C$ (\cite{YY93} notation) we give all the rational functions explicitly, As  for  the series $D$ we give the results for some specific examples because the description of the matrices engaged, in full generality doesn't have a simply presentable form.\\

\noindent For the exceptional finite subgroups of $\mathbf{SL}_3\mathbb{C}$ the numerators of the rational functions tend to be very long and we give them explicitly only for the cases where they are reasonably presentable. In all the cases we give the denominators explicitly. This is done in \textbf{Section 4}.\\


\section{\textsf{Branching law for the finite subgroups of $\mathbf{SL}_2\mathbb{C}$}}

\subsection{\textsf{The formal power series of the multiplicities is a rational function }}

\noindent $\bullet$ Let $\Gamma$ be a finite subgroup of
$\mathbf{SL}_2\mathbb{C}$ and
$\{\gamma_0,\dots,\,\gamma_l\}$ the set of
equivalence classes of irreducible finite dimensional complex
representations of $\Gamma$, where $\gamma_0$ is the trivial
representation. We denote by
 $\chi_j$ the character  associated to $\gamma_i$.\\
Consider $\gamma : \Gamma\rightarrow
\mathbf{SL}_2\mathbb{C}$, the natural $2-$dimensional
representation. Its character is denoted by $\chi$.
 We have then the  decomposition
$\gamma_j\otimes\gamma=\bigoplus_{i=0}^l a_{ij}\gamma_i$ for every $j\in\llbr 0,\,l\rrbr$. This defines an
$(l+1)\times (l+1)$ square matrix
$A:=\left(a_{ij}\right)_{(i,j)\in\llbr 0,\, l\rrbr^2}$ .\\

\noindent $\bullet$  Let $\go{h}$ be a Cartan
subalgebra of $\go{sl}_2\mathbb{C}$ and
$\varpi_1$ be the corresponding fundamental weight, and
 $V(n\varpi_1)$ be the simple $\go{sl}_2$-module of highest weight $n\varpi_1$.
 This give rise to an irreducible representation $\pi_n:\mathbf{SL}_2\mathbb{C}\longrightarrow V(n\varpi_1)$.\\
The restriction of
$\pi_{n}$ to the subgroup $\Gamma$, is a representation
of~$\Gamma$, and by complete reducibility, we have a
decomposition ${\pi_{n}}|_{\Gamma}=\bigoplus_{i=0}^l
m_i(n)\gamma_i,$ where the $m_i(n)$'s are non negative integers.
Let $\mathcal{E}:=(e_0,\dots,\,e_l)$ be the canonical basis of
$\mathbb{C}^{l+1}$, and $$\displaystyle{v_{n}:=\sum_{i=0}^l
m_i(n)e_i\in \mathbb{C}^{l+1}}.$$\\
As $\gamma_0$ is the trivial representation, we have
$v_{0}=e_0$. Let us consider the vector (with elements of
$\mathbb{C}\llbr t \rrbr$ as coefficients)
$$\displaystyle{P_\Gamma(t):=\sum_{n=0}^\infty
v_{n}t^n\in\left(\mathbb{C}\llbr t\rrbr\right)^{l+1}},$$
and denote by $\displaystyle{P_\Gamma(t)_j}$ its $j-$th coordinate in the basis $\mathcal{E}$. The series $P_\Gamma(t)_0$ is the Poincar\'e series of the invariant ring. Note also that $P_\Gamma(t)$ can also be seen as a formal power
series with coefficients in $\mathbb{C}^{l+1}$.\\
We proceed to calculate $P_\Gamma(t)$ . \\

\noindent $\bullet$ We get by Clebsch-Gordan formula that : $\pi_n\otimes\pi_1=\pi_{n+1}\oplus\pi_{n-1}$, so we have  $Av_n=v_{n+1}+v_{n-1}.$ From this  we deduce the relation $$(1-tA+t^2)P_\Gamma(t)=v_0.$$
Let us denote by $\{C_0,\dots,\,C_l\}$ the set of conjugacy
classes of $\Gamma$, and for any $j\in\llbr 0,\,l\rrbr$, let $g_j$
be an element of $C_j$. So the character table of $\Gamma$ is the
matrix $T_\Gamma\in \mathbf{M}_{l+1}\mathbb{C}$ defined by
$({T_\Gamma})_{i,j}:=\chi_i(g_j)$.\\
For all the finite subgroups of $\mathbf{SL}_2\mathbb{C}$ we have that,
$T_\Gamma$ is invertible, and $\Lambda:=T_\Gamma^{-1}\,A\,T_\Gamma$ is diagonal, with $\Lambda_{jj}=\overline{\chi(g_j)}$.\\
Set $\Theta:=(\Lambda_{00},\dots,\,\Lambda_{ll})$. We deduce from the preceding formula that $$T_\Gamma(1-t\Lambda+t^2)T_\Gamma^{-1}P_\Gamma(t)=v_0.$$
Let us define the rational function $$\begin{array}{rlcl}
  f\ : & \mathbb{C}^2 & \rightarrow & \mathbb{C}(t) \\
   & d & \mapsto & \displaystyle{\frac{1}{1-td+t^2}}. \\
\end{array}$$ Then
$$P_\Gamma(t)=T_\Gamma\,\Delta(t)\,T_\Gamma^{-1}v_0=(T_\Gamma\,\Delta(t)\,T_\Gamma)\,(T_\Gamma^{-2}v_0),$$ where $\Delta(t)\in
\mathbf{M}_{l+1}\mathbb{C}(t)$ is the diagonal matrix with
coefficients in $\mathbb{C}(t)$, defined by $\Delta_{jj}(t)=f(\Lambda_{jj}).$
Consequently, the coefficients of the vector
$P_\Gamma(t)$ are rational fractions in $t$.\\
Hence we get:

\begin{Pro}$\\$
For each $i\in\llbr0,\,l\rrbr$, the formal power series $P_\Gamma(t)_i$ is a rational function.\\
\end{Pro}

\subsection{\textsf{The results for the finite subgroups of $\mathbf{SL}_2\mathbb{C}$ }}

\noindent $\bullet$ The complete classification up to conjugation of all
finite subgroups of $\mathbf{SL}_2\mathbb{C}$ is given in \cite{Sp77}. It
consists of two infinite series (types $A,\,D$) and
three exceptional cases (types $E_6,\,E_7,\,E_8$).\\
We set $\displaystyle{\zeta_j:=e^{\frac{2i\pi}{j}}}$. For $\sigma\in\go{S}_{\llbr 0,\,j-1\rrbr}$, we then define the matrix $Q^\sigma:=\left(\delta_k^{\sigma(l)}\right)_{(k,l)\in\llbr 0,\,j-1\rrbr}.$\\

\subsubsection{\textsf{Type $A$ --- Cyclic groups}}

\noindent $\bullet$ Here, we take $\Gamma=\mathbb{Z}/j\mathbb{Z}$. The natural representation and the natural character of $\Gamma$ are $$\begin{array}{rcl|rcl}
  \gamma\ :\ \mathbb{Z}/j\mathbb{Z} & \rightarrow & \mathbf{SL}_2\mathbb{C} & \chi\ :\ \mathbb{Z}/j\mathbb{Z} & \rightarrow & \mathbf{SL}_2\mathbb{C}\\
  \overline{k} & \mapsto & \left(%
\begin{array}{cc}
  \zeta_j^k & 0 \\
   0 & \zeta_j^{-k} \\
\end{array}%
\right) & \overline{k} & \mapsto & \zeta_j^{k}+\zeta_j^{-k}.\\
\end{array}$$
The character table is the Vandermonde matrix
$T_\Gamma=\left(\zeta_j^{kl}\right)_{(k,l)\in\llbr
0,\,j-1\rrbr}$. Let
$\sigma$ be the permutation $\sigma\in\go{S}_{\llbr 0,\,j-1\rrbr}$ defined by
$\sigma(0)=0$ and $\forall\ i\in \llbr 0,\,j-1\rrbr,\
\sigma(i)=j-i$. Then $T_\Gamma^2=j\,Q^\sigma$, i.e.
$T_\Gamma^{-1}=\frac{1}{j}\,T_\Gamma Q^\sigma.$ The eigenvalues of $A$ are the
numbers $\overline{\chi(k)}=\zeta_j^k+\zeta_j^{-k},$ for
$k\in\llbr 0,\,j-1\rrbr$. Then
$$P_\Gamma(t)_i=\frac{1}{j}\,\left(T_\Gamma\Delta(t)T_\Gamma
Q^\sigma\right)_{i0}=\frac{1}{j}\,\sum_{p=0}^{j-1}\frac{\zeta_j^{ip}}{(1-t\zeta_j^p)(1-t\zeta_j^{-p})}.$$
Note that $(1-t^j)(1-t^2)$ is a common denominator of all the
terms of the preceding sum.\\

\subsubsection{\textsf{Type $D$ --- Binary dihedral groups}}

\noindent The binary dihedral group is the subgroup $\langle a_n,\,b\rangle$ of $\mathbf{SL}_2\mathbb{C}$, with
$$a_n:=\left(
       \begin{array}{cc}
         \zeta_{2n} & 0 \\
         0 & \zeta_{2n}^{-1} \\
       \end{array}
     \right)
,\
b:=\left(
     \begin{array}{cc}
       0 & i \\
       i & 0 \\
     \end{array}
   \right)
.$$
The order of $\Gamma$ is $4n$. The $n+3$ conjugacy classes of $\Gamma$ are
$$\begin{array}{||l||c|c|c|c|c|c|c|c|c||}\hline\hline & & & & & & & & & \\[-.4cm]
   \textrm{Class}  & id & a_nb & b & a_n^n & a_n & a_n^2 & a_n^3 & \cdots & a_n^{n-1} \\ \hline
  \textrm{Cardinality} & 1 & n & n & 1 & 2 & 2 & 2 & \cdots & 2 \\ \hline\hline
  \end{array}
$$
The character table of $\Gamma$ is
\begin{small}
$$T_\Gamma:=\left(
              \begin{array}{cccc|cccc}
                1 & 1 & 1 & 1 & 1 & 1 & \cdots & 1 \\
                1 & i^n & -i^n & (-1)^n & -1 & 1 & \cdots & (-1)^{n-1} \\
                1 & -i^n & i^n & (-1)^n & -1 & 1 & \cdots & (-1)^{n-1} \\
                1 & -1 & -1 & 1 & 1 & 1 & \cdots & 1 \\ \hline
                2 & 0 & 0 & -2 & \zeta_{2n}+\zeta_{2n}^{-1} & \zeta_{2n}^{2}+\zeta_{2n}^{-2} & \cdots & \zeta_{2n}^{n-1}+\zeta_{2n}^{-(n-1)} \\
                2 & 0 & 0 & 2 & \zeta_{2n}^{2}+\zeta_{2n}^{-2} & \zeta_{2n}^{4}+\zeta_{2n}^{-4} & \cdots & \zeta_{2n}^{2(n-1)}+\zeta_{2n}^{-2(n-1)} \\
                \vdots & \vdots & \vdots & \vdots & \vdots & \vdots &  & \vdots \\
                2 & 0 & 0 & (-1)^{n-1}2 & \zeta_{2n}^{n-1}+\zeta_{2n}^{-(n-1)} & \zeta_{2n}^{2(n-1)}+\zeta_{2n}^{-2(n-1)} & \cdots & \zeta_{2n}^{(n-1)(n-1)}+\zeta_{2n}^{-(n-1)(n-1)} \\
              \end{array}
            \right)
$$
\end{small}
The natural character $\chi$ of $\Gamma$ is given by $(\chi(g_0),\dots,\,\chi(g_l))=\overline{\Theta}$, with
$$\begin{array}{rcl}
    \Theta\ =\ \overline{\Theta} & = & [tr(id),\ tr(a_nb),\ tr(b),\ tr(a_n^n),\ tr(a_n),\ tr(a_n^2),\ tr(a_n^3),\dots,\ tr(a_n^{n-1})]  \\
   & = & [2,\ 0,\ 0,\ -2,\ \zeta_{2n}+\zeta_{2n}^{-1},\ \zeta_{2n}^2+\zeta_{2n}^{-2},\ \zeta_{2n}^3+\zeta_{2n}^{-3},\dots,\ \zeta_{2n}^{n-1}+\zeta_{2n}^{-(n-1)}].
  \end{array}
$$
Set $Diag\left(d_1,\ d_2,\ d_3,\ d_4,\ \delta_1,\ \delta_2,\dots,\ \delta_{n-1}\right):=\Delta(t)$. We deduce the formula for the series $P_\Gamma(t)$ that :
$$\begin{array}{rcl}
    P_\Gamma(t)_0 & = & \displaystyle{\frac{3n-1}{8n^2}\left(d_1+d_2+d_3+d_4+2\sum_{k=1}^{n-1}\delta_k\right)
    +(-1)^n\,\frac{n-1}{8n^2}\left(d_1+(-1)^n(d_2+d_3)+d_4+2\sum_{k=1}^{n-1}(-1)^k\delta_k\right)} \\
     &  & \displaystyle{+\sum_{l=1}^{n-1}(-1)^l\,\frac{n-1}{4n^2}\left(d_1+(-1)^l(d_2+d_3)+d_4+\sum_{k=1}^{n-1}(\zeta_{2n}^{lk}+\zeta_{2n}^{-lk})\delta_k\right)},
  \end{array}
$$

$$\begin{array}{rcl}
    P_\Gamma(t)_1 & = & \displaystyle{\frac{3n-1}{8n^2}\left(d_1+i^nd_2-i^nd_3+(-1)^nd_4+2\sum_{k=1}^{n-1}(-1)^k\delta_k\right)}
     \\
     & & \displaystyle{+(-1)^n\,\frac{n-1}{8n^2}\left(d_1+(-1)^n(i^nd_2-i^nd_3+d_4)+2\sum_{k=1}^{n-1}\delta_k\right)} \\
     &  & \displaystyle{+\sum_{l=1}^{n-1}(-1)^l\,\frac{n-1}{4n^2}\left(d_1+(-1)^l(i^nd_2-i^nd_3)+(-1)^nd_4+\sum_{k=1}^{n-1}(-1)^k(\zeta_{2n}^{lk}+\zeta_{2n}^{-lk})\delta_k\right)},
  \end{array}
$$
and then $P_\Gamma(t)_{3}$ (resp. $P_\Gamma(t)_{2}$) is obtained by replacing in $P_\Gamma(t)_{0}$ (resp. $P_\Gamma(t)_{1}$) $d_2$ by $-d_2$ and $d_3$ by $-d_3$.\\
Finally, for $i\in\llbr 1,\,n-1\rrbr$, we have

$$\begin{array}{rcl}
    P_\Gamma(t)_{i+3} & = & \displaystyle{\frac{3n-1}{8n^2}\left(2d_1+2(-1)^id_4+2\sum_{k=1}^{n-1}(\zeta_{2n}^{ki}+\zeta_{2n}^{-ki})\delta_k\right)}
     \\
     & & \displaystyle{+(-1)^n\,\frac{n-1}{8n^2}\left(2d_1+2(-1)^id_4+2\sum_{k=1}^{n-1}(-1)^k(\zeta_{2n}^{ik}+\zeta_{2n}^{-ik})\delta_k\right)} \\
     &  & \displaystyle{+\sum_{l=1}^{n-1}(-1)^l\,\frac{n-1}{4n^2}\left(2d_1+2(-1)^{i}d_4+\sum_{k=1}^{n-1}(\zeta_{2n}^{ki}+\zeta_{2n}^{-ki})\delta_k(\zeta_{2n}^{kl}+\zeta_{2n}^{-kl})\right)}.
  \end{array}
$$

\subsection{\textsf{Exceptional cases}}

\subsubsection{\textsf{Type $E_6$ --- Binary tetrahedral group}}

\noindent The binary tetrahedral group is the subgroup $\langle a^2,\,b,\,c\rangle$ of $\mathbf{SL}_2\mathbb{C}$, with
$$a:=\left(
       \begin{array}{cc}
         \zeta_8 & 0 \\
         0 & \zeta_8^{7} \\
       \end{array}
     \right)
,\
b:=\left(
     \begin{array}{cc}
       0 & i \\
       i & 0 \\
     \end{array}
   \right)
,\
c:=\frac{1}{\sqrt{2}}\left(
     \begin{array}{cc}
       \zeta_8^7 & \zeta_8^7 \\
       \zeta_8^5 & \zeta_8 \\
     \end{array}
   \right)
.$$
The order of $\Gamma$ is $24$. The $7$ conjugacy classes of $\Gamma$ are
$$\begin{array}{||l||c|c|c|c|c|c|c||}\hline\hline & & & & & & & \\[-.4cm]
   \textrm{Class}  & id & a^4=-id & b & c & c^2 & -c & -c^2 \\ \hline
  \textrm{Cardinality} & 1 & 1 & 6 & 4 & 4 & 4 & 4 \\ \hline\hline
  \end{array}
$$

\noindent The character table $T_\Gamma$ of $\Gamma$ and the matrix $A$ are
\begin{small}
$$T_\Gamma=\left( \begin {array}{ccccccc} 1&1&1&1&1&1&1\\\noalign{\medskip}1&1&1&j&{j}^{2}&j&{j}^{2}\\\noalign{\medskip}1&1&1&{j}^{2}&j&{j}^{2}&j
\\\noalign{\medskip}2&-2&0&1&-1&-1&1\\\noalign{\medskip}2&-2&0&j&-{j}^
{2}&-j&{j}^{2}\\\noalign{\medskip}2&-2&0&{j}^{2}&-j&-{j}^{2}&j
\\\noalign{\medskip}3&3&-1&0&0&0&0\end {array} \right),\ A=\left( \begin {array}{ccccccc} 0&0&0&1&0&0&0\\\noalign{\medskip}0&0&0&0&1&0&0\\\noalign{\medskip}0&0&0&0&0&1&0\\\noalign{\medskip}1&0&0&0&0
&0&1\\\noalign{\medskip}0&1&0&0&0&0&1\\\noalign{\medskip}0&0&1&0&0&0&1
\\\noalign{\medskip}0&0&0&1&1&1&0\end {array} \right),$$
\end{small}
and the eigenvalues are $\Theta=(2,\ -2,\ 0,\ 1,\ -1,\ -1,\ 1).$\\

\noindent The series $P_\Gamma(t)_i=\frac{N_\Gamma(t)_i}{D_\Gamma(t)}$ are given by $D_\Gamma(t)=(1-t^6)(1-t^8),$ and
$$\begin{array}{rcl|rcl}
    N_\Gamma(t)_0 & = & t^{12}+1, & N_\Gamma(t)_4 & = & t^9+t^7+t^5+t^3, \\
    N_\Gamma(t)_1 & = & t^8+t^4, & N_\Gamma(t)_5 & = & t^9+t^7+t^5+t^3, \\
    N_\Gamma(t)_2 & = & t^8+t^4, & N_\Gamma(t)_6 & = & t^{10}+t^8+2t^6+t^4+t^2. \\
    N_\Gamma(t)_3 & = & t^{11}+t^7+t^5+t, &  &  &
  \end{array}
$$

\subsubsection{\textsf{Type $E_7$ --- Binary octahedral group}}

\noindent The binary octahedral group is the subgroup $\langle a,\,b,\,c\rangle$ of $\mathbf{SL}_2\mathbb{C}$, with $a,\ b,\ c$
defined as in the preceding section. The order of $\Gamma$ is $48$. The $8$ conjugacy classes of $\Gamma$ are
$$\begin{array}{||l||c|c|c|c|c|c|c|c||}\hline\hline & & & & & & & & \\[-.4cm]
   \textrm{Class}  & id & a^4=-id & ab & b & c^2 & c & a & a^3 \\ \hline
  \textrm{Cardinality} & 1 & 1 & 12 & 6 & 8 & 8 & 6 & 6 \\ \hline\hline
  \end{array}
$$

\noindent The character table $T_\Gamma$ of $\Gamma$ and the matrix $A$ are
\begin{small}
$$T_\Gamma=\left( \begin {array}{cccccccc} 1&1&1&1&1&1&1&1\\\noalign{\medskip}1&1&-1&1&1&1&-1&-1\\\noalign{\medskip}2&2&0&2&-1&-1&0&0
\\\noalign{\medskip}2&-2&0&0&-1&1&\sqrt{2} &-\sqrt{2} \\\noalign{\medskip}2&-2&0&0&-1&1&-\sqrt{2} &\sqrt{2}
\\\noalign{\medskip}3&3&-1&-1&0&0&
1&1\\\noalign{\medskip}3&3&1&-1&0&0&-1&-1\\\noalign{\medskip}4&-4&0&0&
1&-1&0&0\end {array} \right),\ A=\left( \begin {array}{cccccccc} 0&0&0&1&0&0&0&0\\\noalign{\medskip}0&0&0&0&1&0&0&0\\\noalign{\medskip}0&0&0&0&0&0&0&1\\\noalign{\medskip}1&0
&0&0&0&1&0&0\\\noalign{\medskip}0&1&0&0&0&0&1&0\\\noalign{\medskip}0&0
&0&1&0&0&0&1\\\noalign{\medskip}0&0&0&0&1&0&0&1\\\noalign{\medskip}0&0
&1&0&0&1&1&0\end {array} \right),$$
\end{small}
and the eigenvalues are $\Theta=(2,\ -2,\ 0,\ 0,\ -1,\ 1,\ \sqrt{2},\ -\sqrt{2}).$\\

\noindent The series $P_\Gamma(t)_i=\frac{N_\Gamma(t)_i}{D_\Gamma(t)}$ are given by $D_\Gamma(t)=(1-t^8)(1-t^{12}),$ and
$$\begin{array}{rcl|rcl}
    N_\Gamma(t)_0 & = & t^{18}+1, & N_\Gamma(t)_4 & = & t^{13}+t^{11}+t^7+t^5, \\
    N_\Gamma(t)_1 & = & t^{12}+t^6, & N_\Gamma(t)_5 & = & t^{16}+t^{12}+t^{10}+t^8+t^6+t^2, \\
    N_\Gamma(t)_2 & = & t^{14}+t^{10}+t^8+t^4, & N_\Gamma(t)_6 & = & t^{14}+t^{12}+t^{10}+t^8+t^6+t^4, \\
    N_\Gamma(t)_3 & = & t^{17}+t^{11}+t^7+t, & N_\Gamma(t)_7 & = & t^{15}+t^{13}+t^{11}+2t^9+t^7+t^5+t^3.
  \end{array}
$$

\subsubsection{\textsf{Type $E_8$ --- Binary icosahedral group}}

\noindent The binary icosahedral group is the subgroup $\langle a,\,b,\,c\rangle$ of $\mathbf{SL}_2\mathbb{C}$, with
$$a:=\left(
     \begin{array}{cc}
       -\zeta_5^3 & 0 \\
       0 & -\zeta_5^2 \\
     \end{array}
   \right)
,\
b:=\left(
       \begin{array}{cc}
         0 & 1 \\
         -1 & 0 \\
       \end{array}
     \right)
,\
c:=\frac{1}{\zeta_5^2+\zeta_5^{-2}}\left(
     \begin{array}{cc}
       \zeta_5+\zeta_5^{-1} & 1 \\
       1 & -\zeta_5-\zeta_5^{-1} \\
     \end{array}
   \right)
.$$
The order of $\Gamma$ is $120$. The $9$ conjugacy classes of $\Gamma$ are
$$\begin{array}{||l||c|c|c|c|c|c|c|c|c||}\hline\hline & & & & & & & & & \\[-.4cm]
   \textrm{Class}  & id & b^2=-id & a & a^2 & a^3 & a^4 & abc & (abc)^2 & b \\ \hline
  \textrm{Cardinality} & 1 & 1 & 12 & 12 & 12 & 12 & 20 & 20 & 30 \\ \hline\hline
  \end{array}
$$

\noindent The character table $T_\Gamma$ of $\Gamma$ and the matrix $A$ are
\begin{small}
$$T_\Gamma:=\left( \begin {array}{ccccccccc} 1&1&1&1&1&1&1&1&1\\\noalign{\medskip}2&-2&\frac{1-\sqrt {5}}{2}&-\frac{1+\sqrt {5}}{2}&\frac{1+\sqrt {5}}{2}&\frac{-1+\sqrt {5}}{2}&1&-1&0\\\noalign{\medskip}2&-2&\frac{1+\sqrt {5}}{2}&\frac{-1+\sqrt {5}}{2}&\frac{1-\sqrt {5}}{2}&-\frac{1+\sqrt {5}}{2}&
1&-1&0\\\noalign{\medskip}3&3&\frac{1+\sqrt {5}}{2}&\frac{1-\sqrt {5}}{2}&\frac{1-\sqrt {5}}{2}&\frac{1+\sqrt {5}}{2}&0&0&-1\\\noalign{\medskip}3&3&\frac{1-\sqrt {5}}{2}&\frac{1+\sqrt {5}}{2}&\frac{1+\sqrt {5}}{2}&\frac{1-\sqrt {5
}}{2}&0&0&-1\\\noalign{\medskip}4&4&-1&-1&-1&-1&1&1&0\\\noalign{\medskip}4
&-4&1&-1&1&-1&-1&1&0\\\noalign{\medskip}5&5&0&0&0&0&-1&-1&1
\\\noalign{\medskip}6&-6&-1&1&-1&1&0&0&0\end {array} \right),\ A=\left( \begin {array}{ccccccccc} 0&0&1&0&0&0&0&0&0\\\noalign{\medskip}0&0&0&0&0&1&0&0&0\\\noalign{\medskip}1&0&0&1&0&0&0
&0&0\\\noalign{\medskip}0&0&1&0&0&0&1&0&0\\\noalign{\medskip}0&0&0&0&0
&0&0&0&1\\\noalign{\medskip}0&1&0&0&0&0&0&0&1\\\noalign{\medskip}0&0&0
&1&0&0&0&1&0\\\noalign{\medskip}0&0&0&0&0&0&1&0&1\\\noalign{\medskip}0
&0&0&0&1&1&0&1&0\end {array} \right)
,$$
\end{small}
and the eigenvalues are
$\Theta=\left(2,\ -2,\ \frac{1+\sqrt{5}}{2},\ \frac{-1+\sqrt{5}}{2},\ \frac{1-\sqrt{5}}{2},\ \frac{-1-\sqrt{5}}{2},\ 1,\ -1,
\ 0\right).$\\

\noindent The series $P_\Gamma(t)_i=\frac{N_\Gamma(t)_i}{D_\Gamma(t)}$ are given by $D_\Gamma(t)=(1-t^{12})(1-t^{20}),$ and
$$\begin{array}{rcl|rcl}
    N_\Gamma(t)_0 & = & t^{30}+1, & N_\Gamma(t)_4 & = & t^{24}+t^{20}+t^{16}+t^{14}+t^{10}+t^6, \\
    N_\Gamma(t)_1 & = & t^{23}+t^{17}+t^{13}+t^7, & N_\Gamma(t)_5 & = & t^{24}+t^{22}+t^{18}+t^{16}+t^{14}+t^{12}+t^8+t^6, \\
    N_\Gamma(t)_2 & = & t^{29}+t^{19}+t^{11}+t, & N_\Gamma(t)_6 & = & t^{27}+t^{21}+T^{19}+t^{17}+t^{13}+t^{11}+t^9+t^3, \\
    N_\Gamma(t)_3 & = & t^{28}+t^{20}+t^{18}+t^{12}+t^{10}+t^2, &  &  &
  \end{array}$$
$$\begin{array}{rcl}
    N_\Gamma(t)_7 & = & t^{26}+t^{22}+t^{20}+t^{18}+t^{16}+t^{14}+t^{12}+t^{10}+t^8+t^4, \\
    N_\Gamma(t)_8 & = & t^{25}+t^{23}+t^{21}+t^{19}+t^{17}+2t^{15}+t^{13}+t^{11}+t^9+t^7+t^5.
  \end{array}
$$


\section{\textsf{Branching law for the finite subgroups of $\mathbf{SL}_3\mathbb{C}$}}

\noindent $\bullet$ Let $\Gamma$ be a finite subgroup of
$\mathbf{SL}_3\mathbb{C}$ and
$\{\gamma_0,\dots,\,\gamma_l\}$ the set of
equivalence classes of irreducible finite dimensional complex
representations of $\Gamma$, where $\gamma_0$ is the trivial
representation. The character
associated to $\gamma_j$ is denoted by $\chi_j$.\\
Consider $\gamma : \Gamma\rightarrow
\mathbf{SL}_3\mathbb{C}$ the natural $3-$dimensional
representation, and $\gamma^*$ its contragredient representation. The character of $\gamma$ is denoted by $\chi$.
By complete reducibility  we get the decompositions
$$\forall\ j\in\llbr 0,\,l\rrbr,\,
\ \gamma_j\otimes\gamma=\bigoplus_{i=0}^l a^{(1)}_{ij}\gamma_i\ \
\textrm{and}\ \ \gamma_j\otimes\gamma^*=\bigoplus_{i=0}^l
a^{(2)}_{ij}\gamma_i.$$ This defines two square matrices
$A^{(1)}:=\left(a^{(1)}_{ij}\right)_{(i,j)\in\llbr 0,\, l\rrbr^2}$
and $A^{(2)}:=\left(a^{(2)}_{ij}\right)_{(i,j)\in\llbr 0,\,
l\rrbr^2}$ of $\mathbf{M}_{l+1}\mathbb{N}$.\\

\noindent $\bullet$ Let $\underline{\mathbf{h}}$ be a Cartan subalgebra of $\mathbf{sl_3}\mathbb{C}$ and let $\varpi_1,\varpi_2$ be the corresponding fundamental weights, and $V(m\varpi_1+n\varpi_2)$ the simple
$\mathbf{sl_3}\mathbb{C}$ module of highest weight $m\varpi_1+n\varpi_2$ with $(m,\,n)\in \mathbb{N}^2$. Then we get an irreducible representation
$\pi_{m,n}:\mathbf{SL}_3\mathbb{C}\rightarrow
\mathbf{GL}(V(m\varpi_1+n\varpi_2))$.
 The restriction of
$\pi_{m,n}$ to the subgroup $\Gamma$ is a representation
of~$\Gamma$, and by complete reducibility, we get the
decomposition $${\pi_{m,n}}|_{\Gamma}=\bigoplus_{i=0}^l
m_i(m,n)\gamma_i,$$ where the $m_i(m,n)$'s are non negative integers.
Let $\mathcal{E}:=(e_0,\dots,\,e_l)$ be the canonical basis of
$\mathbb{C}^{l+1}$, and $$\displaystyle{v_{m,n}:=\sum_{i=0}^l
m_i(m,n)e_i\in \mathbb{C}^{l+1}}.$$\\
As $\gamma_0$ is the trivial representation, we have
$v_{0,0}=e_0$. Let us consider the vector (with elements of
$\mathbb{C}\llbr t,\,u\rrbr$ as coefficients)
$$\displaystyle{P_\Gamma(t,\,u):=\sum_{m=0}^\infty\sum_{n=0}^\infty
v_{m,n}t^mu^n\in\left(\mathbb{C}\llbr t,\,u\rrbr\right)^{l+1}},$$
and denote by $\displaystyle{P_\Gamma(t,\,u)_j}$ its $j-$th coordinate in the basis $\mathcal{E}$.
Note that $P_\Gamma(t,\,u)$ can also be seen as a formal power
series with coefficients in $\mathbb{C}^{l+1}$. The aim of this article is to compute $P_\Gamma(t,\,u)$.\\

\subsection{\textsf{The formal power series of the multiplicies is a rational function }}

\noindent Here we establish some properties of
the series $P_\Gamma(t,\,u)$, in order to give an explicit formula
for it. The first proposition follows from  the uniqueness of the decomposition
of a representation as sum of irreducible representations.\\

\begin{Pro}$\\$
$\bullet$ $A^{(2)}=\,^tA^{(1)}$.\\
$\bullet$ $A^{(1)}$ and $A^{(2)}$ commute, i.e. $A^{(1)}$ is a
normal matrix.
\end{Pro}

\noindent Since $A^{(1)}$ is normal, we know that it is
diagonalizable with eigenvectors forming an orthogonal basis. Now we will
diagonalize the matrix $A^{(1)}$ by using the character table of
the group $\Gamma$.
Let us denote by $\{C_0,\dots,\,C_l\}$ the set of conjugacy
classes of $\Gamma$, and for any $j\in\llbr 0,\,l\rrbr$, let $g_j$
be an element of $C_j$. So the character table of $\Gamma$ is the
matrix $T_\Gamma\in \mathbf{M}_{l+1}\mathbb{C}$ defined by~$(T_\Gamma)_{i,j}:=\chi_i(g_j)$.\\

\begin{Pro}\label{vp}$\\$
For $k\in\llbr 0,\,l\rrbr$, set
$w_k:=(\chi_0(g_k),\dots,\,\chi_l(g_k))\in\mathbb{C}^{l+1}$. Then
$w_k$ is an eigenvector of $A^{(2)}$ associated to the eigenvalue
$\chi(g_k)$. Similarly, $w_k$ is an eigenvector of $A^{(1)}$ associated to the
eigenvalue $\overline{\chi(g_k)}$.\\
\end{Pro}

\noindent We will see in the sequel that
$\mathcal{W}:=(w_0,\dots,\,w_l)$ is always a basis of eigenvectors
of $A^{(1)}$ and $A^{(2)}$, so that $T_\Gamma^{-1}A^{(1)}T_\Gamma$
and $T_\Gamma^{-1}A^{(2)}T_\Gamma$
are diagonal matrices.\\
Now, we make use of the Clebsch-Gordan formula
\begin{equation}\label{CG}
    \pi_{1,0}\otimes\pi_{m,n}=\pi_{m+1,n}\oplus\pi_{m,n-1}\oplus\pi_{m-1,n+1},\ \
  \pi_{0,1}\otimes\pi_{m,n}=\pi_{m,n+1}\oplus\pi_{m-1,n}\oplus\pi_{m+1,n-1}.
\end{equation}

\begin{Pro}$\\$
The vectors $v_{m,n}$ satisfy the following recurrence relations
$$\begin{array}{c}
  A^{(1)}v_{m,n}=v_{m+1,n}+v_{m,n-1}+v_{m-1,n+1}, \\
  A^{(2)}v_{m,n}=v_{m,n+1}+v_{m-1,n}+v_{m+1,n-1}. \\
\end{array}$$
\end{Pro}

\underline{Proof:}\\
The definition of $v_{m,n}$ reads
$v_{m,n}=\sum_{i=0}^l m_i(m,n)e_i$, thus
$A^{(1)}v_{m,n}=\sum_{i=0}^l\left(\sum_{j=0}^l
m_j(m,n)a^{(1)}_{ij}\right)e_i$.\\
Now
$(\pi_{m,n}\otimes\pi_{1,0})|_\Gamma=\pi_{m,n}|_\Gamma\otimes\gamma=\sum_{j=0}^l
m_j(m,n)\gamma_j\otimes\gamma=\sum_{i=0}^l \left(\sum_{j=0}^l
m_j(m,n)a^{(1)}_{ij}\right)\gamma_i$,\\and
$\pi_{m+1,n}|_\Gamma+\pi_{m,n-1}|_\Gamma+\pi_{m-1,n+1}|_\Gamma = \sum_{i=0}^l
\left(m_i(m+1,n)+m_i(m,n-1)+m_i(m-1,n+1)\right)\gamma_i$.\\
By uniqueness, $\sum_{j=0}^l
m_j(m,n)a^{(1)}_{ij}=m_i(m+1,n)+m_i(m,n-1)+m_i(m-1,n+1).$ $\blacksquare$\\

\begin{Pro}\label{form}$\\$
The series $P_\Gamma(t,\,u)$ satisfies the following relation
$$\left(1-tA^{(1)}+t^2A^{(2)}-t^3\right)\left(1-uA^{(2)}+u^2A^{(1)}-u^3\right)P_\Gamma(t,\,u)=(1-tu)v_{0,0}.$$
\end{Pro}

\underline{Proof:}\\
$\bullet$ Set $x:=P_\Gamma(t,\,u)$. Set also $v_{m,-1}:=0$ and
$v_{-1,n}:=0$ for $(m,\,n)\in\mathbb{N}$, such that, according to
the Clebsch-Gordan formula, the formulae of the preceding
corollary are still true for $(m,\,n)\in\mathbb{N}$. We have
$$tuA^{(1)}x=tu\sum_{m=0}^\infty\sum_{n=0}^\infty A^{(1)}v_{m,n}t^m
u^n=\sum_{m=0}^\infty\sum_{n=0}^\infty
(v_{m+1,n}+v_{m,n-1}+v_{m-1,n+1})t^{m+1} u^{n+1}.$$
$$\textrm{Now}\ \ \ \ \ \ \ \ \ \ \ \ \sum_{m=0}^\infty\sum_{n=0}^\infty v_{m+1,n}t^{m+1}
u^{n+1}=u\sum_{m=1}^\infty\sum_{n=0}^\infty v_{m,n}t^m
u^n=ux-u\sum_{n=0}^\infty v_{0,n}u^n,$$
$$\sum_{m=0}^\infty\sum_{n=0}^\infty v_{m,n-1}t^{m+1}
u^{n+1}=tu^2\sum_{m=0}^\infty\sum_{n=1}^\infty v_{m,n-1}t^m
u^{n-1}=tu^2x,$$ $$\begin{array}{rcl}
  \textrm{and}\ \ \ \ \displaystyle{\sum_{m=0}^\infty\sum_{n=0}^\infty v_{m-1,n+1}t^{m+1}
u^{n+1}} & = & \displaystyle{t^2\sum_{m=1}^\infty\sum_{n=0}^\infty
v_{m-1,n+1}t^{m-1} u^{n+1}=t^2\sum_{m=0}^\infty\sum_{n=0}^\infty
v_{m,n+1}t^{m}u^{n+1}} \\
 &   = & \displaystyle{t^2\sum_{m=0}^\infty\sum_{n=0}^\infty v_{m,n}t^{m}
u^{n}-t^2\sum_{m=0}^\infty v_{m,0}t^{m}=t^2x-t^2\sum_{m=0}^\infty
v_{m,0}t^{m}.} \\
\end{array}$$
\begin{equation}\label{form1}
\textrm{Therefore}\ \ \ \ \ \ \ \ \ \ \ \ tuA^{(1)}x=(u+tu^2+t^2)x-u\sum_{n=0}^\infty
v_{0,n}u^n-t^2\sum_{m=0}^\infty v_{m,0}t^m.
\end{equation}
\begin{equation}\label{form2}
\textrm{We\ proceed\ likewise\ to\ obtain}\ \ \ \ tuA^{(2)}x=(t+tu^2+u^2)x-t\sum_{m=0}^\infty
v_{m,0}t^m-u^2\sum_{n=0}^\infty v_{0,n}u^n.
\end{equation}
$\bullet$ By using Equations (\ref{form1}) and (\ref{form2}), we
have $\displaystyle{tuA^{(2)}x-tu^2A^{(1)}x=t(1-u^3)x+(t^2u-t)\sum_{m=0}^\infty
v_{m,0}t^m,}$
\begin{equation}\label{form3}
\textrm{i.e.}\ \ \ \ \left(1-uA^{(2)}+u^2A^{(1)}-u^3\right)x=(1-tu)\sum_{m=0}^\infty
v_{m,0}t^m.
\end{equation}
Besides $A^{(1)}v_{m,0}=v_{m+1,0}+v_{m-1,1}$, and
$A^{(2)}v_{m-1,0}=v_{m-1,1}+v_{m-2,0}$, hence
$$A^{(1)}v_{m,0}=v_{m+1,0}+A^{(2)}v_{m-1,0}-v_{m-2,0}.$$
Set $\displaystyle{y:=\sum_{m=0}^\infty v_{m,0}t^m}$. Then
$$\begin{array}{rcl}
  tA^{(1)}y & = & \displaystyle{\sum_{m=0}^\infty
v_{m+1,0}t^{m+1}+A^{(2)}\sum_{m=1}^\infty
v_{m-1,0}t^{m+1}-\sum_{m=2}^\infty
v_{m-2,0}t^{m+1}} \\
  &  = & \displaystyle{\sum_{m=1}^\infty
v_{m,0}t^{m}+t^2A^{(2)}\sum_{m=0}^\infty
v_{m,0}t^{m}-t^3\sum_{m=0}^\infty
v_{m,0}t^{m}}\ =\ y-v_{0,0}+t^2A^{(2)}y-t^3y. \\
\end{array}$$
\begin{equation}\label{form4}
 \textrm{So}\ \ \ \ \left(1-tA^{(1)}+t^2A^{(2)}-t^3\right)y=v_{0,0}.
\end{equation}
Combining Eq. \ref{form3} and \ref{form4}, we have
$\left(1-tA^{(1)}+t^2A^{(2)}-t^3\right)\left(1-uA^{(2)}+u^2A^{(1)}-u^3\right)x=(1-tu)v_{0,0}.\ \blacksquare$\\

\noindent We may inverse the relation obtained in Proposition \ref{form}
and obtain an explicit expression\footnote{$P_\Gamma(t,\,u)=(1-tu)\left(\sum_{p=0}^\infty\left(u^3+uA^{(2)}-u^2A^{(1)}\right)^p\right)
\left(\sum_{q=0}^\infty\left(t^3+tA^{(1)}-t^2A^{(2)}\right)^q\right)v_{0,0}.$\\For $z\in\mathbb{R}$, let $\lceil z\rceil$ be the
smallest integer that is greater or equal to $z$, and set $\{r,\,s\}:=\{1,\,2\}$. For $m\in\mathbb{N}$, set
$$\displaystyle{\alpha_m^{(r)}:=\sum_{q=\lceil\frac{m}{3}\rceil}^m\left(\sum_{j=\lceil\frac{3q-m}{2}\rceil}^{\min(3q-m,\,q)}C_q^jC_j^{3q-m-j}
  (-1)^{3q-m}{A^{(r)}}^{3q-m-j}{A^{(s)}}^{2j-3q+m}\right)}$$
Then $v_{m,n}=v_{0,0}$ if $m=n=0$; $\alpha_n^{(2)}v_{0,0}$ if $m=0,\ n\neq0$; $\alpha_m^{(1)}v_{0,0}$ if $n=0,\ m\neq0$; $(\alpha_n^{(2)}\alpha_m^{(1)}-\alpha_{n-1}^{(2)}\alpha_{m-1}^{(1)})v_{0,0}$ otherwise.} for $P_\Gamma(t,\,u)$ as well as
an explicit formula for the vector $v_{m,n}$. But, for the explicit
calculations of $P_\Gamma(t,\,u)$, we will use an other fundamental formula (we inverse complex numbers instead of matrices). We need the rational function $f$ defined by $$\begin{array}{rlcl}
  f\ : & \mathbb{C}^2 & \rightarrow & \mathbb{C}(t,\,u) \\
   & (d_1,\,d_2) & \mapsto & \displaystyle{\frac{1-tu}{(1-td_1+t^2d_2-t^3)(1-ud_2+u^2d_1-u^3)}}. \\
\end{array}$$
The complete classification up to conjugation of all
finite subgroups of $\mathbf{SL}_3\mathbb{C}$ is given in \cite{YY93}. It
consists in four infinite series (types $A,\, B,\, C,\, D$) and eight exceptional cases (types $E,\ F,\ G,\ H,\ I,\ J,\ K,\ L$).\\
In all the cases, the character table $T_\Gamma$ is invertible, and $\Lambda^{(1)}:=T_\Gamma^{-1}\,A^{(1)}\,T_\Gamma$ and $\Lambda^{(2)}:=T_\Gamma^{-1}\,A^{(2)}\,T_\Gamma$ are diagonal matrices, with $\Lambda^{(1)}_{jj}=\overline{\chi(g_j)}$ and $\Lambda^{(2)}_{jj}=\chi(g_j)$. According to Proposition \ref{form}, we may write
 $$T_\Gamma\left(1-t\Lambda^{(1)}+t^2\Lambda^{(2)}-t^3\right)\left(1-u\Lambda^{(2)}+u^2\Lambda^{(1)}-u^3\right)T_\Gamma^{-1}P_\Gamma(t,\,u)=(1-tu)v_{0,0}.$$
We deduce that
\begin{equation}\label{fond}
P_\Gamma(t,\,u)=T_\Gamma\,\Delta(t,\,u)\,T_\Gamma^{-1}v_{0,0}=(T_\Gamma\,\Delta(t,\,u)\,T_\Gamma)\,(T_\Gamma^{-2}v_{0,0}),
\end{equation}
where $\Delta(t,\,u)\in\mathbf{M}_{l+1}\mathbb{C}(t,\,u)$ is the diagonal matrix defined by $\Delta(t,\,u)_{jj}=f(\Lambda_{jj},\,\overline{\Lambda_{jj}})=f(\overline{\chi(g_j)},\,\chi(g_j)).$
Let $\Theta:=(\Lambda^{(1)}_{00},\dots,\,\Lambda^{(1)}_{ll})$ be the list of eigenvalues of $A^{(1)}$.\\

\noindent As a corollary  of the preceding formula we get:\\

\begin{Pro}$\\$
The coefficients of the vector $P_\Gamma(t,\,u)$ are rational fractions in $t$ and $u$, hence the formal power series of the multiplicities is a rational function.\\
\end{Pro}

\noindent We will denote them by
$$P_\Gamma(t,\,u)_i:=\frac{N_\Gamma(t,\,u)_i}{D_\Gamma(t,\,u)_i},\ i\in\llbr 0,\,l\rrbr$$
where $N_\Gamma(t,\,u)_i$ and $D_\Gamma(t,\,u)_i$ are elements of
$\mathbb{C}[t,\,u]$ that will be explicitly computed in the sequel.\\
Finally, we introduce a generalized Cartan matrix that we will study for every finite subgroup of $\mathbf{SL}_3\mathbb{C}$.\\

\begin{Def}$\\$
For every finite subgroup of $\mathbf{SL}_3\mathbb{C}$, we define a generalized Cartan matrix by the following formula:
$$C_\Gamma:=2\,I-A^{(1)}-\,^tA^{(1)}+2\,Diag(A^{(1)}).$$
For $k\in\llbr 0,\,l\rrbr$, the matrix of the reflection $s_k$ associated to the $k-$th root of $g(C_\Gamma)$ the
Kac Moody algebra attached to $C_\Gamma$ is defined by
$$(s_k)_{ij}=\delta_i^j-(C_\Gamma)_{k,j}\delta_i^k.$$
\end{Def}

\noindent For each finite subgroup, we will give a decomposition of the set of simple reflections $\mathcal{S}=\{s_0,\dots,\,s_l\}$ in $p$ sets (with $p$ minimal), i.e. $$\mathcal{S}=S_0\sqcup\dots\sqcup S_{p-1},$$ such that roots corresponding to reflections in those sets form a partition of the set of simple roots to mutually orthogonal sets. We denote by $\tau_l$ the (commutative) product of the elements of $S_l$. Then we deduce the following decomposition of $C_\Gamma$:
$$C_\Gamma=p\,I-\sum_{k=0}^{p-1}\tau_k.$$

\begin{Rk}$\\$
Along this section we will present matrices that have
only $-2,\,-1,\,0,\,1,\,2$ as entries. For a clearer exposition, we represent the non-zero entries by
colored points. The correspondence is the following: dark grey
$=-2$, light grey $=-1$, white $=1$, black $=2$, empty $=0$.\\
\end{Rk}


\subsection{\textsf{Explicit results for the infinite series --- Types $A,\, B,\, C,\, D$}}

\subsubsection{\textsf{The $A$ Series}}

\noindent In this section, we consider $\Gamma$ a finite diagonal
abelian subgroup of $\mathbf{SL}_3\mathbb{C}$. Then $\Gamma$ is
isomorphic to a product of cyclic groups:
$$\Gamma\simeq\mathbb{Z}/j_1\mathbb{Z}\times\dots\times
\mathbb{Z}/j_k\mathbb{Z}.$$

\noindent If $\Gamma$ is a finite subgroup of
$\mathbf{SL}_r\mathbb{C}$, then $\Gamma$ is a \emph{small}
subgroup of $\mathbf{GL}_r\mathbb{C}$, i.e. no element of
$\Gamma$ has an eigenvalue $1$ of multiplicity $r-1$. In fact, if
$g\in\Gamma$ has an eigenvalue $1$ of multiplicity $r-1$, then the
last eigenvalue of $g$ is different from $1$ and the
determinant of $g$ is also different from $1$, which is impossible.\\
Then, according to a lemma of \cite{DHZ05} (p.13), $\Gamma$ has at most $r-1$ generators.
So, for a subgroup $\Gamma$ of type $A$, we may assume
that $k\leq 2$, i.e. we have two cases:\\
$(A_1)$ $\Gamma\simeq \mathbb{Z}/j\mathbb{Z}$,\\
$(A_2)$ $\Gamma\simeq \mathbb{Z}/j_1\mathbb{Z}\times \mathbb{Z}/j_2\mathbb{Z}$, with $j_1\geq j_2\geq 2$.\\

\subsubsection{\textsf{Type $A_1$}}

\noindent $\bullet$ Here, we take $\Gamma=\mathbb{Z}/j\mathbb{Z}$.
The natural representation and the natural character of $\Gamma$ are $$\begin{array}{rcl|rcl}
  \gamma\ :\ \mathbb{Z}/j\mathbb{Z} & \rightarrow & \mathbf{SL}_3\mathbb{C} & \chi\ :\ \mathbb{Z}/j\mathbb{Z} & \rightarrow & \mathbf{SL}_3\mathbb{C}\\
  \overline{k} & \mapsto & \left(%
\begin{array}{ccc}
  \zeta_j^k & 0 & 0 \\
  0 & 1 & 0 \\
  0 & 0 & \zeta_j^{-k} \\
\end{array}%
\right) & \overline{k} & \mapsto & 1+\zeta_j^{k}+\zeta_j^{-k}.\\
\end{array}$$
The character table of $\Gamma$ is $T_\Gamma=\left(\zeta_j^{kl}\right)_{(k,l)\in\llbr
0,\,j-1\rrbr}$. Let $\sigma\in\go{S}_{\llbr 0,\,j-1\rrbr}$ be the permutation defined by
$\sigma(0)=0$ and $\forall\ i\in \llbr 0,\,j-1\rrbr,\
\sigma(i)=j-i$. Then $T_\Gamma^{-1}=\frac{1}{j}\,T_\Gamma Q^\sigma.$ The eigenvalues of $A^{(1)}$ are the
numbers $\overline{\chi(k)}=1+\zeta_j^k+\zeta_j^{-k}$, for
$k\in\llbr 0,\,j-1\rrbr$. According to Formula \ref{fond},
$$P_\Gamma(t,\,u)_i=\frac{1}{j}\,\left(T_\Gamma\Delta(t,\,u)T_\Gamma
Q^\sigma\right)_{i0}=\frac{1}{j}\,\sum_{p=0}^{j-1}\frac{\zeta_j^{ip}(1-tu)}{(1-t)(1-t\zeta_j^p)(1-t\zeta_j^{-p})
(1-u)(1-u\zeta_j^p)(1-u\zeta_j^{-p})}.$$ Note that
$(1-t^j)(1-t^2)(1-u^j)(1-u^2)$ is a common denominator of all the
terms of the preceding sum, which is independent of $i$.\\

\noindent $\bullet$ The matrix $A^{(1)}\in\mathbf{M}_j\mathbb{C}$ is $$A^{(1)}=\left(%
\begin{array}{cc}
  1 & 2 \\
  2 & 1 \\
\end{array}%
\right)\ \textrm{if}\ j=2,\ \ A^{(1)}=\left(%
\begin{array}{ccccc}
  1 & 1      &        & 1  \\
  1 & \ddots & \ddots &    \\
    & \ddots & \ddots & 1  \\
  1 &        & 1      & 1 \\
\end{array}%
\right)\ \textrm{if}\ j\geq 3.$$
Then the set of reflections $\mathcal{S}$ may be
decomposed in two (resp. three) sets if $j$ is even (resp. odd).\\
$\triangleright$ If $j$ is even, we have $\tau_0=s_0s_4\dots s_{j-2}$,
$\tau_1=s_1s_3\dots s_{j-1}$, and
$C_{A_1}(j)=2I_{j-1}-(\tau_0+\tau_1)$.\\
$\triangleright$ If $j$ is odd, we have $\tau_0=s_2s_4\dots s_{j-1}$,
$\tau_1=s_1s_3\dots s_{j-2}$, $\tau_2=s_0$, and
$C_{A_1}(j)=3I_j-(\tau_0+\tau_1+\tau_2)$.\\
The graph associated to $C_{A_1}(j)$ is a cyclic graph with $j$ vertices and $j$ edges.\\

\subsubsection{\textsf{Type $A_2$}}

\noindent $\bullet$ We now consider the case
$\Gamma=\mathbb{Z}/j_1\mathbb{Z}\times \mathbb{Z}/j_2\mathbb{Z}$, with $j_1\geq j_2\geq 2$. The natural representation and
 the natural character of $\Gamma$ are
$$\begin{array}{rcl|rcl}
  \gamma\ :\ \mathbb{Z}/j_1\mathbb{Z}\times \mathbb{Z}/j_2\mathbb{Z} & \rightarrow & \mathbf{SL}_3\mathbb{C} & \chi\ :\ \mathbb{Z}/j_1\mathbb{Z}\times \mathbb{Z}/j_2\mathbb{Z} & \rightarrow & \mathbf{SL}_3\mathbb{C}\\
  (\overline{k_1},\,\overline{k_1}) & \mapsto & \left(%
\begin{array}{ccc}
  \zeta_{j_1}^{k_1} & 0 & 0 \\
  0 & \zeta_{j_2}^{k_2} & 0 \\
  0 & 0 & \zeta_{j_1}^{-k_1}\zeta_{j_2}^{-k_2} \\
\end{array}%
\right) & (\overline{k_1},\,\overline{k_2}) & \mapsto & \zeta_{j_1}^{k_1}+\zeta_{j_2}^{k_2}+\zeta_{j_1}^{-k_1}\zeta_{j_2}^{-k_2}.\\
\end{array}$$
The irreducible characters of $\Gamma$ are the elements of the form
$\chi_1\otimes\chi_2$, where $\chi_1$ and $\chi_2$ are irreducible characters of
$\mathbb{Z}/j_1\mathbb{Z}$ and $\mathbb{Z}/j_2\mathbb{Z}$, i.e.
the irreducible characters of $\Gamma$ are, for $(l_1,\,l_2)\in\llbr 0,\,j_1-1\rrbr\times\llbr
0,\,j_2-1\rrbr$,
$$\begin{array}{rcl}
  \chi_{l_1,l_2}\ :\ \mathbb{Z}/j_1\mathbb{Z}\times \mathbb{Z}/j_2\mathbb{Z} & \rightarrow & \mathbf{SL}_3\mathbb{C}\\
  (\overline{k_1},\,\overline{k_2}) & \mapsto & \zeta_{j_1}^{k_1l_1}\zeta_{j_2}^{k_2l_2}\\
\end{array}.$$
For $k\in\{1,\,2\}$, let us denote by $T_k$ the character table of
the group $\mathbb{Z}/j_k\mathbb{Z}$. Then the character table of
$\Gamma=\mathbb{Z}/j_1\mathbb{Z}\times \mathbb{Z}/j_2\mathbb{Z}$
is the Kronecker product\footnote{Recall that the Kronecker
product of two matrices $A\in\mathbf{M}_m\mathbb{C}$ and
$B\in\mathbf{M}_n\mathbb{C}$ is the block-matrix $A\otimes
B\in\mathbf{M}_{mn}\mathbb{C}$ defined by the formula: $$\forall\
(i,\,j)\in\llbr 1,\,m\rrbr,\ (A\otimes B)_{ij}=a_{ij}B.$$ An
important property of the Kronecker product is the relation
$$\textrm{tr}(A\otimes B)=\textrm{tr}(A)\textrm{tr}(B).$$ The
equality $T_\Gamma=T_1\otimes T_2$ is implied by this relation} $T_\Gamma=T_1\otimes T_2.$ Let
$\sigma_k\in\go{S}_{\llbr 0,\,j_k-1\rrbr}$ be the permutation defined by
$\sigma_k(0)=0$ and $\forall\ i\in \llbr 0,\,j_k-1\rrbr,\
\sigma_k(i)=j_k-i$. We have
$$T_\Gamma^{-1}=(T_1\otimes T_2)^{-1}=\frac{1}{j_1j_2}(T_1\otimes T_2)(Q^{\sigma_1}\otimes
Q^{\sigma_2})=\frac{1}{j_1j_2}(T_1Q^{\sigma_1})\otimes(T_2Q^{\sigma_2}).$$
The eigenvalues of $A^{(1)}$ are the numbers
$\overline{\chi(k_1,\,k_2)}=\zeta_{j_1}^{-k_1}+\zeta_{j_2}^{-k_2}+\zeta_{j_1}^{k_1}\zeta_{j_2}^{k_2},$
for $(k_1,\,k_2)\in\llbr 0,\,j_1-1\rrbr\times \llbr
0,\,j_2-1\rrbr$.\\
Let us denote by $\Lambda^{(1)}:=\textrm{Diag}\left(\Lambda^{(1)}_0,\dots,\,\Lambda^{(1)}_{j_1}\right)$ the diagonal block-matrix defined
by $$(\Lambda^{(1)}_{k_1})_{k_2k_2}=\overline{\chi(k_1,\,k_2)}=\zeta_{j_1}^{-k_1}+\zeta_{j_2}^{-k_2}+\zeta_{j_1}^{k_1}\zeta_{j_2}^{k_2}.$$
According to Formula \ref{fond}, for $(m,\,n)\in\llbr 0,\,j_1-1\rrbr\times\llbr
0,\,j_2-1\rrbr$, we have
$$\begin{array}{rl}
    & P_\Gamma(t,\,u)_{mj_2+n}\vspace{.2cm} \\
  = & \displaystyle{\frac{1}{j_1j_2}\,\sum_{k=0}^{j_1-1}\sum_{l=0}^{j_2-1}\zeta_{j_1}^{mk}\zeta_{j_2}^{nl}(1-tu)
\left(1-t(\zeta_{j_1}^{-k}+\zeta_{j_2}^{-l}+\zeta_{j_1}^{k}\zeta_{j_2}^{l})
+t^2(\zeta_{j_1}^{k}+\zeta_{j_2}^{l}+\zeta_{j_1}^{-k}\zeta_{j_2}^{-l})-t^3\right)^{-1}}\\
  & \displaystyle{\ \ \ \ \ \ \ \ \ \ \ \ \ \ \ \ \ \ \ \ \ \ \ \ \
  \left(1-u(\zeta_{j_1}^{k}+\zeta_{j_2}^{l}+\zeta_{j_1}^{-k}\zeta_{j_2}^{-l})
  +u^2(\zeta_{j_1}^{-k}+\zeta_{j_2}^{-l}+\zeta_{j_1}^{k}\zeta_{j_2}^{l})-u^3\right)^{-1}.} \\
\end{array}$$

\noindent $\bullet$ The matrix $A^{(1)}$ is a block-matrix with $j_1^2$ blocs of size
$j_2$, and we have
\begin{small}
$$A^{(1)}=\mathbf{1}_{4,4}-I_4,
\ \textrm{if}\ j_1=j_2=2,\ \textrm{and}\ A^{(1)}=\left(%
\begin{array}{ccccc}
  Q^{\rho_2}   & ^tQ^{\rho_2} &        & I_{j_2}  \\
  I_{j_2}      & \ddots       & \ddots &          \\
               & \ddots       & \ddots & ^tQ^{\rho_2} \\
  ^tQ^{\rho_2} &              & I_{j_2}& Q^{\rho_2} \\
\end{array}%
\right)\ \textrm{if}\ j_1\geq 2,\ Q^{\rho_k}:=\left(%
\begin{array}{cccc}
  0  &        &        & 1  \\
  1  & \ddots &        &    \\
     & \ddots & \ddots &    \\
     &        &   1    & 0 \\
\end{array}%
\right).$$
\end{small}
So, we may write $A^{(1)}=I_{j_1}\otimes
Q^{\rho_2}+Q^{\rho_1}\otimes I_{j_2}+\,^tQ^{\rho_1}\otimes\,
^tQ^{\rho_2}.$\\
Note that $Diag(A^{(1)})=0$. Then
$C_{A_2}(j_1,\,j_2)=I_{j_1}\otimes W+Q^{\tau_1}\otimes
S+\,^tQ^{\tau_1}\otimes\, ^tS.$\\

\noindent $\bullet$ Now, let us decompose the matrix $A^{(1)}$:\\

\noindent If $j_1=j_2=2$, then the decomposition of $C_{A_2}(2,\,2)$ is $C_{A_2}(2,\,2)=4\,I_4-(s_0+s_1+s_2+s_3)$.\\
Now, we assume that $j_1\geq 3$. For $(i_1,\,i_2)\in\llbr0,\,j_1-1\rrbr\times\llbr0,\,j_2-1\rrbr$, let $s_{i_1,i_2}$ be the reflection
associated to the $(i_1j_2+i_2)-$th root. Then the set
$\mathcal{S}$ may be decomposed into $p$ sets where $p\in\{4,\,6,\,9\}$.\\
For $i_1\in\llbr 0,\, j_1-1\rrbr$, define $\widetilde{S_{i_1}}:=\{s_{i_1,0},\dots,\,s_{i_1,j_2-1}\}$ .\\
$\triangleright$ If $j_1$ is odd, we set
$$\widehat{I_0}:=\{0,\,2,\dots,\,j_1-3\},\ \widehat{I_1}:=\{1,\,3,\dots,\,j_1-2\},\ \widehat{I_0}:=\{j_1-1\}.$$
$\triangleright$ If $j_1$ is even, we set
$$\widehat{I_0}:=\{0,\,2,\dots,\,j_1-2\},\ \widehat{I_1}:=\{1,\,3,\dots,\,j_1-1\}.$$
\noindent Then, the roots associated to the reflections of distinct
$\widetilde{S_{i_1}}$'s for $i_1$ belonging to a same $\widehat{I_k}$ are orthogonal.\\
Now, we decompose each $\widetilde{S_{i_1}}$, i.e. $\widetilde{S_{i_1}}=\widehat{S_{i_1,0}}\sqcup\dots\sqcup\widehat{S_{i_1,q-1}}$,
such that $q\in\{2,\,3\}$ and for every $k\in\llbr 0,\,q-1\rrbr$, the roots associated to the reflections belonging to $\widehat{S_{i_1,k}}$ are orthogonal:\\
$\triangleright$ If $j_2$ is odd, then $\widetilde{S_{i_1}}=\widehat{S_{i_1,0}}\sqcup\widehat{S_{i_1,1}}\sqcup\widehat{S_{i_1,2}}$, with
$$\widehat{S_{i_1,0}}=\{s_{i_1,0},\,s_{i_1,2},\dots,\,s_{i_1,j_2-3}\},\
\widehat{S_{i_1,1}}=\{s_{i_1,1},\,s_{i_1,3},\dots,\,s_{i_1,j_2-2}\},\
\widehat{S_{i_1,2}}=\{s_{i_1,j_2-1}\}.$$
$\triangleright$ If $j_2$ is even, then $\widetilde{S_{i_1}}=\widehat{S_{i_1,0}}\sqcup\widehat{S_{i_1,1}}$, with
$$\widehat{S_{i_1,0}}=\{s_{i_1,0},\,s_{i_1,2},\dots,\,s_{i_1,j_2-2}\},\
\widehat{S_{i_1,1}}=\{s_{i_1,1},\,s_{i_1,3},\dots,\,s_{i_1,j_2-1}\}.$$

\noindent Finally, for $(k,\,l)\in\{0,\,1,\,2\}^2$, we set $$S_{k,l}:=\coprod_{i_1\in \widehat{I_k}}\widehat{S_{i_1,l}},$$
and we denote by $p\in\{4,\,6,\,9\}$ the number of non-empty sets $S_{k,l}$, and by $\tau_{k,l}$ the commutative product of the reflections of $S_{k,l}$. Then, we have $$C_{A_2}(j_1,\,j_2):=2\,I_{j_1j_2}-A^{(1)}-\,^tA^{(1)}+2\,Diag(A^{(1)})=p\,I_{j_1j_2}-\sum_{(k,l)\in\{0,\,1,\,2\}^2}\tau_{k,l}.$$\\

\noindent $\bullet$ If $j_1\geq 3$, the graph associated to $\Gamma$ is  a $j_1-$gone, such that every vertex
of this $j_1-$gone is a $j_2-$gone, and every vertex of each $j_2-$gone is connected with exactly $2$ vertices of
both adjacent $j_2-$gones (for $j_1=j_2=2$, see Remark \ref{A1caspart}).\\

\begin{Rk}\label{A1caspart}$\\$
In the cases $j_1=2,\,j_2=2$ and $j_1=3,\,j_2=2$, we obtain full
matrices and complete graphs. Moreover the complete graphs with
$4$ and $6$ vertices are the unique complete graphs that we can
obtain for the type $A_2$ (the complete graphs with $2$ and $3$
vertices are the unique complete graphs that we can obtain for the
type $A_1$).\\
\end{Rk}

\begin{Ex}$\\$
We consider the case where $j_1=6$ and $j_2=5$. Then the
decomposition of the Cartan matrix is
$$C_{A_2}(6,\,5)=6\,I_{30}-(\tau_{0,0}+\tau_{0,1}+\tau_{0,2}+\tau_{1,0}+\tau_{1,1}+\tau_{1,2}),$$
where the matrix $C_{A_2}(6,\,5)$, the $\tau_{i,j}$'s and the graph associated to
$C_{A_2}(6,\,5)$ are given by Figure \ref{ExA65}.
\begin{figure}[h]
\begin{center}
\begin{tabular}{ccc}
$\begin{array}{l}\\[-4cm]
  \tau_{0,0}:=(s_{0,0}s_{0,2})(s_{2,0}s_{2,2})(s_{4,0}s_{4,2})\\
  \tau_{0,1}:=(s_{0,1}s_{0,3})(s_{2,1}s_{2,3})(s_{4,1}s_{4,3})\\
  \tau_{0,2}:=(s_{0,4})(s_{2,4})(s_{4,4})\\
  \tau_{1,0}:=(s_{1,0}s_{1,2})(s_{3,0}s_{3,2})(s_{5,0}s_{5,2}),\\
  \tau_{1,1}:=(s_{1,1}s_{1,3})(s_{3,1}s_{3,3})(s_{5,1}s_{5,3}),\\
  \tau_{1,2}:=(s_{1,4})(s_{3,4})(s_{5,4}).\\
\end{array}$ & \includegraphics[width=4cm,height=4cm]{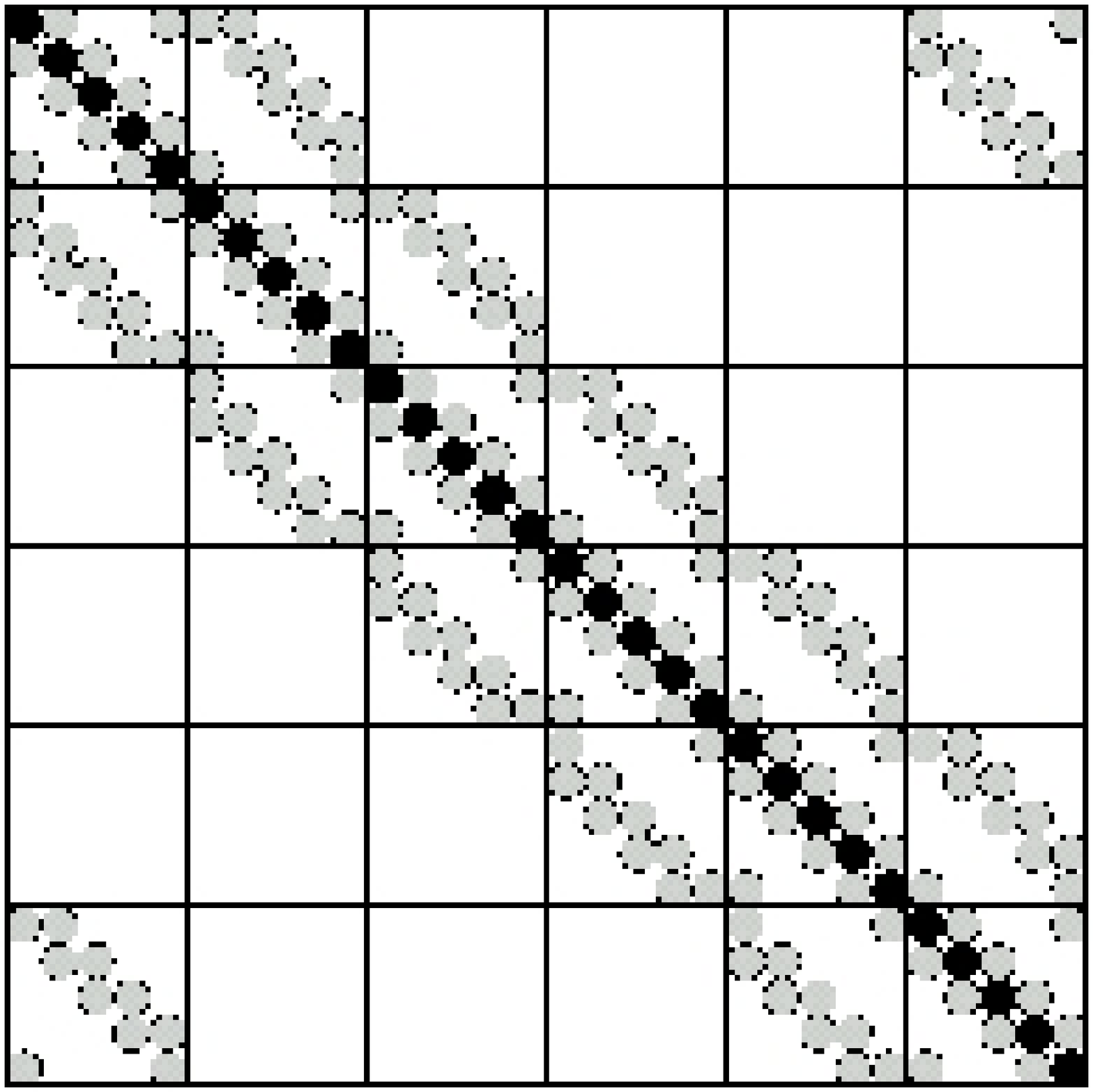} & \includegraphics[width=4cm,height=4cm]{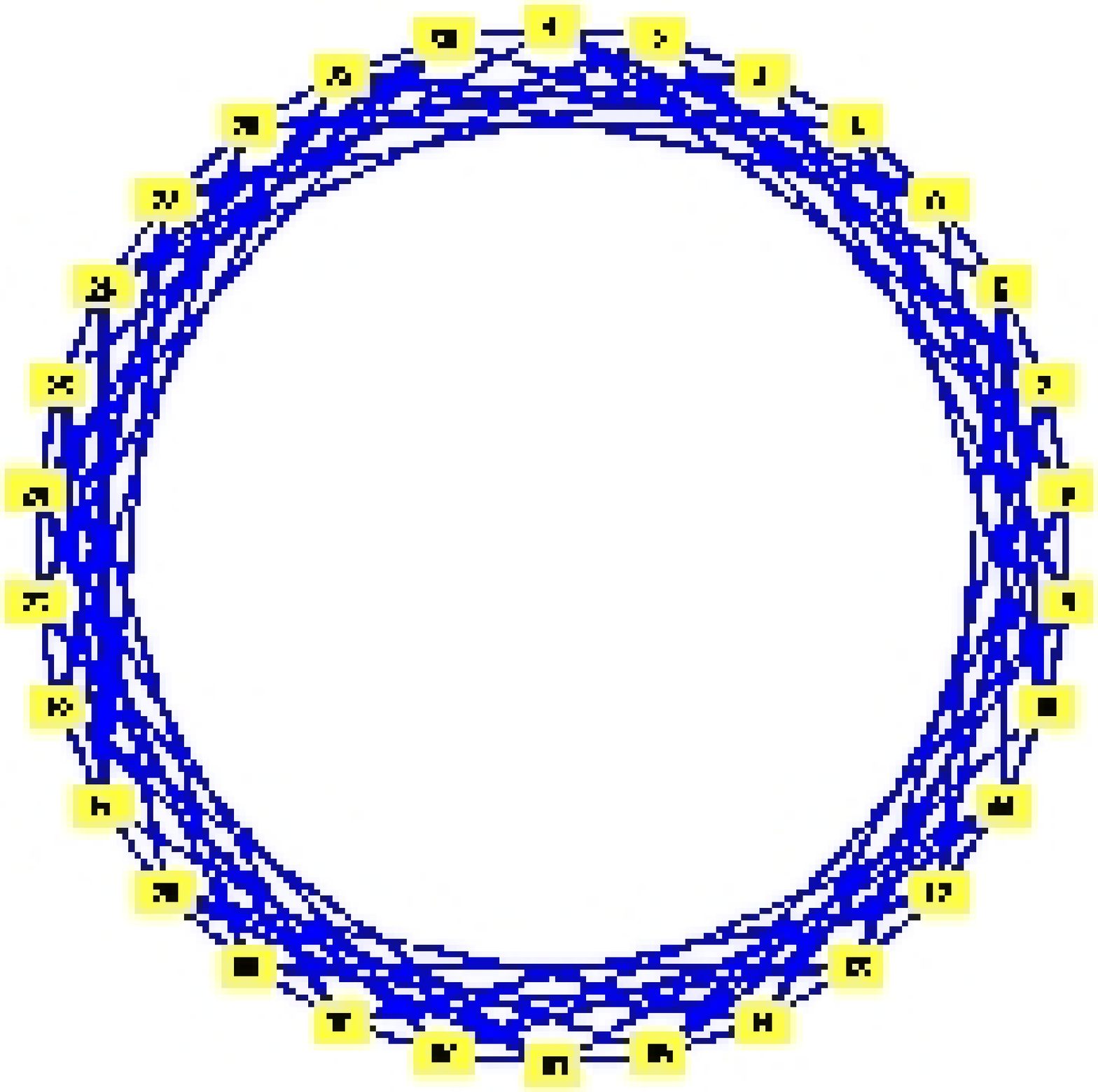}\\
 \end{tabular}
\end{center}
\caption{Matrix $C_{A_2}(6,\,5)$ and corresponding graph.}\label{ExA65}
\end{figure}
\end{Ex}


\subsection{\textsf{The $B$ series}}\label{section542}

\noindent In this section, we study the binary groups of $\mathbf{SL}_3\mathbb{C}$.\\
 We give a general formula for the types $BDa$, $BTa$, $BO$ and $BI$. In all these cases, the group $\Gamma$ contains two normal subgroups
 $\Gamma_1$ and $\Gamma_2$ such that $\Gamma_1\cap\Gamma_2=\{id\}$, and $|\Gamma_1|\cdot|\Gamma_2|=|\Gamma|$,
so that $\Gamma=\Gamma_2\Gamma_1$ and $\Gamma\simeq\Gamma_2\times\Gamma_1$. The group $\Gamma_1$ is isomorphic to a binary group of $\mathbf{SL}_2\mathbb{C}$ and $\Gamma_2$ is isomorphic to $\mathbb{Z}/m\mathbb{Z}$. So, we can deduce the results for $\Gamma$ of the results obtained for the group $\Gamma_1$.\\
If we denote by $T_k$ the character table of
the group $\Gamma_k$, the character table of the direct product
$\Gamma=\Gamma_2\times \Gamma_1$ is the Kronecker product $T_\Gamma=T_2\otimes T_1$. The matrix $T_2$ is given in Section 2.2.1 (Type $A$ --- Cyclic groups), and the matrix $T_1$ is given in the section dealing with the corresponding binary group of $\mathbf{SL}_2\mathbb{C}$. We also have $T_\Gamma^{-1}=T_2^{-1}\otimes T_1^{-1}=\frac{1}{m}(T_2Q^{\sigma_2})\otimes T_1^{-1},$
where $\sigma_2$ is the permutation matrix defined by
$\sigma_2(0)=0$, and $\forall\ i\in \llbr 0,\,m-1\rrbr,\
\sigma_2(i)=m-i$.\\
Let us denote by $h$ the number of conjugacy classes of $\Gamma_1$. The columns of $T_\Gamma$ give a basis of eigenvectors and the eigenvalues of $A^{(1)}$ are the
numbers $\overline{\chi_{i,j}},\ (i,\,j)\in\llbr 0,\,m-1\rrbr\times\llbr 0,\,h-1\rrbr$ where $\chi_{i,j}$ is the value of the natural character of $\Gamma$ on the $(i,\,j)-$th conjugacy class.\\
Let us also denote by $\Lambda^{(1)}$ the diagonal block-matrix with $m\times m$ blocks of size $h\times h$ defined by
$$\Lambda^{(1)}:=\textrm{Diag}\left(\Lambda^{(1)}_0,\dots,\,\Lambda^{(1)}_{m-1}\right),\ (\Lambda^{(1)}_{i})_{jj}=\overline{\chi_{i,\,j}}.$$
According to Formula \ref{fond},
$$P_\Gamma(t,\,u)=T_\Gamma\Delta(t,\,u)T_\Gamma^{-1}v_{0,0}=
\frac{1}{m}\,(T_2\otimes
T_1)\Delta(t,\,u)((T_2Q^{\sigma_2})\otimes
(T_1T_1^{-2}))v_{0,0},$$ where $\Delta(t,\,u)$ is the diagonal
block-matrix defined by
$$\Delta(t,\,u)=\textrm{Diag}\left(\Delta(t,\,u)^{(0)},\dots,\,\Delta(t,\,u)^{(m-1)}\right),$$
with $\begin{array}{rcllcl} \Delta(t,\,u)^{(i)}_{jj} & = & f(\overline{\chi_{i,j}},\,\chi_{i,j}) & = &
\displaystyle{\frac{1-tu}{(1-t\overline{\chi_{i,j}}+t^2\chi_{i,j}-t^3)
(1-u\chi_{i,j}+u^2\overline{\chi_{i,j}}-u^3)}}=:f_{ij}. \\
\end{array}$\\

\noindent The decomposition of the matrix $A^{(1)}$ and the description of the associated graph are made in the same way for the binary tetrahedral, octahedral and icosahedral groups: the results are collected in \ref{decomposition}.\\

\subsubsection{\textsf{The  $BDa$ subseries --- Binary dihedral groups}}

\noindent $\bullet$ For $(q,\,m)\in\mathbb{N}^2$, let $\psi_{2q}$, $\tau$ and $\phi_{2m}$ be the
 following elements of $\mathbf{SL}_3\mathbb{C}$:

$$\psi_{2q}=\left(
              \begin{array}{ccc}
                1 & 0 & 0\\
                0 & \zeta_{2q} & 0 \\
                0 & 0 & \zeta_{2q}^{-1} \\
              \end{array}
            \right)
,\ \tau=\left(
          \begin{array}{ccc}
            1 & 0 & 0 \\
            0 & 0 & i \\
            0 & i & 0 \\
          \end{array}
        \right)
,\ \phi_{2m}=\left(
               \begin{array}{ccc}
                 \zeta_{2m}^{-2} & 0 & 0 \\
                 0 & \zeta_{2m} & 0 \\
                 0 & 0 & \zeta_{2m} \\
               \end{array}
             \right).
$$
\noindent In this section\footnote{The other case --- the type $BDb$ --- is $1<q<n$, $n\wedge q=1$, and $m:=n-q\equiv 0 \mod 2$. This group is not a direct product, the general expression for groupes in this subseries is unclear.}, we assume that $1<q<n$, $n\wedge q=1$, and $m:=n-q\equiv 1 \mod 2$, and we consider the subgroup
$\Gamma:=\langle\psi_{2q},\ \tau,\ \phi_{2m}\rangle$ of $\mathbf{SL}_3\mathbb{C}$. Note that $\phi_{2m}=\psi_{2q}^q\phi_m^{-\frac{m-1}{2}}$, so that $\Gamma=\langle\psi_{2q},\ \tau,\ \phi_{m}\rangle$. Set~$\Gamma_1:=\langle\psi_{2q},\ \tau\rangle$ and $\Gamma_2:=\langle\phi_{m}\rangle\simeq \mathbb{Z}/m\mathbb{Z}$. Then $\Gamma\simeq \Gamma_2\times\Gamma_1$. With the notations used for $\mathbf{SL}_2\mathbb{C}$, $\psi_{2q}$ (resp. $\tau$) represents $a_q$ (resp. $b$),
where $\langle a_q,\,b\rangle$ is the binary dihedral subgroup of $\mathbf{SL}_2\mathbb{C}$.\\

\noindent The natural character of $\Gamma$ is given by $\chi=(\chi_i)_{i=0\dots m-1}$, with
$$\begin{array}{rcl}
\chi_{i} & = & \left[\chi_{i,0},\ \chi_{i,1},\ \chi_{i,2},\ \chi_{i,3},\ \chi_{i,4},\ \chi_{i,5},\ \chi_{i,6},\dots,\ \chi_{i,q+2}\right]\\
\\
     & = & \left[tr(\phi_m^i),\ tr(\phi_m^i\psi_{2q}\tau),\ tr(\phi_m^i\tau),\ tr(\phi_m^i\psi_{2q}^q)
    ,\ tr(\phi_m^i\psi_{2q}),\ tr(\phi_m^i\psi_{2q}^2),\ tr(\phi_m^i\psi_{2q}^3),\dots,\ tr(\phi_m^i\psi_{2q}^{q-1})\right] \\
    \\
     & = &  \Big[\zeta_m^{-2i}+2\zeta_m^i,\ \zeta_m^{-2i},\ \zeta_m^{-2i},\ \zeta_m^{-2i}-2\zeta_m^i
    ,\ \zeta_m^{-2i}+\zeta_m^i(\zeta_{2q}+\zeta_{2q}^{-1}),\ \zeta_m^{-2i}+\zeta_m^i(\zeta_{2q}^2+\zeta_{2q}^{-2}),\\
     & & \zeta_m^{-2i}+\zeta_m^i(\zeta_{2q}^3+\zeta_{2q}^{-3}),\dots,\ \zeta_m^{-2i}+\zeta_m^i(\zeta_{2q}^{q-1}+\zeta_{2q}^{-(q-1)})\Big].
  \end{array}
$$
We deduce the following formula for the series $P_\Gamma(t,\,u)$:
$$\begin{array}{rcl}
    P_\Gamma(t,\,u)_{i_1(q+3)} & = & \displaystyle{\frac{1}{m}\sum_{r=0}^{m-1}\zeta_m^{i_1r}\Bigg[\frac{3q-1}{8q^2}\left(d_1^{(r)}+d_2^{(r)}+d_3^{(r)}+d_4^{(r)}
    +2\sum_{k=1}^{q-1}\delta_k^{(r)}\right)}\\
     & & \displaystyle{+(-1)^q\,\frac{q-1}{8q^2}\left(d_1^{(r)}+(-1)^q(d_2^{(r)}+d_3^{(r)})+d_4^{(r)}+2\sum_{k=1}^{q-1}(-1)^k\delta_k^{(r)}\right)} \\
     &  & \displaystyle{+\sum_{l=1}^{q-1}(-1)^l\,\frac{q-1}{4q^2}\left(d_1^{(r)}+(-1)^l(d_2^{(r)}+d_3^{(r)})+d_4^{(r)}
     +\sum_{k=1}^{q-1}(\zeta_{2q}^{lk}+\zeta_{2q}^{-lk})\delta_k^{(r)}\right)\Bigg]},
  \end{array}
$$
$$\begin{array}{rcl}
    P_\Gamma(t,\,u)_{i_1(q+3)+1} & = & \displaystyle{\frac{1}{m}\sum_{r=0}^{m-1}\zeta_m^{i_1r}\Bigg[\frac{3q-1}{8q^2}\left(d_1^{(r)}+i^nd_2^{(r)}-i^nd_3^{(r)}+(-1)^nd_4^{(r)}
    +2\sum_{k=1}^{q-1}(-1)^k\delta_k^{(r)}\right)}
     \\
     & & \displaystyle{+(-1)^q\,\frac{q-1}{8q^2}\left(d_1^{(r)}+(-1)^q(i^nd_2^{(r)}-i^nd_3^{(r)}+d_4^{(r)})+2\sum_{k=1}^{q-1}\delta_k^{(r)}\right)} \\
     &  & \displaystyle{+\sum_{l=1}^{q-1}(-1)^l\,\frac{q-1}{4q^2}\left(d_1^{(r)}+(-1)^l(i^nd_2^{(r)}-i^nd_3^{(r)})+(-1)^nd_4^{(r)}
     +\sum_{k=1}^{q-1}(-1)^k(\zeta_{2q}^{lk}+\zeta_{2q}^{-lk})\delta_k^{(r)}\right)\Bigg]},
  \end{array}
$$
and $P_\Gamma(t)_{i_1(q+3)+3}$ (resp. $P_\Gamma(t)_{i_1(q+3)+2}$) is obtained by replacing in $P_\Gamma(t)_{i_1(q+3)}$ (resp. $P_\Gamma(t)_{i_1(q+3)+1}$) $d_2$ by $-d_2$ and $d_3$ by $-d_3$. Finally, for $i_2\in\llbr 1,\,q-1\rrbr$, we have
$$\begin{array}{rcl}
    P_\Gamma(t,\,u)_{i_1(q+3)+i_2+3} & = & \displaystyle{\frac{1}{m}\sum_{r=0}^{m-1}\zeta_m^{i_1r}\Bigg[\frac{3q-1}{8q^2}\left(2d_1^{(r)}+2(-1)^{i_2}d_4^{(r)}
    +2\sum_{k=1}^{q-1}(\zeta_{2q}^{ki_2}+\zeta_{2q}^{-ki_2})\delta_k^{(r)}\right)}
     \\
     & & \displaystyle{+(-1)^q\,\frac{q-1}{8q^2}\left(2d_1^{(r)}+2(-1)^{i_2}d_4^{(r)}
     +2\sum_{k=1}^{q-1}(-1)^k(\zeta_{2q}^{i_2k}+\zeta_{2q}^{-i_2k})\delta_k^{(r)}\right)} \\
     &  & \displaystyle{+\sum_{l=1}^{q-1}(-1)^l\,\frac{q-1}{4q^2}\left(2d_1^{(r)}+2(-1)^{i_2}d_4^{(r)}
     +\sum_{k=1}^{q-1}(\zeta_{2q}^{ki_2}+\zeta_{2q}^{-ki_2})\delta_k^{(r)}(\zeta_{2q}^{kl}+\zeta_{2q}^{-kl})\right)\Bigg]}.
  \end{array}
$$

\noindent $\bullet$ We now make the matrix $A^{(1)}$ explicit: $A^{(1)}$ is a block-matrix with $m\times m$ blocks of size $(q+3)\times(q+3)$.\\
$\triangleright$ If $m\geq 5$, then the matrices $A^{(1)}$ and $C_\Gamma:=2I-A^{(1)}-\,^tA^{(1)}$ are defined by
$$A^{(1)}=\left(
    \begin{array}{ccccc}
       0 &   & I   &   & B \\
      B & \ddots  &   & \ddots  &   \\
        & \ddots  &  \ddots  &    & I  \\
      I &   & \ddots  &  \ddots  &   \\
       & I &   & B & 0  \\
    \end{array}
  \right),\ \ C_\Gamma=\left(
    \begin{array}{ccccccc}
       2I & -B & -I  &   &   & -I & -B \\
      -B & 2I  & -B & -I  &   &   & -I \\
       -I & -B & 2I  & -B &  -I &   &   \\
        &  -I & -B & \ddots  & \ddots &  \ddots &   \\
        &   &  -I & \ddots & \ddots  & -B & -I  \\
      -I &   &   & \ddots  & -B &  2I & -B \\
      -B & -I &   &   & -I  & -B &  2I \\
    \end{array}
  \right),
$$
with
$$B=\left(
                                     \begin{array}{ccccc}
                                       0 & 0 & 0 & 0 & 1 \\
                                       0 & 0 & 0 & 0 & 1 \\
                                       0 & 0 & 0 & 0 & 1 \\
                                       0 & 0 & 0 & 0 & 1 \\
                                       1 & 1 & 1 & 1 & 0 \\
                                     \end{array}
                                   \right),\ \textrm{if}\ q=2,\ \textrm{and}\
B=\left(
                                        \begin{array}{ccc|ccccc}
                                          0 & 0 & 0 & 0 & 1 & 0 & 0 & 0 \\
                                           0 & 0 & 0 & 0 & 0 & 0 & 0 & 1 \\
                                          0 & 0 & 0 & 0 & 0 & 0 & 0 & 1 \\ \hline
                                          0 & 0 & 0 & 0 & 1 & 0 & 0 & 0 \\
                                          1 & 0 & 0 & 1 & 0 & 1 & 0 & 0 \\
                                          0 & 0 & 0 & 0 & 1 & 0 & \ddots & 0 \\
                                          0 & 0 & 0 & 0 & 0 & \ddots & \ddots & 1 \\
                                          0 & 1 & 1 & 0 & 0 & 0 & 1 & 0 \\
                                        \end{array}
                                      \right)\ \textrm{if}\ q\geq 3
$$\\
$\triangleright$ If $m=3$, then the matrices $A^{(1)}$ and $C_\Gamma:=2I-A^{(1)}-\,^tA^{(1)}$ are defined by
$$A^{(1)}=\left(
            \begin{array}{ccc}
              0 & 0 & B \\
              B & 0 & 0 \\
              0 & B & 0 \\
            \end{array}
          \right),\ \ C_\Gamma=
          \left(
            \begin{array}{ccc}
              0 & B & B \\
              B & 0 & B \\
              B & B & 0 \\
            \end{array}
          \right),
$$
with $$B=\left(
                                     \begin{array}{ccccc}
                                       1 & 0 & 0 & 0 & 1 \\
                                       0 & 1 & 0 & 0 & 1 \\
                                       0 & 0 & 1 & 0 & 1 \\
                                       0 & 0 & 0 & 1 & 1 \\
                                       1 & 1 & 1 & 1 & 1 \\
                                     \end{array}
                                   \right),\ \textrm{if}\ q=2,\ \textrm{and}\
B=\left(
                                        \begin{array}{ccc|ccccc}
                                          1 & 0 & 0 & 0 & 1 & 0 & 0 & 0 \\
                                           0 & 1 & 0 & 0 & 0 & 0 & 0 & 1 \\
                                          0 & 0 & 1 & 0 & 0 & 0 & 0 & 1 \\ \hline
                                          0 & 0 & 0 & 1& 1 & 0 & 0 & 0 \\
                                          1 & 0 & 0 & 1 & 1 & 1 & 0 & 0 \\
                                          0 & 0 & 0 & 0 & 1 & 1 & \ddots & 0 \\
                                          0 & 0 & 0 & 0 & 0 & \ddots & \ddots & 1 \\
                                          0 & 1 & 1 & 0 & 0 & 0 & 1 & 1 \\
                                        \end{array}
                                      \right)\ \textrm{if}\ q\geq 3.
$$\\
$\triangleright$ If $m=1$, then the matrices $A^{(1)}$ and $C_\Gamma$ are
$A^{(1)}=B$ and  $C_\Gamma:=2I-A^{(1)}-\,^tA^{(1)}+2\,Diag(A^{(1)})$, with $B$ defined as in the case $m=3$.\\

\noindent $\bullet$ For $i_1\in\llbr0,\,m-1\rrbr$ and
$i_2\in\llbr 0,\,q+2\rrbr$, let $s_{i_1,i_2}$ be the reflection
associated to the $(i_1(q+3)+i_2)-$th root. Then the set
$\mathcal{S}$ may be decomposed in $p$ sets where $p\in\{2,\,3,\,4,\,5\}$:\\
\noindent --- If $m\geq 3$, then:\\
\noindent $\triangleright$ If $m\equiv 0 \mod 3$, set
$S_{l}:=\{s_{3k+l,i_2}\ /\ (k,\,i_2)\in\llbr
0,\, \frac{m}{3}-1\rrbr\times\llbr
0,\, q+2\rrbr\}$ for $l\in\llbr
0,\, 2\rrbr$. Then $\mathcal{S}=S_0\sqcup S_1\sqcup
S_2,$ $p=3$, and $C_\Gamma=3I-\tau_0-\tau_1-\tau_2$.\\
\noindent $\triangleright$ If $m\equiv 1 \mod 3$, set
$S_{l}:=\{s_{3k+l,i_2}\ /\ (k,\,i_2)\in\llbr
0,\, \frac{m-1}{3}-1\rrbr\times\llbr
0,\, q+2\rrbr\}$ for $l\in\llbr
0,\, 2\rrbr$, and $S_{3}:=\{s_{m-1,i_2}\ /\ i_2\in\llbr
0,\, q+2\rrbr\}$. Then $\mathcal{S}=S_0\sqcup S_1\sqcup
S_2\sqcup S_3,$ $p=4$, and $C_\Gamma=4I-\tau_0-\tau_1-\tau_2-\tau_3$.\\
\noindent $\triangleright$ If $m\equiv 2 \mod 3$, set
$S_{l}:=\{s_{3k+l,i_2}\ /\ (k,\,i_2)\in\llbr
0,\, \frac{m-2}{3}-1\rrbr\times\llbr
0,\, q+2\rrbr\}$ for $l\in\llbr
0,\, 2\rrbr$, and $S_{3}:=\{s_{m-2,i_2}\ /\ i_2\in\llbr
0,\, q+2\rrbr\}$, $S_{4}:=\{s_{m-1,i_2}\ /\ i_2\in\llbr
0,\, q+2\rrbr\}$. Then $\mathcal{S}=S_0\sqcup S_1\sqcup
S_2\sqcup S_3\sqcup S_4,$ $p=5$, and $C_\Gamma=5I-\tau_0-\tau_1-\tau_2-\tau_3-\tau_4$.\\

\noindent --- If $m=1$, then:\\
\noindent $\triangleright$ If $q=2$, set
$S_{0}:=\{s_{0,0},\,s_{0,1},\,s_{0,2}\,s_{0,3}\}$ and $S_{1}:=\{s_{0,4}\}$.
Then $\mathcal{S}=S_0\sqcup S_1$, $p=2$, and $C_\Gamma=2I-\tau_0-\tau_1$.\\
\noindent $\triangleright$ If $q\geq 3$ and $q$ is even, set
$S_{0}:=\{s_{0,0},\,s_{0,1},\,s_{0,2}\}$ and $S_{1}:=\{s_{0,3},\dots,\,s_{0,\,q+1}\}$, $S_{2}:=\{s_{0,4},\dots,\,s_{0,\,q+2}\}$.
Then $\mathcal{S}=S_0\sqcup S_1 \sqcup S_2$, $p=3$, and $C_\Gamma=2I-\tau_0-\tau_1-\tau_1$.\\
\noindent $\triangleright$ If $q\geq 3$ and $q$ is odd, set
$S_{0}:=\{s_{0,0},\,s_{0,1},\,s_{0,2}\}$ and $S_{1}:=\{s_{0,3},\dots,\,s_{0,\,q+2}\}$, $S_{2}:=\{s_{0,4},\dots,\,s_{0,\,q+1}\}$.
Then $\mathcal{S}=S_0\sqcup S_1 \sqcup S_2$, $p=3$, and $C_\Gamma=2I-\tau_0-\tau_1-\tau_1$.\\

\noindent $\bullet$ If $m=1$, the graph associated to $\Gamma$ is the following:\\

\begin{scriptsize}
\def\arrow(#1,#2){\ncline{-}{#1}{#2}}
$$\begin{array}{c@{\hskip .5cm}c@{\hskip .2cm}c@{\hskip .2cm}c@{\hskip .5cm}c}
      \rnode{1}{\Tcircle{1}} & & & & \rnode{0}{\Tcircle{0}}\\[.1cm]
      & \rnode{q2}{\Tcircle{q+2}} & \rnode{pts}{\cdots}  & \rnode{4}{\Tcircle{4}} \\[.1cm]
      \rnode{2}{\Tcircle{2}} & & & & \rnode{3}{\Tcircle{3}} \\[.1cm]
    \end{array}
    \everypsbox{\scriptstyle}
    \psset{nodesep=0pt,arm=.6,linearc=.4,angleA=0,angleB=90}
     \arrow(1,q2) \arrow(2,q2)    \arrow(4,0) \arrow(4,3)
$$
\end{scriptsize}

\noindent If $m\geq 2$, the graph associated to $\Gamma$ consists in $q+3$ $m-$gones that are linked
together.\\

\subsubsection{\textsf{The $BTa$ subseries --- Binary tetrahedral groups}}

\noindent Let $\psi_4,\ \tau,\ \eta,\ \phi_{2m}$ be the elements $$\psi_4:=\left(
                                                                             \begin{array}{ccc}
                                                                               1 & 0 & 0 \\
                                                                               0 & \zeta_4 & 0 \\
                                                                               0 & 0 & \zeta_4^{-1} \\
                                                                             \end{array}
                                                                           \right),\
\tau:=\left(
        \begin{array}{ccc}
          1 & 0 & 0 \\
          0 & 0 & i \\
          0 & i & 0 \\
        \end{array}
      \right),\
\eta:=\frac{1}{\sqrt{2}}\left(
        \begin{array}{ccc}
          \sqrt{2} & 0 & 0 \\
          0 & \zeta_8^7 & \zeta_8^7 \\
          0 & \zeta_8^5 & \zeta_8 \\
        \end{array}
      \right),\
\phi_{2m}:=\left(
             \begin{array}{ccc}
               \zeta_{2m}^{-2} & 0 & 0 \\
               0 & \zeta_{2m} & 0 \\
               0 & 0 & \zeta_{2m} \\
             \end{array}
           \right).
$$
In this section\footnote{The other case --- the type $BTb$ --- is $m\equiv 3 \mod 6$. This group is not a direct product.}, we assume that $m\equiv 1\, \textrm{or}\, 5 \mod 6$, and we consider the subgroup
$\Gamma:=\langle\psi_{4},\ \tau,\ \eta,\ \phi_{2m}\rangle$ of $\mathbf{SL}_3\mathbb{C}$. Note that $\phi_{2m}=\psi_4^2\phi_m^{-\frac{m-1}{2}}$, so that $\Gamma:=\langle\psi_{4},\ \tau,\ \eta,\ \phi_{m}\rangle$. Set $\Gamma_1:=\langle\psi_{4},\ \tau,\ \eta\rangle$ and $\Gamma_2:=\langle\phi_{m}\rangle\simeq \mathbb{Z}/m\mathbb{Z}$. Then we have $\Gamma\simeq \Gamma_2\times\Gamma_1$.\\
With the notations of the binary tetrahedral subgroup of $\mathbf{SL}_2\mathbb{C}$, $\psi_4$ (resp. $\tau$, $\eta$) represents $a^2$ (resp. $b$, $c$). $\Gamma_1\simeq\langle a^2,\,b,\,c\rangle$ is the binary tetrahedral subgroup of $\mathbf{SL}_2\mathbb{C}$. Representatives to its $7$ conjugacy classes are
$$\{id,\ a^4=-id,\ b,\ c,\ c^2,\ -c,\ -c^2\},$$
and its character table is the matrix given in Section 2.3.1 (Type $E_6$ --- Binary tetrahedral group).

\noindent The natural character of $\Gamma$ is given by $$\begin{array}{lcrcl|lcrcl}
                                                                \chi_{i,0} & = & \chi(\phi_m^i) & = & \zeta_m^{-2i}+2\zeta_m^{i} & \chi_{i,4} & = & \chi(\phi_m^i\eta^2) & = & \zeta_m^{-2i}-\zeta_m^{i} \\
                                                                \chi_{i,1} & = & \chi(\phi_m^i\psi_4^2) & = & \zeta_m^{-2i}-2\zeta_m^{i} & \chi_{i,5} & = & \chi(\phi_m^i\psi_4^2\eta) & = & \zeta_m^{-2i}-\zeta_m^{i} \\
                                                                \chi_{i,2} & = & \chi(\phi_m^i\tau) & = & \zeta_m^{-2i} & \chi_{i,6} & = & \chi(\phi_m^i\psi_4^2\eta^2) & = & \zeta_m^{-2i}+\zeta_m^{i} \\
                                                                \chi_{i,3} & = & \chi(\phi_m^i\eta) & = & \zeta_m^{-2i}+\zeta_m^{i} &   &   &   \end{array}
$$
Finally, we obtain the series $P_\Gamma(t,\,u)$: for $p\in\llbr 0,\,m-1\rrbr$, we get

$$P_\Gamma(t,\,u)_{7p} = \displaystyle{\frac{1}{24\,m}\sum_{k=0}^{m-1}\zeta_m^{pk}(f_{{k,0}}+f_{{k,1}}+6\,f_{{k,2}}+4\,f_{{k,3}}+4\,f_{{k,4}}+4\,f_{{k,5}}+4\,f_{{k,6}}),}$$
$$P_\Gamma(t,\,u)_{7p+1} = \displaystyle{\frac{1}{24\,m}\sum_{k=0}^{m-1}\zeta_m^{pk}(f_{{k,0}}+
f_{{k,1}}+6\,f_{{k,2}}+ \left( 4\,f_{{k,3}}+4\,f_{{k,5}} \right) j+
 \left( 4\,f_{{k,4}}+4\,f_{{k,6}} \right) {j}^{2}),}$$
$$P_\Gamma(t,\,u)_{7p+3}  =  \displaystyle{\frac{1}{24\,m}\sum_{k=0}^{m-1}\zeta_m^{pk}(2\,f_{{k,0}}-2\,f_{{k,1}}+4\,f_{{k,3}}-4\,f_{{k,4}
}-4\,f_{{k,5}}+4\,f_{{k,6}}),}$$
$$P_\Gamma(t,\,u)_{7p+4}  =  \displaystyle{\frac{1}{24\,m}\sum_{k=0}^{m-1}\zeta_m^{pk}(2\,f_{{k,0}}-2\,f_{{k,1
}}+ \left( 4\,f_{{k,3}}-4\,f_{{k,5}} \right) j+ \left( -4\,f_{{k,4}}+4
\,f_{{k,6}} \right) {j}^{2}),}$$
$$P_\Gamma(t,\,u)_{7p+6}  =  \displaystyle{\frac{1}{24\,m}\sum_{k=0}^{m-1}\zeta_m^{pk}(3\,f_{{k,0}}+3\,f_{{k,1
}}-6\,f_{{k,2}}),}$$
and $P_\Gamma(t,\,u)_{7p+2}$ (resp. $P_\Gamma(t,\,u)_{7p+5}$) is obtained by exchanging $j$ and $j^2$ in $P_\Gamma(t,\,u)_{7p+1}$ (resp. $P_\Gamma(t,\,u)_{7p+4}$).\\

\subsubsection{\textsf{The $BO$ subseries --- Binary octahedral groups}}

\noindent For $m\in \mathbb{N}$ such that $m\wedge 6=1$, let $\psi_8,\ \tau,\ \eta,\ \phi_{2m}$ be the elements $$\psi_8:=\left(
                                                                             \begin{array}{ccc}
                                                                               1 & 0 & 0 \\
                                                                               0 & \zeta_8 & 0 \\
                                                                               0 & 0 & \zeta_8^{7} \\
                                                                             \end{array}
                                                                           \right),\
\tau:=\left(
        \begin{array}{ccc}
          1 & 0 & 0 \\
          0 & 0 & i \\
          0 & i & 0 \\
        \end{array}
      \right),\
\eta:=\frac{1}{\sqrt{2}}\left(
        \begin{array}{ccc}
          \sqrt{2} & 0 & 0 \\
          0 & \zeta_8^7 & \zeta_8^7 \\
          0 & \zeta_8^5 & \zeta_8 \\
        \end{array}
      \right),\
\phi_{2m}:=\left(
             \begin{array}{ccc}
               \zeta_{2m}^{-2} & 0 & 0 \\
               0 & \zeta_{2m} & 0 \\
               0 & 0 & \zeta_{2m} \\
             \end{array}
           \right),
$$
and consider the subgroup $\Gamma=\langle\psi_8,\ \tau,\ \eta,\ \phi_{2m}\rangle$ of $\mathbf{SL}_3\mathbb{C}$.
Note that $\phi_{2m}=\psi_8^4\phi_m^{-\frac{m-1}{2}}$, so that $\Gamma:=\langle\psi_{8},\ \tau,\ \eta,\ \phi_{m}\rangle$.
Set $\Gamma_1:=\langle\psi_{8},\ \tau,\ \eta\rangle$ and $\Gamma_2:=\langle\phi_{m}\rangle\simeq \mathbb{Z}/m\mathbb{Z}$.\\
Then $\Gamma\simeq\Gamma_2\times\Gamma_1$.\\
With the notations of the binary octahedral subgroup of $\mathbf{SL}_2\mathbb{C}$, $\psi_8$ (resp. $\tau$, $\eta$) represents $a$ (resp. $b$, $c$). $\Gamma_1\simeq\langle a,\,b,\,c\rangle$ is the binary octahedral subgroup of $\mathbf{SL}_2\mathbb{C}$. Reperesentatives of its $8$ conjugacy classes are
$$\{id,\ a^4=-id,\ ab,\ b,\ c^2,\ c,\ a,\ a^3\},$$
and its character table is the matrix given in Section 2.3.2 (Type $E_7$ --- Binary octahedral group).\\
The natural character of $\Gamma$ is given by $$\begin{array}{lcrcl|lcrcl}
                                                                \chi_{i,0} & = & \chi(\phi_m^i) & = & \zeta_m^{-2i}+2\zeta_m^{i} & \chi_{i,4} & = & \chi(\phi_m^i\eta^2) & = & \zeta_m^{-2i}-\zeta_m^{i} \\
                                                                \chi_{i,1} & = & \chi(\phi_m^i\psi_8^4) & = & \zeta_m^{-2i}-2\zeta_m^{i} & \chi_{i,5} & = & \chi(\phi_m^i\eta) & = & \zeta_m^{-2i}+\zeta_m^{i} \\
                                                                \chi_{i,2} & = & \chi(\phi_m^i\psi_8\tau) & = & \zeta_m^{-2i} & \chi_{i,6} & = & \chi(\phi_m^i\psi_8) & = & \zeta_m^{-2i}+\zeta_m^{i}\sqrt{2} \\
                                                                \chi_{i,3} & = & \chi(\phi_m^i\tau) & = & \zeta_m^{-2i} & \chi_{i,7} & = & \chi(\phi_m^i\psi_8^3)  & =  &  \zeta_m^{-2i}-\zeta_m^{i}\sqrt{2} \end{array}
$$
\noindent Finally, we obtain the series $P_\Gamma(t,\,u)$: for $p\in\llbr 0,\,m-1\rrbr$, we have

$$P_\Gamma(t,\,u)_{8p}  =  \displaystyle{\frac{1}{48\,m}\sum_{k=0}^{m-1}\zeta_m^{pk}(f_{{k,0}}+f_{{k,1}}+12\,f_{{k,2}}+6\,f_{{k,3}}+8\,f_{{k,4}}+8\,f_{{k,5}}+6\,f_{{k,6}}+6\,f_{{k,7}}),}$$
$$P_\Gamma(t,\,u)_{8p+2}  =  \displaystyle{\frac{1}{48\,m}\sum_{k=0}^{m-1}\zeta_m^{pk}(2
\,f_{{k,0}}+2\,f_{{k,1}}+12\,f_{{k,3}}-8\,f_{{k,4}}-8\,f_{{k,5}}),}$$
$$P_\Gamma(t,\,u)_{8p+3}  =  \displaystyle{\frac{1}{48\,m}\sum_{k=0}^{m-1}\zeta_m^{pk}(2\,f_{{k,0}}-2\,f_{{k,1}}-8\,f_{{k,4}}+8\,f_{{k,5}
}+6\,\sqrt {2}f_{{k,6}}-6\,\sqrt {2}f_{{k,7}}),}$$
$$P_\Gamma(t,\,u)_{8p+5}  =  \displaystyle{\frac{1}{48\,m}\sum_{k=0}^{m-1}\zeta_m^{pk}(3\,f_{{k,0}}+3\,f_{{k,1}}-12\,
f_{{k,2}}-6\,f_{{k,3}}+6\,f_{{k,6}}+6\,f_{{k,7}}),}$$
$$P_\Gamma(t,\,u)_{8p+7}  =  \displaystyle{\frac{1}{48\,m}\sum_{k=0}^{m-1}\zeta_m^{pk}(4\,f_{{k,0}}-4\,f_{{k,1}}+8\,f_{{k,4}}-8\,f
_{{k,5}}),}$$
and $P_\Gamma(t,\,u)_{8p+1}$ (resp. $P_\Gamma(t,\,u)_{8p+4}$, $P_\Gamma(t,\,u)_{8p+6}$) is obtained by replacing $f_{k,2},\ f_{k,6},\ f_{k,7}$
by their opposite in $P_\Gamma(t,\,u)_{8p}$ (resp. $P_\Gamma(t,\,u)_{8p+3}$, $P_\Gamma(t,\,u)_{8p+5}$).\\

\subsubsection{\textsf{The $BI$ subseries --- Binary icosahedral groups}}

\noindent For $m\in \mathbb{N}$ such that $m\wedge 30=1$, let $\mu,\ \tau,\ \eta$ be the elements
$$ \mu:=\left(
        \begin{array}{ccc}
          1 & 0 & 0 \\
          0 & -\zeta_5^3 & 0 \\
          0 & 0 & -\zeta_5^2 \\
        \end{array}
      \right),\
\tau:=\left(
\begin{array}{ccc}
1 & 0 & 0 \\
0 & 0 & 1 \\
0 & -1 & 0 \\
\end{array}
\right),\
\eta:=\frac{1}{\zeta_5^2-\zeta_5^{-2}}\left(
        \begin{array}{ccc}
          \zeta_5^2-\zeta_5^{-2} & 0 & 0 \\
          0 & \zeta_5-\zeta_5^{-1} & 1 \\
          0 & 1 & -\zeta_5-\zeta_5^{-1} \\
        \end{array}
      \right),$$
and consider the subgroup $\Gamma=\langle\mu,\ \tau,\ \eta,\ \phi_{2m}\rangle$ of $\mathbf{SL}_3\mathbb{C}$.\\
Note that $\phi_{2m}=\eta^2\phi_m^{-\frac{m-1}{2}}$, so that $\Gamma:=\langle\mu,\ \tau,\ \eta,\ \phi_{m}\rangle$.
Set $\Gamma_1:=\langle\mu,\ \tau,\ \eta \rangle$ and $\Gamma_2:=\langle\phi_{m}\rangle\simeq \mathbb{Z}/m\mathbb{Z}$.
Then $\Gamma\simeq \Gamma_2\times\Gamma_1$.\\
With the notations of the binary icosahedral subgroup of $\mathbf{SL}_2\mathbb{C}$, $\mu$ (resp. $\tau$, $\eta$) represents $a$
(resp. $b$, $c$). $\Gamma_1\simeq\langle a,\,b,\,c\rangle$ is the binary icosahedral subgroup of $\mathbf{SL}_2\mathbb{C}$. Representatives of its $9$ conjugacy classes are
$$\{id,\ b^2=-id,\ a,\ a^2,\ a^3,\ a^4,\ abc,\ (abc)^2,\ b\},$$
and its character table is given in Section 2.3.3 (Type $E_8$ --- Binary icosahedral group).\\
The natural character of $\Gamma$ is given by $$\begin{array}{lcrcl|lcrcl}
                                                                \chi_{i,0} & = & \chi(\phi_m^i) & = & \zeta_m^{-2i}+2\zeta_m^{i} & \chi_{i,5} & = & \chi(\phi_m^i\mu^4) & = & \zeta_m^{-2i}-\frac{1+\sqrt{5}}{2}\,\zeta_m^{i} \\
                                                                \chi_{i,1} & = & \chi(\phi_m^i\tau^2) & = & \zeta_m^{-2i}-2\zeta_m^{i} & \chi_{i,6} & = & \chi(\phi_m^i\mu\tau\eta) & = & \zeta_m^{-2i}+\zeta_m^{i} \\
                                                                \chi_{i,2} & = & \chi(\phi_m^i\mu) & = & \zeta_m^{-2i}+\frac{1+\sqrt{5}}{2}\,\zeta_m^{i} & \chi_{i,7} & = & \chi(\phi_m^i(\mu\tau\eta)^2) & = & \zeta_m^{-2i}-\zeta_m^{i} \\
                                                                \chi_{i,3} & = & \chi(\phi_m^i\mu^2) & = & \zeta_m^{-2i}+\frac{-1+\sqrt{5}}{2}\,\zeta_m^{i} & \chi_{i,8} & = & \chi(\phi_m^i\tau)  & =  &  \zeta_m^{-2i}  \\
                                                                \chi_{i,4} & = & \chi(\phi_m^i\mu^3) & = & \zeta_m^{-2i}+\frac{1-\sqrt{5}}{2}\,\zeta_m^{i} &  &  &
                                                                \end{array}
$$
\noindent Finally, we obtain the series $P_\Gamma(t,\,u)$: for $p\in\llbr 0,\,m-1\rrbr$, we have

$$\begin{array}{rcl}
P_\Gamma(t,\,u)_{9p} & = & \displaystyle{\frac{1}{120\,m}\sum_{k=0}^{m-1}\zeta_m^{pk}
(f_{{k,0}}+f_{{k,1}}+12\,f_{{k,2}}+12\,f_{{k,3}}+12\,f_{{k,4}}+12\,f_{{k,5}}+20\,f_{{k,6}}+20\,f_{{k,7}}+30\,f_{{k,
8}}),}
\end{array}$$
$$\begin{array}{rcl}
P_\Gamma(t,\,u)_{9p+1} & = & \displaystyle{\frac{1}{120\,m}\sum_{k=0}^{m-1}\zeta_m^{pk}(2\,f_{{k,0}}-2\,f_{{k,1}}+ ( 6-6\,\sqrt {5
} ) f_{{k,2}}+ ( -6-6\,\sqrt {5} ) f_{{k,3}}+(
6+6\,\sqrt {5} ) f_{{k,4}} }\\
 & & \ \ \ \ \ \ \ \ \ \ \ \ \ \ \ \displaystyle{+ ( -6+6\,\sqrt {5} ) f_{{k
,5}}+20\,f_{{k,6}}-20\,f_{{k,7}}),}
  \end{array}
$$
$$\begin{array}{rcl}
P_\Gamma(t,\,u)_{9p+3} & = & \displaystyle{\frac{1}{120\,m}\sum_{k=0}^{m-1}\zeta_m^{pk}(3\,f_{{k,0}}+3\,f_{{k,1}}+ ( 6+6\,\sqrt {5}
 ) f_{{k,2}}+ ( 6-6\,\sqrt {5} ) f_{{k,3}}+ ( 6-
6\,\sqrt {5} ) f_{{k,4}} }\\
 & & \ \ \ \ \ \ \ \ \ \ \ \ \ \ \ \displaystyle{+ ( 6+6\,\sqrt {5} ) f_{{k,5}
}-30\,f_{{k,8}}),}
  \end{array}
$$
$$\begin{array}{rcl}
P_\Gamma(t,\,u)_{9p+5} & = & \displaystyle{\frac{1}{120\,m}\sum_{k=0}^{m-1}\zeta_m^{pk}(4\,f_{{k,0}}+4\,f_
{{k,1}}-12\,f_{{k,2}}-12\,f_{{k,3}}-12\,f_{{k,4}}-12\,f_{{k,5}}+20\,f_
{{k,6}}+20\,f_{{k,7}}) }
\end{array}
$$
$$\begin{array}{rcl}
P_\Gamma(t,\,u)_{9p+7} & = & \displaystyle{\frac{1}{120\,m}\sum_{k=0}^{m-1}\zeta_m^{pk}(5\,f_{{k,0}}+5\,f_{{k,1}}-20\,f_{{k,6}}
-20\,f_{{k,7}}+30\,f_{{k,8}}) }
\end{array}
$$
$$\begin{array}{rcl}
P_\Gamma(t,\,u)_{9p+8} & = & \displaystyle{\frac{1}{120\,m}\sum_{k=0}^{m-1}\zeta_m^{pk}(6\,f_{{k,0}}-6\,f_{{k,
1}}-12\,f_{{k,2}}+12\,f_{{k,3}}-12\,f_{{k,4}}+12\,f_{{k,5}}), }
\end{array}
$$
and $P_\Gamma(t,\,u)_{9p+4}$ (resp. $P_\Gamma(t,\,u)_{9p+2}$) is obtained by replacing $6\sqrt{5}$ by its opposite in $P_\Gamma(t,\,u)_{9p+3}$ (resp. $P_\Gamma(t,\,u)_{9p+1}$), and $P_\Gamma(t,\,u)_{9p+6}$ is obtained by replacing $f_{k,1},\ f_{k,2},\ f_{k,4},\ f_{k,6}$ by their opposite in $P_\Gamma(t,\,u)_{9p+5}$.\\

\subsubsection{\textsf{Decomposition of $C_\Gamma$ for the subseries $BTa$, $BO$, $BI$}}\label{decomposition}

\noindent $\bullet$ For the subseries $BTa$ (resp. $BO$, $BI$), we set $n=7$ (resp. $n=8$, $n=9$). We now make the matrix $A^{(1)}$ explicit:
 $A^{(1)}$ is a block-matrix with $m\times m$ blocks of size $n\times n$.\\
$\triangleright$ If $m\geq 5$, then the matrices $A^{(1)}$ and $C_\Gamma:=2I-A^{(1)}-\,^tA^{(1)}$ are defined by
\begin{small}
$$A^{(1)}=\left(
    \begin{array}{ccccc}
       0 &   & I   &   & B_n \\
      B_n & \ddots  &   & \ddots  &   \\
        & \ddots  &  \ddots  &    & I  \\
      I &   & \ddots  &  \ddots  &   \\
       & I &   & B_n & 0  \\
    \end{array}
  \right),\ \ C_\Gamma=\left(
    \begin{array}{ccccccc}
       2I & -B_n & -I  &   &   & -I & -B_n \\
      -B_n & 2I  & -B_n & -I  &   &   & -I \\
       -I & -B_n & 2I  & -B_n &  -I &   &   \\
        &  -I & -B_n & \ddots  & \ddots &  \ddots &   \\
        &   &  -I & \ddots & \ddots  & -B_n & -I  \\
      -I &   &   & \ddots  & -B_n &  2I & -B_n \\
      -B_n & -I &   &   & -I  & -B_n &  2I \\
    \end{array}
  \right),
$$
\end{small}
with
\begin{small}
$$B_7=\left(
            \begin{array}{ccccccc}
              0 & 0 & 0 & 1 & 0 & 0 & 0 \\
              0 & 0 & 0 & 0 & 1 & 0 & 0 \\
              0 & 0 & 0 & 0 & 0 & 1 & 0 \\
              1 & 0 & 0 & 0 & 0 & 0 & 1 \\
              0 & 1 & 0 & 0 & 0 & 0 & 1 \\
              0 & 0 & 1 & 0 & 0 & 0 & 1 \\
              0 & 0 & 0 & 1 & 1 & 1 & 0 \\
            \end{array}
          \right),\ B_8=\left( \begin {array}{cccccccc} 0&0&0&1&0&0&0&0\\\noalign{\medskip}0&0&0&0&1&0&0&0\\\noalign{\medskip}0&0&0&0&0&0&0&1\\\noalign{\medskip}1&0
&0&0&0&1&0&0\\\noalign{\medskip}0&1&0&0&0&0&1&0\\\noalign{\medskip}0&0
&0&1&0&0&0&1\\\noalign{\medskip}0&0&0&0&1&0&0&1\\\noalign{\medskip}0&0
&1&0&0&1&1&0\end {array} \right),\ B_9=\left( \begin {array}{ccccccccc} 0&0&1&0&0&0&0&0&0\\\noalign{\medskip}0&0&0&0&0&1&0&0&0\\\noalign{\medskip}1&0&0&1&0&0&0
&0&0\\\noalign{\medskip}0&0&1&0&0&0&1&0&0\\\noalign{\medskip}0&0&0&0&0
&0&0&0&1\\\noalign{\medskip}0&1&0&0&0&0&0&0&1\\\noalign{\medskip}0&0&0
&1&0&0&0&1&0\\\noalign{\medskip}0&0&0&0&0&0&1&0&1\\\noalign{\medskip}0
&0&0&0&1&1&0&1&0\end {array} \right).$$
\end{small}

\noindent $\triangleright$ If $m=1$, then the matrix $A^{(1)}$ is $A^{(1)}:=B_n+I$,
and $C_\Gamma$ is defined by  $C_\Gamma:=2I-A^{(1)}-\,^tA^{(1)}+2\,Diag(A^{(1)})$, i.e. $C_\Gamma=2I-2B_n$.\\

\noindent $\bullet$ The decomposition of $C_\Gamma$ is the following: for $i_1\in\llbr0,\,m-1\rrbr$ and
$i_2\in\llbr 0,\,n-1\rrbr$, let $s_{i_1,i_2}$ be the reflection
associated to the $(ni_1+i_2)-$th root. Then the set
$\mathcal{S}$ may be decomposed in $p$ sets where $p\in\{2,\,4,\,5\}$.\\

\noindent --- If $m\geq 5$, then:\\
\noindent $\triangleright$ If $m\equiv 1 \mod 3$, set
$S_{l}:=\{s_{3k+l,i_2}\ /\ (k,\,i_2)\in\llbr
0,\, \frac{m-1}{3}-1\rrbr\times\llbr
0,\, n-1\rrbr\}$ for $l\in\llbr
0,\, 2\rrbr$, and $S_{3}:=\{s_{m-1,i_2}\ /\ i_2\in\llbr
0,\, n-1\rrbr\}$. Then $S=S_0\sqcup S_1\sqcup
S_2\sqcup S_3,$ $p=4$, and $C_\Gamma=4I-\tau_0-\tau_1-\tau_2-\tau_3$.\\
\noindent $\triangleright$ If $m\equiv 2 \mod 3$, set
$S_{l}:=\{s_{3k+l,i_2}\ /\ (k,\,i_2)\in\llbr
0,\, \frac{m-2}{3}-1\rrbr\times\llbr
0,\, n-1\rrbr\}$ for $l\in\llbr
0,\, 2\rrbr$, and $S_{3}:=\{s_{m-2,i_2}\ /\ i_2\in\llbr
0,\, n-1\rrbr\}$, $S_{4}:=\{s_{m-1,i_2}\ /\ i_2\in\llbr
0,\, n-1\rrbr\}$. Then $S=S_0\sqcup S_1\sqcup
S_2\sqcup S_3\sqcup S_4,$ $p=5$, and $C_\Gamma=5I-\tau_0-\tau_1-\tau_2-\tau_3-\tau_4$.\\

\noindent --- If $m=1$, then $S=S_0\sqcup S_1$, $p=2$, and $C_\Gamma=2I-\tau_0-\tau_1$, with\\
$\triangleright$ If $n=7$,
$S_{0}:=\{s_{0,0},\,s_{0,1},\,s_{0,2},\,s_{0,6}\}$ and $S_{1}:=\{s_{0,3},\,s_{0,4},\,s_{0,5}\}$.\\
$\triangleright$ If $n=8$,
$S_{0}:=\{s_{0,0},\,s_{0,1},\,s_{0,2},\,s_{0,5},\,s_{0,6}\}$ and $S_{1}:=\{s_{0,3},\,s_{0,4},\,s_{0,7}\}$.\\
$\triangleright$ If $n=9$,
$S_{0}:=\{s_{0,0},\,s_{0,3},\,s_{0,4},\,s_{0,5},\,s_{0,7}\}$ and $S_{1}:=\{s_{0,1},\,s_{0,2},\,s_{0,6},\,s_{0,8}\}$.\\

\noindent $\bullet$ If $m=1$, the graph associated to $\Gamma$ is the following:
\def\arrow(#1,#2){\ncline{-}{#1}{#2}}
\begin{scriptsize}
$$\begin{array}{||c|c|c||}
\hline\hline & &   \\[-.2cm]
  n=7 & n=8 & n=9\\ \hline
  \begin{array}{c@{\hskip .3cm}c@{\hskip .3cm}c@{\hskip .3cm}c@{\hskip .3cm}c}
      \rnode{1}{\Tcircle{1}} & \rnode{4}{\Tcircle{4}}\\[.2cm]
      & & \rnode{6}{\Tcircle{6}} & \rnode{3}{\Tcircle{3}}  & \rnode{0}{\Tcircle{0}} \\[.2cm]
      \rnode{2}{\Tcircle{2}} & \rnode{5}{\Tcircle{5}} \\[.2cm]
    \end{array}
    \everypsbox{\scriptstyle}
    \psset{nodesep=0pt,arm=.6,linearc=.4,angleA=0,angleB=90}
     \arrow(1,4) \arrow(4,6)   \arrow(2,5) \arrow(5,6) \arrow(6,3) \arrow(3,0) &
\begin{array}{c@{\hskip .3cm}c@{\hskip .3cm}c@{\hskip .3cm}c@{\hskip .3cm}c}
      \rnode{0}{\Tcircle{0}} & \rnode{3}{\Tcircle{3}} & \rnode{5}{\Tcircle{5}}  \\[.2cm]
      & & & \rnode{7}{\Tcircle{7}} & \rnode{2}{\Tcircle{2}}  \\[.2cm]
      \rnode{1}{\Tcircle{1}} & \rnode{4}{\Tcircle{4}} & \rnode{6}{\Tcircle{6}} \\[.2cm]
    \end{array}
    \everypsbox{\scriptstyle}
    \psset{nodesep=0pt,arm=.6,linearc=.4,angleA=0,angleB=90}
     \arrow(0,3) \arrow(3,5)   \arrow(5,7) \arrow(7,2) \arrow(1,4) \arrow(4,6) \arrow(6,7) &
\begin{array}{c@{\hskip .3cm}c@{\hskip .3cm}c@{\hskip .3cm}c@{\hskip .3cm}c@{\hskip .3cm}c@{\hskip .3cm}c@{\hskip .3cm}c}
      \rnode{1}{\Tcircle{1}} & \rnode{5}{\Tcircle{5}} \\[.2cm]
      & & \rnode{8}{\Tcircle{8}} & \rnode{7}{\Tcircle{7}} & \rnode{6}{\Tcircle{6}} & \rnode{3}{\Tcircle{3}} & \rnode{2}{\Tcircle{2}} & \rnode{0}{\Tcircle{0}}  \\[.2cm]
      & \rnode{4}{\Tcircle{4}}  \\[.2cm]
    \end{array}
    \everypsbox{\scriptstyle}
    \psset{nodesep=0pt,arm=.6,linearc=.4,angleA=0,angleB=90}
     \arrow(1,5) \arrow(5,8)   \arrow(4,8) \arrow(8,7) \arrow(7,6) \arrow(6,3) \arrow(3,2) \arrow(2,0) \\
  \hline \hline
\end{array}
$$
\end{scriptsize}

\noindent If $m\geq 5$, the graph associated to $\Gamma$ is a graph of type $BTa$ (resp. $BO$, $BI$), $m=1$, such that
every vertex is a $m-$gone.\\


\subsection{\textsf{The $C$ series}}

\noindent Let
$H\simeq\mathbb{Z}/j_1\mathbb{Z}\times\mathbb{Z}/j_2\mathbb{Z}$
be a group of the series $A$, with eventually $j_1=1$ or $j_2=1$,
and consider the matrix $T:=\left(%
\begin{array}{ccc}
  0 & 1 & 0 \\
  0 & 0 & 1 \\
  1 & 0 & 0 \\
\end{array}%
\right)$, that is the matrix of the permutation $(1,\,2,\,3)$ of
$\go{S}_3$. In this section, we study~$\Gamma:=\langle H,\,T\rangle$, the finite
subgroup of $\mathbf{SL}_3\mathbb{C}$ generated by $H$ and $T$. The subgroup $N$ of $\Gamma$ which consists of all the diagonal
matrices of $\Gamma$ is a normal subgroup of $\Gamma$. By using the Bezout theorem,
\begin{equation}\label{descrN}
N=\left\{g_{k_1,k_2}:=\left(%
\begin{array}{ccc}
  \zeta_m^{k_1} & 0 & 0 \\
  0 & \zeta_m^{k_2} & 0 \\
  0 & 0 & \zeta_m^{-k_1-k_2} \\
\end{array}%
\right)\ /\ (k_1,\,k_2)\in\llbr 0,\,m-1\rrbr^2\right\}.
\end{equation}
Moreover, we have $N\cap\langle T\rangle=\{id\}$ and
$|N\langle T\rangle|=\frac{|N||\langle T\rangle|}{|N\cap\langle
T\rangle|}=3m^2=|\Gamma|$. So, $\Gamma$ is the semi-direct product
$$\Gamma\simeq N\rtimes \langle T\rangle\simeq (\mathbb{Z}/m\mathbb{Z})^2 \rtimes \langle T\rangle.$$
We will obtain all the irreducible characters
of $\Gamma$ by induction; we distinguish
two cases corresponding to the two following subsections.\\

\subsubsection{\textsf{Series $C$ --- $m$ non divisible by $3$}}

\begin{figure}[h]
\begin{center}
\begin{tabular}{ccc}
\includegraphics[width=6cm,height=6cm]{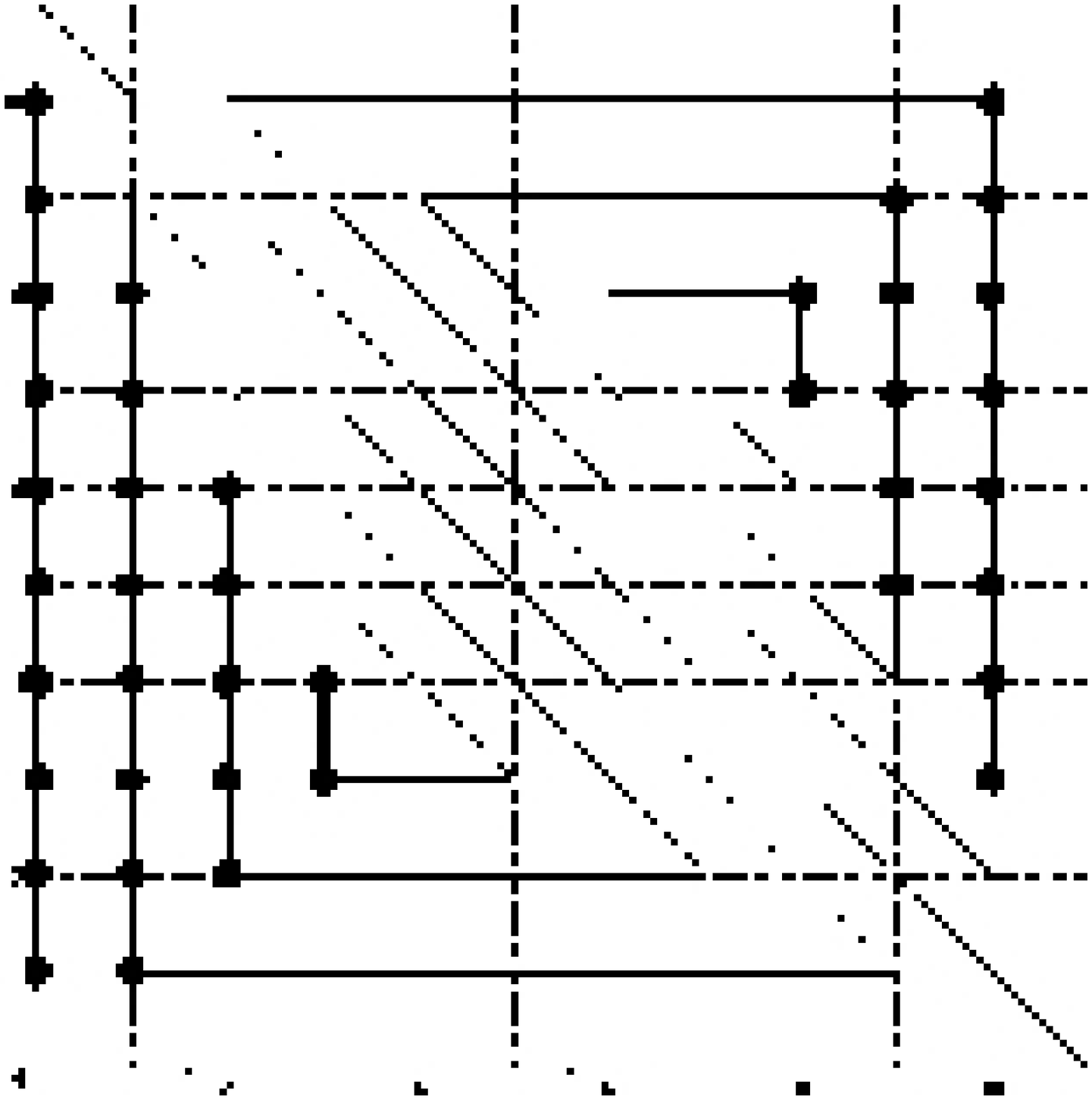} & &  \includegraphics[width=6cm,height=6cm]{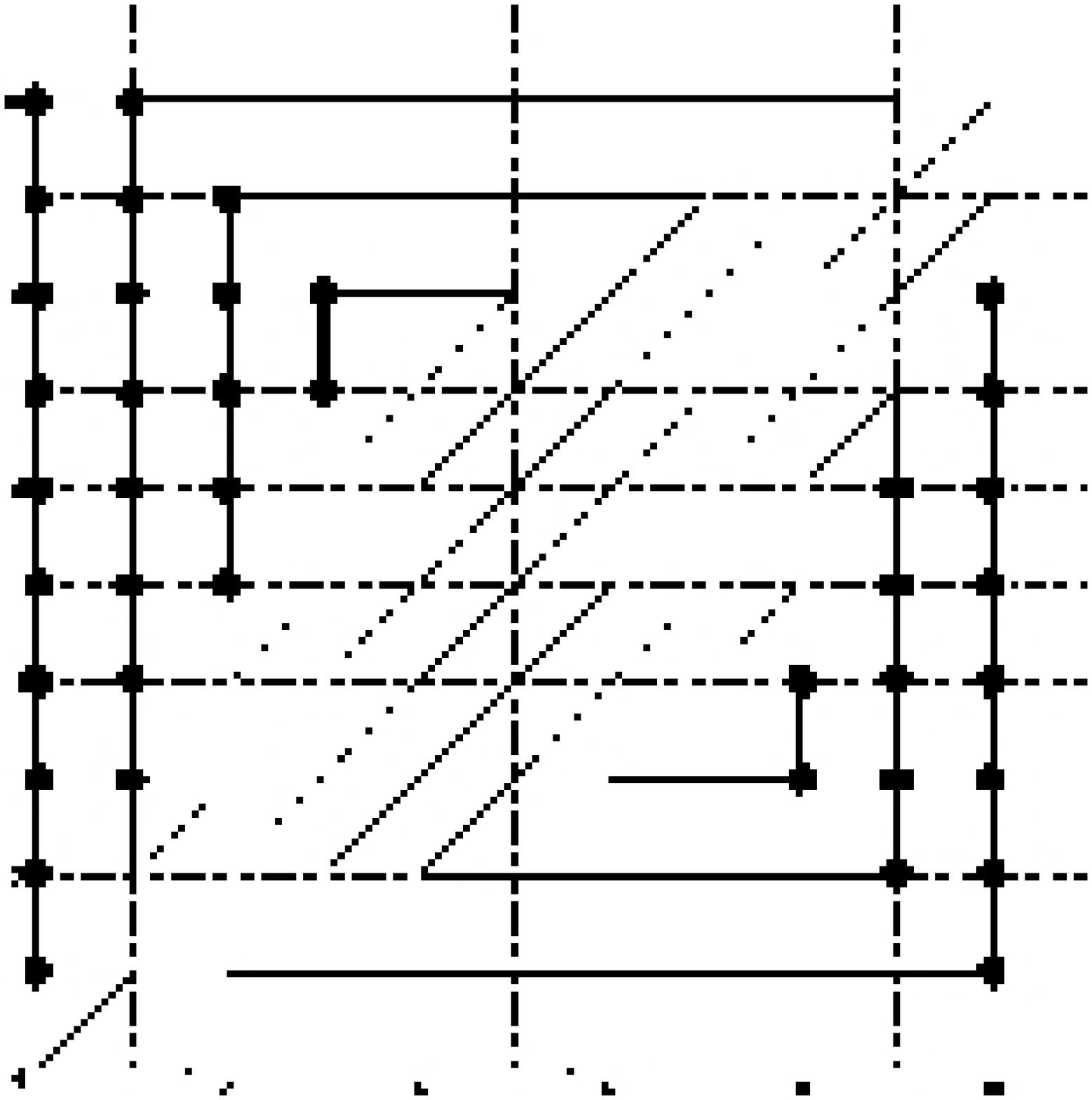} \\
 \end{tabular}
\end{center}
\caption{Conjugacy classes and irreducible characters for $m=11$.}\label{conj&car11}
\end{figure}

\noindent $\bullet$ Set $n':=\frac{m^2-1}{3}$, so that $|N|=3n'+1$ and
$|G|=3m^2=3(3n'+1)$. The conjugacy classes of $\Gamma$ are:
$$\begin{array}{||l||c|c|c|c||}
\hline\hline & & & & \\[-.3cm]
  \textrm{Class} & id & T & T^{-1} & g\in N\backslash\{id\} (n'\ \textrm{classes}) \\
\hline
  \textrm{Cardinality} & 1 & m^2 & m^2 & 3  \\
\hline\hline
\end{array}$$
For each element $g_{k_1,k_2}\in N\backslash\{id\}$, the conjugacy
class of $g$ is the set $\left\{g_{k_1,k_2},\,g_{k_2,-k_1-k_2},\,g_{-k_1-k_2,k_1}\right\}$. In order to obtain a transversal of
$N\backslash\{id\}$, i.e. a set
containing exactly one representant of each conjugacy class of
$N\backslash\{id\}$, we represent the elements of
$N\backslash\{id\}$ by points
$(k_1,\,k_2)$ of $\llbr 0,\,m-1\rrbr^2$.\\
So, we search a
transversal for the set of elements of the form $(k_1,\,k_2),\
(k_2,\,-k_1,-k_2 \mod m)$ and $(-k_1,-k_2 \mod m,\,k_1)$, with
$(k_1,\,k_2)\in\llbr 0,\,m-1\rrbr^2$.\\
A solution is the following: for a given conjugacy class, its
three elements are on the edges of a triangle (see Figure
\ref{conj&car11}), with exactly one element on each edge of the
triangle. Therefore we may take as transversal the set of all points that
belong to the vertical edges minus the nearest point of the
diagonal. More precisely, a transversal for $N\backslash\{id\}$ is
the set $E_{cc}$ defined by
$$\begin{array}{rl}
   & \{(0,\,k_2)\ /\ k_2\in\llbr 1,\,m-1\rrbr\}\vspace{.2cm} \\
  \sqcup & \{(k_1,\,k_2)\ /\
k_1\in\llbr 1,\,\left\lfloor\frac{m}{3}\right\rfloor\rrbr,\
k_2\in\llbr k_1,\,m-1-2k_1\rrbr\}\vspace{.2cm} \\
  \sqcup & \{(k_1,\,k_2)\ /\
k_1\in\llbr m-\lfloor\frac{m}{3}\rfloor,\,m-1\rrbr,\ k_2\in\llbr
2(m-k_1)+1,\,k_1\rrbr\}. \\
\end{array}$$

\noindent We choose the usual lexicographic order on $E_{cc}$, so
that we can number its elements:
$E_{cc}=\{c_1,\,c_2,\dots,\,c_{\frac{m^2-1}{3}}\},$ with
$c_i=(c_i^{(1)},\,c_i^{(2)})\in\llbr 0,\,m-1\rrbr^2.$\\

\noindent $\bullet$ The group $\Gamma$ is generated by $R:=g_{1,0}$ and $T$, which verify the relations
$R^m=(RT)^3=T^3=id.$ As $m$ is not divisible by $3$, the irreducible characters of degree $1$
are $\chi^{0,l}:R \mapsto 1, T \mapsto j^l$ for $l\in\llbr 0,\,2\rrbr$.
We have $[G\,:\,N]=3$ with $N$ abelian, so the possible
degrees of the irreducible characters are $1,\,2,\,3$. The irreducible characters $\chi_{l_1,\,l_2}$
induced by the irreducible characters of $N$ are given by
$$\begin{array}{||l||c|c|c|c||} \hline\hline & & & & \\[-.4cm]
  \textrm{Class} & [id] & [T] & [T^{-1}] & [g],\ g\in N\backslash\{id\}  \\ \hline  & & & & \\[-.3cm]
  \textrm{Value} & 3 & 0 & 0 & \zeta_m^{k_1l_1+k_2l_2}+\zeta_m^{(-k_1-k_2)l_1+k_1l_2}+\zeta_m^{k_2l_1+(-k_1-k_2)l_2} \\ \hline\hline
\end{array}$$
The characters $\chi_{l_1,\,l_2}$ with $(l_1,\,l_2)\neq(0,\,0)$ are represented
by points $(l_1,\,l_2)$ of~$\llbr 0,\,m-1\rrbr^2$. The points that
are associated to the same character are on a triangle or on a
``trident \,\Symbole'', with exactly one point on each edge (see Figure
\ref{conj&car11}). So the set of irreducible characters $\chi_{l_1,\,l_2}$ with
$(l_1,\,l_2)\neq(0,\,0)$ is obtained by taking the following set $E_{ic}$
of indexes:
$$\begin{array}{rl}
   & \{(0,\,k_2)\ /\ k_2\in\llbr 1,\,m-1\rrbr\}\vspace{.2cm} \\
  \sqcup & \{(k_1,\,k_2)\ /\
k_1\in\llbr 1,\,\left\lfloor\frac{m}{3}\right\rfloor\rrbr,\
k_2\in\llbr 2k_1+1,\,m-k_1\rrbr\}\vspace{.2cm} \\
  \sqcup & \{(k_1,\,k_2)\ /\
k_1\in\llbr m-\lfloor\frac{m}{3}\rfloor,\,m-1\rrbr,\ k_2\in\llbr
m-k_1,\,2k_1-m-1\rrbr\}. \\
\end{array}$$

\noindent We choose the usual lexicographic order on $E_{ic}$, so
that we can number its elements:
$E_{ic}=\{d_1,\,d_2,\dots,\,d_{\frac{m^2-1}{3}}\},$ with
$d_i=(d_i^{(1)},\,d_i^{(2)})\in\llbr 0,\,m-1\rrbr^2,$ and the character table $T_\Gamma$ is
$$T_\Gamma=\left(
  \begin{array}{ccc|ccc}
    1 & 1 & 1 & 1 & \cdots & 1 \\
    1 & j & j^2 & 1 & \cdots & 1 \\
    1 & j^2 & j & 1 & \cdots & 1  \\ \hline
    3 & 0 & 0 &  &  &  \\
    \vdots & \vdots & \vdots &  & C &  \\
    3 & 0 & 0 &  &  &  \\
  \end{array}
\right),$$
where the general term of $C\in\mathbf{M}_{\frac{m^2-1}{3}}\mathbb{C}$ is
$$c_{i,\,j}:=\zeta_m^{c_j^{(1)}d_i^{(1)}+c_j^{(2)}d_i^{(2)}}+\zeta_m^{(-c_j^{(1)}-c_j^{(2)})d_i^{(1)}+c_j^{(1)}d_i^{(2)}}+
\zeta_m^{c_j^{(2)}d_i^{(1)}+(-c_j^{(1)}-c_j^{(2)})d_i^{(2)}},\ (i,\,j)\in\llbr 1,\,\frac{m^2-1}{3}\rrbr^2.$$

\noindent $\bullet$ The values of the natural character $\chi$ of
$\Gamma$ are $$\begin{array}{||l||c|c|c|c||} \hline\hline & & & & \\[-.4cm]
  \textrm{Class} & [id] & [T] & [T^{-1}]& [g],\ g\in N\backslash\{id\}  \\ \hline & & & & \\[-.3cm]
  \textrm{Value} & 3 & 0 & 0 & \zeta_m^{k_1}+\zeta_m^{k_2}+\zeta_m^{-k_1-k_2}  \\ \hline\hline
\end{array}$$
Therefore the diagonal matrix $\Delta(t,\,u)$ is $Diag\left(\varepsilon_1,\,\varepsilon_2,\,\varepsilon_3,\,\widetilde{\Delta(t,\,u)}\right)$,
with $\varepsilon_1:=f(3,\,3),\ \varepsilon_2=\varepsilon_3:=f(0,\ 0)$, and the general term of $\widetilde{\Delta(t,\,u)}\in \mathbf{M}_{\frac{m^2-1}{3}}\mathbb{C}$ is $$\gamma_j:=f\left(\zeta_m^{-c_j^{(1)}}+\zeta_m^{-c_j^{(2)}}+\zeta_m^{c_j^{(1)}+c_j^{(2)}},\,\zeta_m^{c_j^{(1)}}+
\zeta_m^{c_j^{(2)}}+\zeta_m^{-c_j^{(1)}-c_j^{(2)}}\right),\ j\in\llbr 1,\,\frac{m^2-1}{3}\rrbr.$$
Then, by setting $$\Sigma:=\sum_{p=1}^{\frac{m^2-1}{3}}\gamma_p,\ \Sigma_{i}:=\sum_{j=1}^{\frac{m^2-1}{3}}\gamma_j\left(\zeta_m^{c_j^{(1)}d_i^{(1)}+c_j^{(2)}d_i^{(2)}}
+\zeta_m^{(-c_j^{(1)}-c_j^{(2)})d_i^{(1)}+c_j^{(1)}d_i^{(2)}}
+\zeta_m^{c_j^{(2)}d_i^{(1)}+(-c_j^{(1)}-c_j^{(2)})d_i^{(2)}}\right),$$ we obtain the formula for $P_\Gamma(t,\,u)$:

$$\begin{array}{rcl}
    P_\Gamma(t,\,u)_0 & = & \displaystyle{\frac{5(m+1)(m-1)+1}{9m^4}(3\Sigma+\varepsilon_1+\varepsilon_2+\varepsilon_3)-\frac{2(m+1)(m-1)}{9m^2}(\varepsilon_1-\varepsilon_2-\varepsilon_3)} \\
     &   & \displaystyle{+\frac{2(m+1)(m-1)}{3m^4}\sum_{q=1}^{\frac{m^2-1}{3}}\left(\varepsilon_1+\varepsilon_2+\varepsilon_3+\sum_{p=1}^{\frac{m^2-1}{3}}\gamma_pc_{p,\,q}\right),}
  \end{array}$$
$$\begin{array}{rcl}
    P_\Gamma(t,\,u)_1 & = & \displaystyle{\frac{5(m+1)(m-1)+1}{9m^4}(3\Sigma+\varepsilon_1+j\varepsilon_2+j^2\varepsilon_3)-\frac{2(m+1)(m-1)}{9m^2}(\varepsilon_1-j\varepsilon_2-j^2\varepsilon_3)} \\
     &   & \displaystyle{+\frac{2(m+1)(m-1)}{3m^4}\sum_{q=1}^{\frac{m^2-1}{3}}\left(\varepsilon_1+j\varepsilon_2+j^2\varepsilon_3+\sum_{p=1}^{\frac{m^2-1}{3}}\gamma_pc_{p,\,q}\right),}
  \end{array}$$
and $P_\Gamma(t,\,u)_2$ is obtained by exchanging $j$ and $j^2$ in $P_\Gamma(t,\,u)_1$, and for $i\in\llbr 1,\,\frac{m^2-1}{3}\rrbr$,
$$\begin{array}{rcl}
   P_\Gamma(t,\,u)_{i+2} & = & \displaystyle{\frac{5(m+1)(m-1)+1}{9m^4}(3\Sigma_{i}+3\varepsilon_1)-\frac{6(m+1)(m-1)}{9m^2}\varepsilon_1} \\
     &   & \displaystyle{+\frac{2(m+1)(m-1)}{3m^4}\sum_{q=1}^{\frac{m^2-1}{3}}\left(3\varepsilon_1+\sum_{p=1}^{\frac{m^2-1}{3}}\gamma_pc_{i,\,p}c_{p,\,q}\right).}
  \end{array}$$

\subsubsection{\textsf{Series $C$ --- $m$ divisible by $3$}}

\noindent $\bullet$ Set $n':=\frac{m^2}{3}$, so that $|N|=3n'$ and
$|G|=3m^2=3(3n')$. Set $a:=R^{\frac{m}{3}}=Diag(j,\,j,\,j)$. The conjugacy classes of $\Gamma$ are:
$$\begin{array}{||l||c|c|c|c|c|c|c|c|c|c||}
\hline\hline & & & & & & & & & & \\[-.4cm]
  \textrm{Class} & id  & a & a^2 &  T &
  T^{-1} & RT & RT^{-1} & R^2T & R^2T^{-1} & g\in N\backslash\{id,\,a,\,a^2\} (n'-1\ \textrm{classes})\\
\hline & & & & & & & & & & \\[-.4cm]
  \textrm{Cardinality} & 1  & 1 & 1 & \frac{m^2}{3} & \frac{m^2}{3} & \frac{m^2}{3} & \frac{m^2}{3} & \frac{m^2}{3} & \frac{m^2}{3} & 3 \\[.1cm]
\hline\hline
\end{array}$$
For each element $g_{k_1,k_2}\in N\backslash\{id,\,a,\,a^2\}$, the conjugacy class of $g$ is the set
$\left\{g_{k_1,k_2},\,g_{k_2,-k_1-k_2},\,g_{-k_1-k_2,k_1}\right\}$. Its
three elements are on the edges of a triangle, with exactly one element of each edge of the
triangle. So, a transversal of $N\backslash\{id,\,a,\,a^2\}$ has the same form as in the case where $3$ does not divide $m$, i.e.
 a transversal for $N\backslash\{id,\,a,\,a^2\}$ is the set $E_{cc}$ defined by
$$\begin{array}{rl}
   & \{(0,\,k_2)\ /\ k_2\in\llbr 1,\,m-1\rrbr\}\vspace{.2cm} \\
  \sqcup & \{(k_1,\,k_2)\ /\
k_1\in\llbr 1,\,\frac{m}{3}-1\rrbr,\
k_2\in\llbr k_1,\,m-1-2k_1\rrbr\}\vspace{.2cm} \\
  \sqcup & \{(k_1,\,k_2)\ /\
k_1\in\llbr m-\frac{m}{3}+1,\,m-1\rrbr,\ k_2\in\llbr
2(m-k_1)+1,\,k_1\rrbr\}. \\
\end{array}$$

\noindent $\bullet$ As $m$ is divisible by $3$, the irreducible characters
of degree $1$ of $\Gamma$ are, for $(k,\,l)\in\llbr 0,\,2\rrbr^2$, the nine elements
$$\begin{array}{rcl}
  \chi^{k,l}:R & \mapsto & j^k \\
  T & \mapsto & j^l. \\
\end{array}$$
As for the case where $m$ is divisible by $3$, the set of irreducible characters $\chi_{l_1,\,l_2}$ with
$(l_1,\,l_2)\neq(0,\,0)$ is obtained by taking the following set $E_{ic}$
of indexes:
$$\begin{array}{rl}
   & \{(0,\,k_2)\ /\ k_2\in\llbr 1,\,m-1\rrbr\}\vspace{.2cm} \\
  \sqcup & \{(k_1,\,k_2)\ /\
k_1\in\llbr 1,\,\frac{m}{3}-1\rrbr,\
k_2\in\llbr 2k_1+1,\,m-k_1\rrbr\}\vspace{.2cm} \\
  \sqcup & \{(k_1,\,k_2)\ /\
k_1\in\llbr m-\frac{m}{3}+1,\,m-1\rrbr,\ k_2\in\llbr
m-k_1,\,2k_1-m-1\rrbr\}. \\
\end{array}$$
We choose the usual lexicographic order on $E_{cc}$ and $E_{ic}$, so
that we can number its elements:
$$E_{cc}=\{c_1,\,c_2,\dots,\,c_{\frac{m^2}{3}-1}\},\ \ E_{ic}=\{d_1,\,d_2,\dots,\,d_{\frac{m^2}{3}-1}\},$$ with
$c_j=(c_j^{(1)},\,c_j^{(2)})\in\llbr 0,\,m-1\rrbr^2$, and $d_i=(d_i^{(1)},\,d_i^{(2)})\in\llbr 0,\,m-1\rrbr^2$, and we deduce the character table~$T_\Gamma$ of
$\Gamma$:
$$T_\Gamma=\left(
             \begin{array}{ccc|ccc|ccc|c}
               1 & 1 & 1 & 1 & 1 & 1 & 1 & 1 & 1 &   \\
               1 & 1 & 1 & j & j^2 & j & j^2 & j & j^2 &  1 \\
               1 & 1 & 1 & j^2 & j & j^2 & j & j^2 & j &   \\ \hline
               1 & 1 & 1 & 1 & 1 & j & j & j^2 & j^2 &   \\
               1 & 1 & 1 & j & j^2 & j^2 & 1 & 1 & j & j^{c_i^{(1)}-c_i^{(2)}} \\
               1 & 1 & 1 & j^2 & j & 1 & j^2 & j & 1 &   \\ \hline
               1 & 1 & 1 & 1 & 1 & j^2 & j^2 & j & j &   \\
               1 & 1 & 1 & j & j^2 & 1 & j & j^2 & 1 &  j^{2(c_i^{(1)}-c_i^{(2)})} \\
               1 & 1 & 1 & j^2 & j & j & 1 & 1 & j^2 &   \\ \hline
               3 & J_1^{(1)} & J_1^{(2)} & 0 & 0 & 0 & 0 & 0 & 0 &   \\
               \vdots & \vdots & \vdots & \vdots & \vdots & \vdots & \vdots & \vdots & \vdots &  C \\
               3 & J_{\frac{m^2}{3}-1}^{(1)} & J_{\frac{m^2}{3}-1}^{(2)} & 0 & 0 & 0 & 0 & 0 & 0 &   \\
             \end{array}
           \right),$$
where $C\in \mathbf{M}_{\frac{m^2}{3}-1}\mathbb{C}$ is a block-matrix with general term
$$c_{i,\,j}:=\zeta_m^{c_j^{(1)}d_i^{(1)}+c_j^{(2)}d_i^{(2)}}+\zeta_m^{(-c_j^{(1)}-c_j^{(2)})d_i^{(1)}+c_j^{(1)}d_i^{(2)}}+
\zeta_m^{c_j^{(2)}d_i^{(1)}+(-c_j^{(1)}-c_j^{(2)})d_i^{(2)}},\ (i,\,j)\in\llbr 1,\,\frac{m^2}{3}-1\rrbr^2,$$
$$J_i^{(1)}:=j^{d_i^{(1)}+d_i^{(2)}}+j^{-2d_i^{(1)}+d_i^{(2)}}+j^{d_i^{(1)}-2d_i^{(2)}},\ \
J_i^{(2)}:=j^{2d_i^{(1)}+2d_i^{(2)}}+j^{-d_i^{(1)}+2d_i^{(2)}}+j^{2d_i^{(1)}-d_i^{(2)}},\ i\in\llbr 1,\,\frac{m^2}{3}-1\rrbr^2.$$

\noindent $\bullet$ The values of the natural character $\chi$ of
$\Gamma$ are $$\begin{array}{||l||c|c|c|c|c|c|c|c|c|c||}
\hline\hline & & & & & & & & & & \\[-.4cm]
  \textrm{Class} & id  & a & a^2 &  T &
  T^{-1} & RT & RT^{-1} & R^2T & R^2T^{-1} & g\in N\backslash\{id,\,a,\,a^2\}\\
\hline & & & & & & & & & & \\[-.3cm]
  \textrm{Value} & 3  & 3j & 3j^2 & 0 & 0 & 0 & 0 & 0 & 0 & \zeta_m^{k_1}+\zeta_m^{k_2}+\zeta_m^{-k_1-k_2} \\
\hline\hline
\end{array}$$

\noindent Therefore the diagonal matrix $\Delta(t,\,u)$ is $Diag(\widetilde{\widetilde{\Delta(t,\,u)}},\,\widetilde{\Delta(t,\,u)})$,
with $$\widetilde{\widetilde{\Delta(t,\,u)}}=Diag(\beta_1,\,\beta_2,\,\beta_3,\,\underbrace{\beta_4,\dots,\,\beta_4}_{6\ \textrm{terms}})=
Diag\left(f(3,\,3),\,f(3j^2,\,3j),\,f(3j,\,3j^2),\,\underbrace{f(0,\,0),\dots,\,f(0,\,0)}_{6\ \textrm{terms}}\right),$$ and the general term of $\widetilde{\Delta(t,\,u)}\in \mathbf{M}_{\frac{m^2}{3}-1}\mathbb{C}$ is $$\gamma_j:=f\left(\zeta_m^{-c_j^{(1)}}+\zeta_m^{-c_j^{(2)}}+\zeta_m^{c_j^{(1)}+c_j^{(2)}},\,\zeta_m^{c_j^{(1)}}+
\zeta_m^{c_j^{(2)}}+\zeta_m^{-c_j^{(1)}-c_j^{(2)}}\right),\ j\in\llbr 1,\,\frac{m^2}{3}-1\rrbr.$$
For $(i,\,r)\in\{1,\,2\}\times\{1,\,2,\,3\}$, and $(s,\,q)\in\{1,\,2\}\times\llbr1,\,\frac{m^2}{3}-1\rrbr$, let us define $$\Sigma^{(r)}:=\sum_{p=1}^{\frac{m^2}{3}-1}\gamma_pj^{(r-1)(c_p^{(1)}-c_p^{(1)})},\ \Phi^{(i,r)}:=\sum_{p=1}^{\frac{m^2}{3}-1}J_p^{(i)}\gamma_pj^{(r-1)(c_p^{(1)}-c_p^{(1)})},$$
$$\Sigma^{(s)}_q:=\sum_{p=1}^{\frac{m^2}{3}-1}J_p^{(s)}\gamma_p\left(\zeta_m^{c_p^{(1)}d_q^{(1)}+c_p^{(2)}d_q^{(2)}}+
\zeta_m^{(-c_p^{(1)}-c_p^{(2)})d_q^{(1)}+c_p^{(1)}d_q^{(2)}}+
\zeta_m^{c_p^{(2)}d_q^{(1)}+(-c_p^{(1)}-c_p^{(2)})d_q^{(2)}}\right),$$
$$\xi:=(\underbrace{-1,\,-1,\,2}_{1},\,\underbrace{-1,\,-1,\,2}_{2},\,\underbrace{-1,\,-1,\,2}_{3},\,\dots,\,\underbrace{-1,\,-1,\,2}_
{\frac{\frac{m}{3}(\frac{m}{3}+1)}{2}-1},\,-1,\,-1,\,\underbrace{-1,\,-1,\,2},\,\underbrace{-1,\,-1,\,2},\,\dots,\,\underbrace{-1,\,-1,\,2}).$$
Then, we may give the expression of the series $P_\Gamma(t,\,u)$:

\begin{small}
$$\begin{array}{rcl}
    P_\Gamma(t,\,u)_0 & = & \displaystyle{\frac{3(5(\frac{m}{3})^2-1)+1}{3m^4}(3\Sigma^{(1)}+\beta_1+\beta_2+\beta_3+6\beta_4)
    +\frac{6(\frac{m}{3})^2-2}{3m^4}(\Phi^{(11)}+\Phi^{(21)}+2\beta_1+2\beta_2+2\beta_3+12\beta_4)} \\
      &   & \displaystyle{+\frac{3(\frac{m}{3})^2-1}{m^4}\sum_{q=1}^{\frac{m^2}{3}-1}\xi_q\left(\beta_1+\beta_2+\beta_3+3j^{c_q^{(1)}-c_q^{(2)}}\beta_4
      +3j^{2(c_q^{(1)}-c_q^{(2)})}\beta_4+\sum_{p=1}^{\frac{m^2}{3}-1}\gamma_pc_{p,q}\right),}
  \end{array}
$$
$$\begin{array}{rcl}
    P_\Gamma(t,\,u)_1 & = & \displaystyle{\frac{3(5(\frac{m}{3})^2-1)+1}{3m^4}(3\Sigma^{(1)}+\beta_1+\beta_2+\beta_3-3\beta_4)
    +\frac{6(\frac{m}{3})^2-2}{3m^4}(\Phi^{(11)}+\Phi^{(21)}+2\beta_1+2\beta_2+2\beta_3-6\beta_4)} \\
      &   & \displaystyle{+\frac{3(\frac{m}{3})^2-1}{m^4}\sum_{q=1}^{\frac{m^2}{3}-1}\xi_q\left(\beta_1+\beta_2+\beta_3+(2j+j^2)j^{c_q^{(1)}-c_q^{(2)}}\beta_4
      +(2j^2+j)j^{2(c_q^{(1)}-c_q^{(2)})}\beta_4+\sum_{p=1}^{\frac{m^2}{3}-1}\gamma_pc_{p,q}\right),}
  \end{array}
$$
$$\begin{array}{rcl}
    P_\Gamma(t,\,u)_3 & = & \displaystyle{\frac{3(5(\frac{m}{3})^2-1)+1}{3m^4}(3\Sigma^{(2)}+\beta_1+\beta_2+\beta_3)
    +\frac{6(\frac{m}{3})^2-2}{3m^4}(\Phi^{(12)}+\Phi^{(22)}+2\beta_1+2\beta_2+2\beta_3)} \\
      &   & \displaystyle{+\frac{3(\frac{m}{3})^2-1}{m^4}\sum_{q=1}^{\frac{m^2}{3}-1}\xi_q\left(\beta_1+\beta_2+\beta_3+(2+j)j^{c_q^{(1)}-c_q^{(2)}}\beta_4
      +(2j^2+j)j^{2(c_q^{(1)}-c_q^{(2)})}\beta_4+\sum_{p=1}^{\frac{m^2}{3}-1}\gamma_pj^{c_p^{(1)}-c_p^{(2)}}c_{p,q}\right),}
  \end{array}
$$
$$\begin{array}{rcl}
    P_\Gamma(t,\,u)_5 & = & \displaystyle{\frac{3(5(\frac{m}{3})^2-1)+1}{3m^4}(3\Sigma^{(2)}+\beta_1+\beta_2+\beta_3)
    +\frac{6(\frac{m}{3})^2-2}{3m^4}(\Phi^{(12)}+\Phi^{(22)}+2\beta_1+2\beta_2+2\beta_3)} \\
      &   & \displaystyle{+\frac{3(\frac{m}{3})^2-1}{m^4}\sum_{q=1}^{\frac{m^2}{3}-1}\xi_q\left(\beta_1+\beta_2+\beta_3
      +\sum_{p=1}^{\frac{m^2}{3}-1}\gamma_pj^{c_p^{(1)}-c_p^{(2)}}c_{p,q}\right),}
  \end{array}
$$
$$\begin{array}{rcl}
    P_\Gamma(t,\,u)_6 & = & \displaystyle{\frac{3(5(\frac{m}{3})^2-1)+1}{3m^4}(3\Sigma^{(3)}+\beta_1+\beta_2+\beta_3)
    +\frac{6(\frac{m}{3})^2-2}{3m^4}(\Phi^{(13)}+\Phi^{(23)}+2\beta_1+2\beta_2+2\beta_3)} \\
      &   & \displaystyle{+\frac{3(\frac{m}{3})^2-1}{m^4}\sum_{q=1}^{\frac{m^2}{3}-1}\xi_q\left(\beta_1+\beta_2+\beta_3+(2+j^2)j^{c_q^{(1)}-c_q^{(2)}}\beta_4
      +(2j+j^2)j^{2(c_q^{(1)}-c_q^{(2)})}\beta_4+\sum_{p=1}^{\frac{m^2}{3}-1}\gamma_pj^{2(c_p^{(1)}-c_p^{(2)})}c_{p,q}\right),}
  \end{array}
$$
\end{small}
\noindent and $P_\Gamma(t,\,u)_2$ (resp. $P_\Gamma(t,\,u)_4$, $P_\Gamma(t,\,u)_8$) is obtained by exchanging the coefficients $2j+j^2$ and $2j^2+j$ (resp. $2+j$ and $2j^2+j$, $2+j^2$ and $2j+j^2$) in $P_\Gamma(t,\,u)_1$ (resp. $P_\Gamma(t,\,u)_3$, $P_\Gamma(t,\,u)_6$); $P_\Gamma(t,\,u)_7$ is obtained by replacing $\Phi^{(1i)}$ by $\Phi^{(2i)}$ and $j^{c_p^{(1)}-c_p^{(2)}}$ by $j^{2(c_p^{(1)}-c_p^{(2)})}$ in $P_\Gamma(t,\,u)_5$.\\
Finally, for $i\in\llbr 1,\,\frac{m^2}{3}-1\rrbr$, we have
\begin{small}
$$\begin{array}{rcl}
    P_\Gamma(t,\,u)_{i+8} & = & \displaystyle{\frac{3(5(\frac{m}{3})^2-1)+1}{3m^4}(3\beta_1+J_i^{(1)}\beta_2+J_i^{(2)}\beta_3+3\Sigma_i)
    +\frac{6(\frac{m}{3})^2-2}{3m^4}(6\beta_1+2J_i^{(1)}\beta_2+2J_i^{(2)}\beta_3+\widetilde{\Sigma_i^{(1)}}+\widetilde{\Sigma_i^{(2)}})} \\
      &   & \displaystyle{+\frac{3(\frac{m}{3})^2-1}{m^4}\sum_{q=1}^{\frac{m^2}{3}-1}\xi_q\left(3\beta_1+J_i^{(1)}\beta_2+J_i^{(2)}\beta_3
      +\sum_{p=1}^{\frac{m^2}{3}-1}\gamma_pc_{i,p}c_{p,q}\right).}
  \end{array}
$$
\end{small}

\subsubsection{\textsf{Decomposition of $C_\Gamma$}}

\noindent We now make the matrix $A^{(1)}$ explicit: the form of the matrix $A^{(1)}$ is nearly the same in the case where
$m$ is divisible by $3$ as in the other case. The main difference between these two cases is due to the fact that in the case where $m$ is divisible by $3$, there are $9$ irreducible characters of degree $1$ instead of $3$.\\
Set $\kappa_m:=\frac{m}{3}-1$ if $3$ divides $m$, and $\kappa_m:=\left\lfloor\frac{m}{3}\right\rfloor$ otherwise. The matrix $A^{(1)}$ is a block-matrix with $\left(2+2\kappa_m\right)\times\left(2+2\kappa_m\right)$ blocks: for example, if $m=16$, the matrices $A^{(1)}$ and $C_\Gamma$ are matrices of size $88$.\\
\begin{figure}[h]
\begin{center}
\begin{tabular}{cc}
\includegraphics[width=8cm,height=8cm]{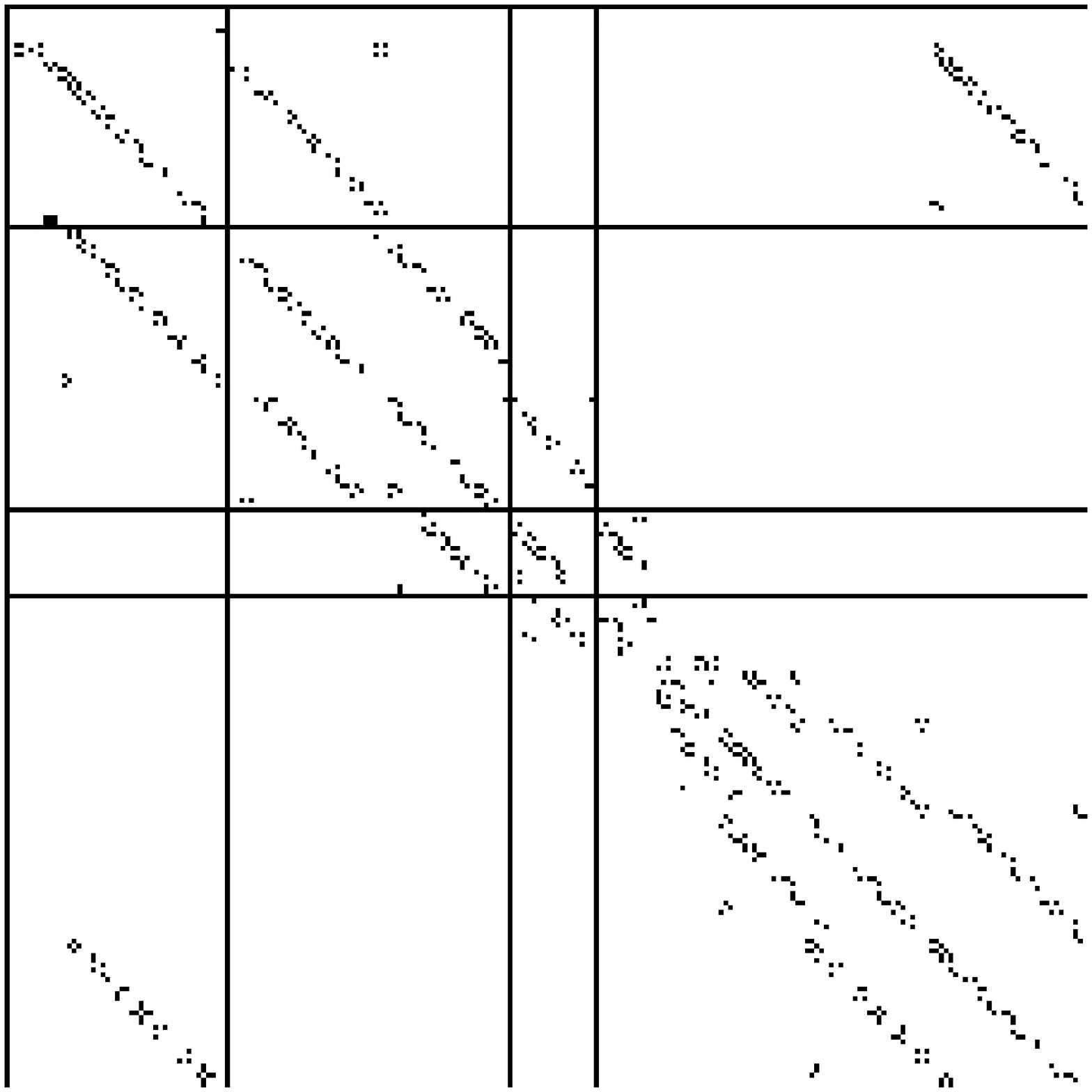}  & \includegraphics[width=8cm,height=8cm]{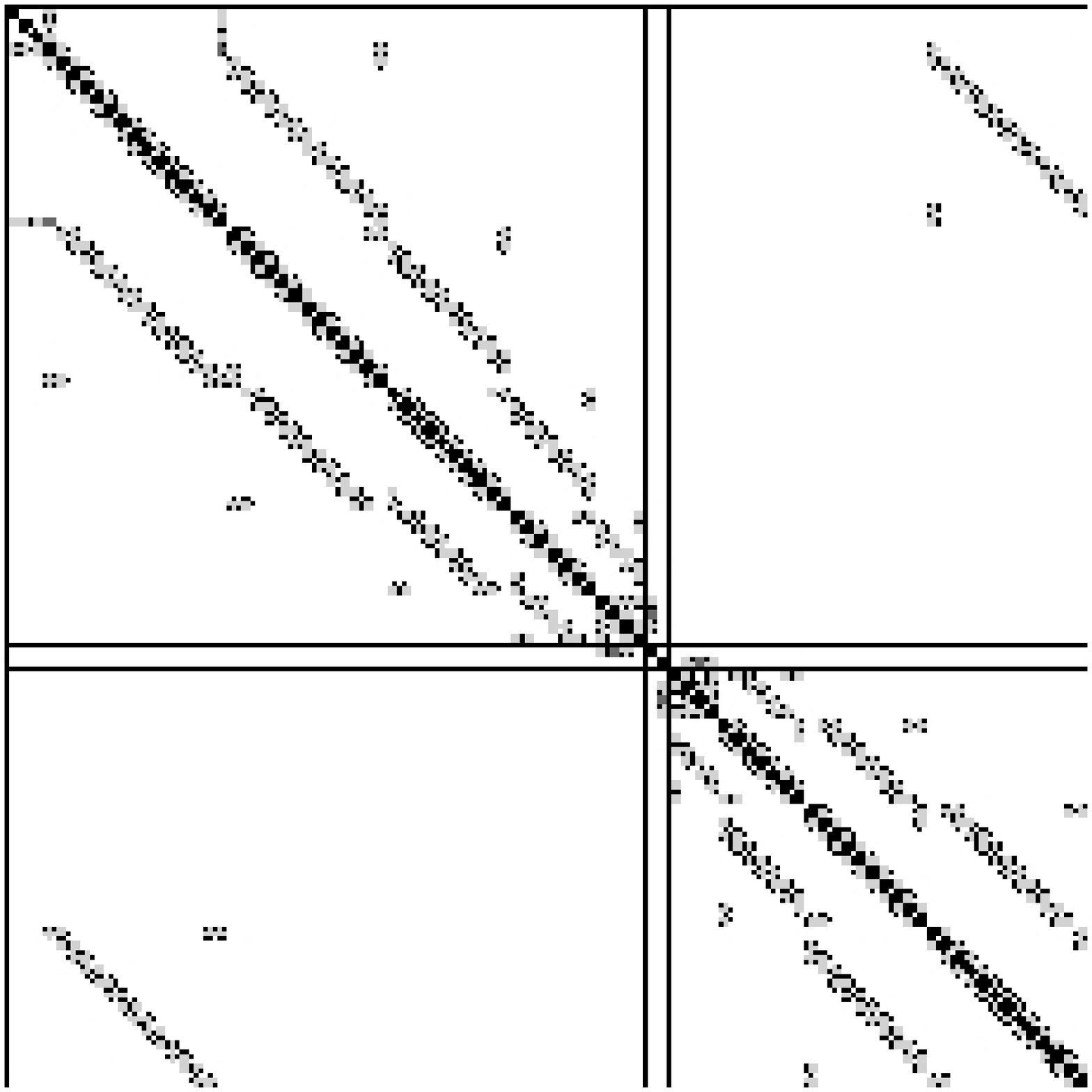} \\
 \end{tabular}
\end{center}
\caption{ The matrices $A^{(1)}$ and $C_\Gamma$ for $m=16$.}\label{mat16}
\end{figure}

\noindent $\bullet$ If $m=2$, then $A^{(1)}=\left(
                                                       \begin{array}{cc}
                                                         \mathbf{0}_{3,3} & \mathbf{1}_{3,1} \\
                                                         \mathbf{1}_{1,3} & 2 \\
                                                       \end{array}
                                                     \right)
$, and $C_\Gamma=3I-\tau_0-\tau_1-\tau_2$, with $\tau_0:=s_0s_2,\ \tau_1:=s_1,\ \tau_2:=s_3$.\\

\noindent $\bullet$ If $m=3$, then $A^{(1)}= \left(
                                              \begin{array}{cc}
                                                \mathbf{0}_{9,9} & \mathbf{1}_{9,2} \\
                                                \mathbf{1}_{2,9} & \widetilde{A} \\
                                              \end{array}
                                            \right),\ \widetilde{A}=\left(
                                                                                     \begin{array}{cc}
                                                                                       0 & 0 \\
                                                                                       3 & 0 \\
                                                                                     \end{array}
                                                                                   \right)
,$ and $C_\Gamma=2I-(s_0s_1\dots s_9)-(s_{10}s_{11})$.\\

\noindent $\bullet$ Now, we assume that $m\geq 4$. For $i_1\in \{0\}\sqcup\llbr 1,\,\kappa_m\rrbr\sqcup\llbr m-\kappa_m,\,m-1\rrbr$,
we define the set $\widetilde{S_{i_1}}$ by:\\
$$\widetilde{S_{i_1}}:=\left\{\begin{array}{ll}
                        \{s_{0,i_2}\ /\ i_2\in\llbr 1,\,m-1\rrbr\} & \textrm{if}\ i_1=0, \\
                        \{s_{i_1,i_2}\ /\ i_2\in\llbr i_1,\,m-1-2i_1\rrbr\} & \textrm{if}\ i_1\in\llbr 1,\,\kappa_m\rrbr, \\
                        \{s_{i_1,i_2}\ /\ i_2\in\llbr 2(m-i_1)+1,\,i_1\rrbr\} & \textrm{if}\ i_1\in\llbr m-\kappa_m,\,m-1\rrbr.
                      \end{array}\right.
$$
Then, we distinguish two cases:\\
$\triangleright$ If $\kappa_m$ is odd, we set
$\widehat{I_0}:=\{0,\,2,\,4,\dots,\,\kappa_m-1,\,
m-\kappa_m,\,m-\kappa_m+2,\dots,\,m-3\},$
\begin{center}
$\widehat{I_1}:=\{1,\,3,\,5,\dots,\,\kappa_m,\,
m-\kappa_m+1,\,m-\kappa_m+3,\dots,\,m-2\},$ $\widehat{I_2}:=\{m-1\}.$
\end{center}
$\triangleright$ If $\kappa_m$ is even, we set
$\widehat{I_0}:=\{0,\,2,\,4,\dots,\,\kappa_m,\,
m-\kappa_m+1,\,m-\kappa_m+3,\dots,\,m-3\},$
\begin{center}$\widehat{I_1}:=\{1,\,3,\,5,\dots,\,\kappa_m-1,\,
m-\kappa_m,\,m-\kappa_m+2,\dots,\,m-2\},$ $\widehat{I_2}:=\{m-1\}.$
\end{center}
\noindent Then, the roots associated to the reflections of distinct
$\widetilde{S_{i_1}}$'s for $i_1$ belonging to a same $\widehat{I_k}$ are orthogonal.\\
Now, we decompose each $\widetilde{S_{i_1}}$, i.e. $\widetilde{S_{i_1}}=\widehat{S_{i_1,0}}\sqcup\dots\sqcup\widehat{S_{i_1,q-1}}$,
such that $q\in\{1,\,2,\,3\}$ and for every $k\in\llbr 0,\,q-1\rrbr$, the roots associated to the reflections belonging to $\widehat{S_{i_1,k}}$ are orthogonal:\\
$\triangleright$ If $i_1=0$, then\begin{description}
                                   \item[$\diamond$] If $m-1$ is odd, then $\widetilde{S_{0}}=\widehat{S_{0,0}}\sqcup\widehat{S_{0,1}}\sqcup\widehat{S_{0,2}}$, with
$$\widehat{S_{0,0}}=\{s_{0,1},\,s_{0,3},\dots,\,s_{0,m-3}\},\
\widehat{S_{0,1}}=\{s_{0,2},\,s_{0,4},\dots,\,s_{0,m-2}\},\
\widehat{S_{0,2}}=\{s_{0,m-1}\}.$$
                                   \item[$\diamond$] If $m-1$ is even, then
                                   $\widetilde{S_{0}}=\widehat{S_{0,0}}\sqcup\widehat{S_{0,1}}$, with
$$\widehat{S_{0,0}}=\{s_{0,1},\,s_{0,3},\dots,\,s_{0,m-2}\},\
\widehat{S_{0,1}}=\{s_{0,2},\,s_{0,4},\dots,\,s_{0,m-1}\}.$$
                                 \end{description}
$\triangleright$ If $i_1\in\llbr 1,\,\kappa_m\rrbr$, then we have $\widetilde{S_{i_1}}=\widehat{S_{i_1,0}}\sqcup\widehat{S_{i_1,1}}\sqcup\widehat{S_{i_1,2}}$, with:
\begin{description}
                                   \item[$\diamond$] If $m-3i_1$ is odd, then $\widehat{S_{i_1,0}}=\{s_{i_1,i_1},\,s_{i_1,i_1+2},\dots,\,s_{i_1,m-2i_1-3}\}$,
$$\widehat{S_{i_1,1}}=\{s_{i_1,i_1+1},\,s_{i_1,i_1+3},\dots,\,s_{i_1,m-2i_1-2}\},\
\widehat{S_{i_1,2}}=\{s_{i_1,m-2i_1-1}\}.$$
                                   \item[$\diamond$] If $m-3i_1$ is even, then $\widehat{S_{i_1,0}}=\{s_{i_1,i_1},\,s_{i_1,i_1+2},\dots,\,s_{i_1,m-2i_1-4}\}$,
$$\widehat{S_{i_1,1}}=\{s_{i_1,i_1+1},\,s_{i_1,i_1+3},\dots,\,s_{i_1,m-2i_1-1}\},\
\widehat{S_{i_1,2}}=\{s_{i_1,m-2i_1-2}\}.$$
                                 \end{description}

\noindent $\triangleright$ If $i_1\in\llbr m-\kappa_m,\,m-1\rrbr$, then we have
$\widetilde{S_{i_1}}=\widehat{S_{i_1,0}}\sqcup\widehat{S_{i_1,1}}\sqcup\widehat{S_{i_1,2}}$, with:
\begin{description}
                                   \item[$\diamond$] If $3i_1-2m$ is odd, then $\widehat{S_{i_1,0}}=\{s_{i_1,2(m-i_1)+1},\,s_{i_1,2(m-i_1)+3},\dots,\,s_{i_1,i_1-2}\}$,
$$\widehat{S_{i_1,1}}=\{s_{i_1,2(m-i_1)+2},\,s_{i_1,2(m-i_1)+4},\dots,\,s_{i_1,i_1-1}\},\
\widehat{S_{i_1,2}}=\{s_{i_1,i_1}\}.$$
                                   \item[$\diamond$] If $3i_1-2m$ is even, then $\widehat{S_{i_1,0}}=\{s_{i_1,2(m-i_1)+1},\,s_{i_1,2(m-i_1)+3},\dots,\,s_{i_1,i_1-1}\}$,
$$\widehat{S_{i_1,1}}=\{s_{i_1,2(m-i_1)+2},\,s_{i_1,2(m-i_1)+4},\dots,\,s_{i_1,i_1-2}\},\
\widehat{S_{i_1,2}}=\{s_{i_1,i_1}\}.$$
                                 \end{description}

\noindent Note that some sets $\widehat{S_{i_1,k}}$ can be empty for $k\in\{1,\,2\}$.\\

\noindent Finally, we set $S_{k,l}:=\coprod_{i_1\in \widehat{I_k}}\widehat{S_{i_1,l}}$, for $(k,\,l)\in\{0,\,1,\,2\}^2\backslash\{(2,\,2)\}$,
and $S_{2,2}:=\left(\coprod_{i_1\in \widehat{I_2}}\widehat{S_{i_1,2}}\right)\sqcup\{s_{-1,0},\dots,\,s_{-1,r}\}$, with $r=8$ if $3$ divides $m$, and $r=2$ otherwise. We denote by $p\in\llbr1,\,9\rrbr$ the number of non-empty sets $S_{k,l}$, and by $\tau_{k,l}$ the commutative product of the reflections of $S_{k,l}$. Then, $C_\Gamma=p\,I-\sum_{(k,l)\in\{0,\,1,\,2\}^2}\tau_{k,l}$.\\

\begin{Ex}$\\$
For $m=16$, we have the following decomposition:
$$\begin{array}{l|l}
  \tau_{0,0}=(s_{0,1}s_{0,3}\dots s_{0,13})(s_{2,2}s_{2,4}\dots s_{2,8})(s_{4,4})(s_{11,11})(s_{13,7}s_{13,9}s_{13,11}) &
  \tau_{0,2}=(s_{0,15})(s_{2,10})(s_{4,6})(s_{13,13}) \\
  \tau_{0,1}=(s_{0,2}s_{0,4}\dots s_{0,14})(s_{2,3}s_{2,5}\dots s_{2,11})(s_{4,5}s_{4,7})(s_{13,8}s_{13,10}s_{13,12}) & \tau_{1,2}=(s_{1,13})(s_{3,9})(s_{12,12})(s_{14,14}) \\
  \tau_{1,0}=(s_{1,1}s_{1,3}\dots s_{1,11})(s_{3,3}s_{3,5}s_{3,7})(s_{5,5})(s_{12,9}s_{12,11})(s_{14,5}s_{14,7}\dots s_{14,13}) & \tau_{2,0}=(s_{15,3}s_{15,5}\dots s_{15,13}) \\
  \tau_{1,1}=(s_{1,2}s_{1,4}\dots s_{1,12})(s_{3,4}s_{3,6}s_{3,8})(s_{12,10})(s_{14,6}s_{14,8}\dots s_{14,12}) &
  \tau_{2,1}=(s_{15,4}s_{15,6}\dots s_{15,14}) \\
  \tau_{2,2}=(s_{15,15})(s_{-1,0}s_{-1,1}s_{-1,2}) &   \\
\end{array}
$$
\end{Ex}

\subsection{\textsf{The $D$ series}}

A group of type $D$ is generated by a group of type $C$ and a matrix $$Q:=\left(
                                                                                                                  \begin{array}{ccc}
                                                                                                                    a & 0 & 0 \\
                                                                                                                    0 & 0 & b \\
                                                                                                                    0 & c & 0 \\
                                                                                                                  \end{array}
                                                                                                                \right),
$$ with $abc=-1$. This group is not a direct product. We can't give a general formula for this group and we only give a simple example.\\

\begin{Ex}$\\$
Consider the group $$\Gamma:=\left\langle
\left(
  \begin{array}{ccc}
    1 & 0 & 0 \\
    0 & -1 & 0 \\
    0 & 0 & -1 \\
  \end{array}
\right),\ \left(
  \begin{array}{ccc}
    0 & 1 & 0 \\
    0 & 0 & 1 \\
    1 & 0 & 0 \\
  \end{array}
\right),\ \left(
  \begin{array}{ccc}
    -1 & 0 & 0 \\
    0 & 0 & 1 \\
    0 & 1 & 0 \\
  \end{array}
\right)
 \right\rangle.$$
This group is isomorphic to the symmetric group $\go{S}_4$, so $\Gamma$ has $5$ conjugacy classes.\\
The series $P_\Gamma(t,\,u)$ verifies
$$\forall\ i\in\llbr 0,\,4\rrbr,\ P_\Gamma(t,\,u)_i=(1-tu)\frac{N(t,\,u)_i}{D(t)D(u)},$$
with $D(t)=(t-1)^3(t^2+t+1)(t^2+1)(t+1)^2$, and\\

\noindent $\begin{array}{rcl}
             N(t,\,u)_0 & = & {t}^{6}{u}^{6}+{t}^{5}{u}^{5}-{t}^{6}{u}^{3}-{t}^{3}{u}^{6}+{t}^{5}{u}^{3}+2\,{t}^{4}{u}^{4}+{t}^{3}{u}^{5}+{t
}^{4}{u}^{3}+{t}^{3}{u}^{4}+{t}^{6}+{t}^{5}u+2\,{t}^{4}{u}^{2}+4\,{t}^
{3}{u}^{3} \\
              &  & +2\,{t}^{2}{u}^{4}+t{u}^{5}+{u}^{6}+{t}^{3}{u}^{2}+{t}^{2}{u
}^{3}+{t}^{3}u+2\,{t}^{2}{u}^{2}+t{u}^{3}-{t}^{3}-{u}^{3}+tu+1,
           \end{array}$\\

\noindent $\begin{array}{rcl}
             N(t,\,u)_1 & = & {t}^{6}{u}^{3}+{t}^{5}{u}^{4}+{t}^{4}{u}^{5}+{t}^{3}{u}^{6}+{t}^{5}{u}
^{3}+{t}^{4}{u}^{4}+{t}^{3}{u}^{5}+{t}^{5}{u}^{2}+{t}^{4}{u}^{3}+{t}^{
3}{u}^{4}+{t}^{2}{u}^{5}+{t}^{4}{u}^{2}+{t}^{2}{u}^{4} \\
              &  & +{t}^{4}u+{t}^{3
}{u}^{2}+{t}^{2}{u}^{3}+t{u}^{4}+{t}^{3}u+{t}^{2}{u}^{2}+t{u}^{3}+{t}^
{3}+{t}^{2}u+t{u}^{2}+{u}^{3},
\end{array}$\\

\noindent $\begin{array}{rcl}
N(t,\,u)_2 & = & \left( {t}^{4}{u}^{2}+{t}^{3}{u}^{3}+{t}^{2}{u}^{4}+{t}^{3}{u}^{2}+{t}^{2}{u}^{3}+{t}^{3}u+t{u}^{3}+{t}^{2}u+t{u}^{2}+{t}^{2}+tu+{u}^{2}
 \right)  \left( {t}^{2}+1 \right)  \left( {u}^{2}+1 \right),
           \end{array}
$\\

\noindent $\begin{array}{rcl}
N(t,\,u)_3 & = & \left( {t}^{4}{u}^{2}+{t}^{3}{u}^{3}+{t}^{2}{u}^{4}+{t}^{3}u+t{u}^{3}+{t}^{2}+tu+{u}^{2} \right)  \left( {t}^{2}+t+1 \right)  \left( {u}^{2
}+u+1 \right),
           \end{array}
$\\

\noindent $\begin{array}{rcl}
N(t,\,u)_4 & = & ( {t}^{4}{u}^{3}+{t}^{3}{u}^{4}-{t}^{4}{u}^{2}-{t}^{3}{u}^{3}-{t}^{2}{u}^{4}+{t}^{4}u+2\,{t}^{3}{u}^{2}+2\,{t}^{2}{u}^{3}+t{u}^{4}-{t}
^{3}u-t{u}^{3}+{t}^{3}+2\,{t}^{2}u\\
 & & +2\,t{u}^{2}+{u}^{3}-{t}^{2}-tu-{u}^
{2}+t+u)  \left( {t}^{2}+t+1 \right)  \left( {u}^{2}+u+1
 \right).
           \end{array}
$\\
\end{Ex}


\section{\textsf{Exceptional subgroups of $\mathbf{SL}_3\mathbb{C}$ --- Types $E,\ F,\ G,\ H,\ I,\ J,\ K,\ L$}}

\noindent For every exceptional subgroup of
$\mathbf{SL}_3\mathbb{C}$, we begin by making the matrix $A^{(1)}$
explicit. Then we give a decomposition of $C_\Gamma:=2\,I-A^{(1)}-A^{(2)}+2\,Diag(A^{(1)})$ as
a sum of $p$ elements, with $p\in\{3,\,4\}$, so that~$C_\Gamma=p\,I-(\tau_0+\dots+\tau_{p-1})$, and
we give the graph associated to $C_\Gamma$. We also write the list $\Theta$ of eigenvalues of $A^{(1)}$.\\
Finally, we compute the sum of the series
$\displaystyle{P_\Gamma(t,\,u)=\frac{N_\Gamma(t,\,u)}{D_\Gamma(t,\,u)}}$.
In all the cases, the denominator is of the form
$D_\Gamma(t,\,u)_i=\widetilde{D_\Gamma}(t)_i\widetilde{D_\Gamma}(u)_i$.
Moreover, we will take the lowest common multiple $D_\Gamma(t)$ of
the $\widetilde{D_\Gamma}(t)_i$'s in order that all the
denominators are the same and have the form
$D_\Gamma(t)D_\Gamma(u)$, i.e.
$$\forall\ i\in \llbr 0,\, l\rrbr,\
P_\Gamma(t,\,u)_i=(1-tu)\frac{M_\Gamma(t,\,u)_i}{D_\Gamma(t)D_\Gamma(u)}.$$
Because of the to big size of the numerators, only the denominator and the relations between the numerators are given in the text: all the numerators may be found on the web.\\
We also give the Poincar\'e series of the invariant ring $\widehat{P}_\Gamma(t):=P_\Gamma(t,\,0)_0=P_\Gamma(0,\,t)_0$.\\

\subsection{\textsf{Type $E$}}

\noindent The group of type $E$ is the group $\langle S,\ T,\ V\rangle$, with
$$S:=\left(
       \begin{array}{ccc}
         1 & 0 & 0 \\
         0 & \zeta_{3} & 0 \\
         0 & 0 & \zeta_3^2 \\
       \end{array}
     \right),\
     T:=\left(
          \begin{array}{ccc}
            0 & 1 & 0 \\
            0 & 0 & 1 \\
            1 & 0 & 0 \\
          \end{array}
        \right),\
        V:=\frac{i}{\sqrt{3}}\left(
             \begin{array}{ccc}
               1 & 1 & 1 \\
               1 & \zeta_3 & \zeta_3^2 \\
               1 & \zeta_3^2 & \zeta_3 \\
             \end{array}
           \right).
$$
Here $l+1=14$, $\textrm{rank}(A^{(1)})=12$, $\Theta=( 3,\ -\zeta_3,\ -\zeta_3^2,\
0,\ 0,\ -1,\ \zeta_3,\ 1,\ \zeta_3^2,\ 1,\ \zeta_3,\ 3\,\zeta_3,\ \zeta_3^2,\ 3\,\zeta_3^2 )$,
$p=3$, and $\tau_0:=s_0s_1s_2s_3s_{12}s_{13},\ \tau_1:=s_5s_7s_{10}s_{11},\ \tau_2:=s_4s_6s_8s_{9}.$\\

\begin{figure}[h]
\begin{center}
\begin{tabular}{ccc}
 \includegraphics[width=4cm,height=4cm]{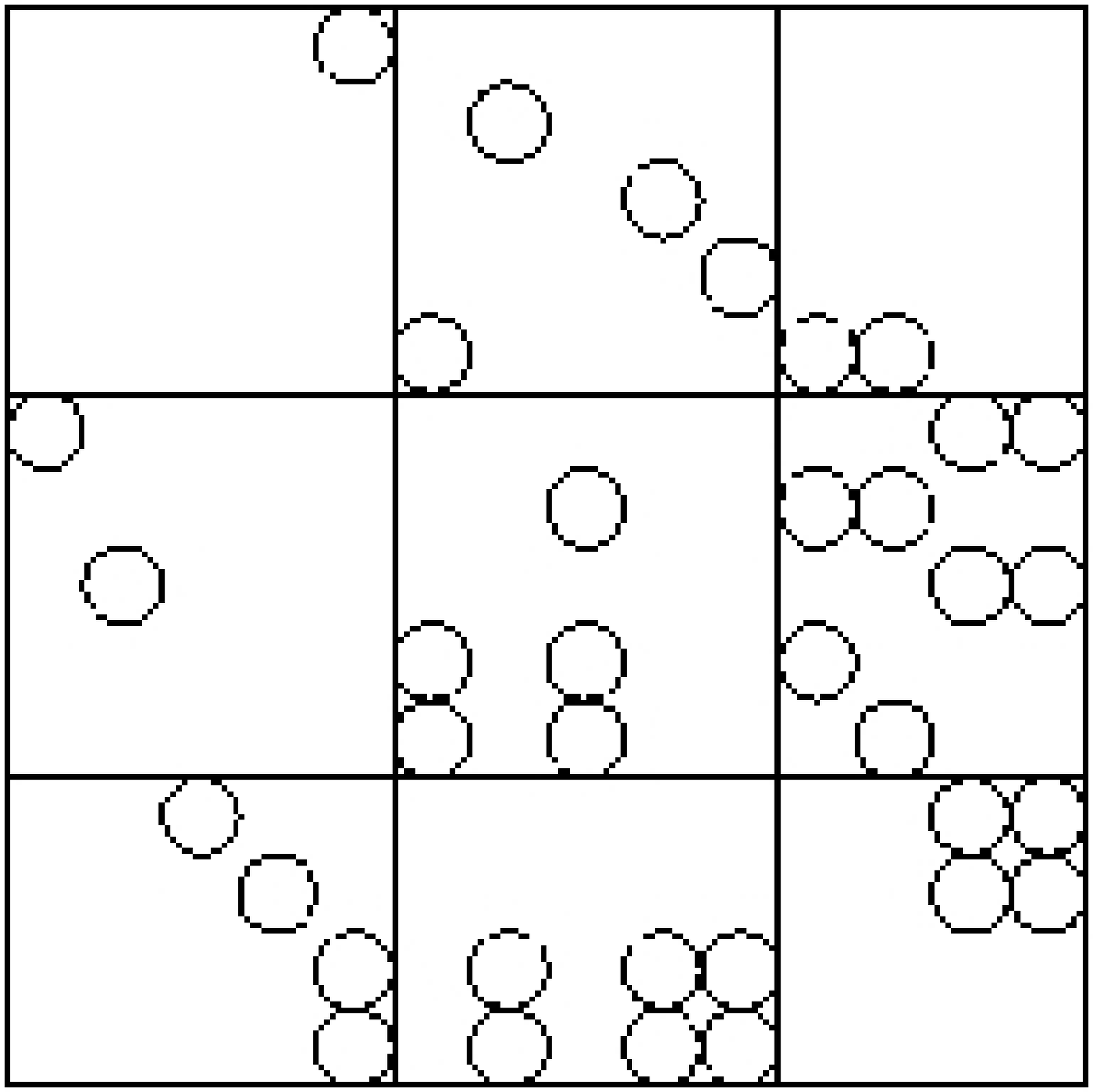} \\[-4cm]&
\begin{tiny}
\def\arrow(#1,#2){\ncline{-}{#1}{#2}}
  $\begin{array}{c@{\hskip .1cm}c@{\hskip .1cm}c@{\hskip .1cm}c@{\hskip .1cm}c@{\hskip .1cm}c@{\hskip .1cm}c}
      & & & \rnode{inv1}{}  \\[.15cm]
      & \rnode{3}{\Tcircle{3}} &   &  & & \rnode{1}{\Tcircle{1}} \\[.15cm]
      \rnode{11}{\Tcircle{11}} & & \rnode{9}{\Tcircle{9}} & &  \rnode{7}{\Tcircle{7}}  & &  \rnode{6}{\Tcircle{6}} \\[.15cm]
      \\[.15cm]
      & &  \rnode{12}{\Tcircle{12}} &   &  \rnode{13}{\Tcircle{13}} \\[.15cm]
      \\[.15cm]
        \rnode{5}{\Tcircle{5}} & & \rnode{4}{\Tcircle{4}} & &  \rnode{10}{\Tcircle{10}}  & &  \rnode{8}{\Tcircle{8}} \\[.15cm]
      & \rnode{0}{\Tcircle{0}} &   &  & & \rnode{2}{\Tcircle{2}} \\[.15cm]
      & & & \rnode{inv2}{}  \\[.15cm]
    \end{array}
    \everypsbox{\scriptstyle}
    \psset{nodesep=0pt,arm=.6,linearc=.4,angleA=0,angleB=90}
     \ncarc[arcangleA=-60, arcangleB=-40]{-}{5}{inv2}      \ncarc[arcangleA=-40, arcangleB=-60]{-}{inv2}{8}
     \ncarc[arcangleA=-60, arcangleB=-40]{-}{6}{inv1}      \ncarc[arcangleA=-40, arcangleB=-60]{-}{inv1}{11}
     \arrow(0,5) \arrow(1,7)   \arrow(2,10) \arrow(3,11) \arrow(4,0) \arrow(4,12) \arrow(4,13) \arrow(5,4)   \arrow(5,9) \arrow(6,1)
     \arrow(6,12) \arrow(6,13)   \arrow(7,6) \arrow(7,8) \arrow(7,9) \arrow(8,2) \arrow(8,12)   \arrow(8,13) \arrow(9,3) \arrow(9,12)
     \arrow(9,13) \arrow(10,4)   \arrow(10,6) \arrow(10,8) \arrow(11,4)  \arrow(11,9)   \arrow(12,5) \arrow(12,7) \arrow(12,10)
     \arrow(12,11) \arrow(13,5)   \arrow(13,7) \arrow(13,10) \arrow(13,11)
$\end{tiny} &
$\begin{array}{l}
    M_E(t,\,u)_2=M_E(t,\,u)_3\\
    M_E(t,\,u)_5=M_E(u,\,t)_4 \\
     M_E(t,\,u)_7=M_E(u,\,t)_6 \\
    M_E(t,\,u)_8=M_E(t,\,u)_9 \\
     M_E(t,\,u)_{10}=M_E(t,\,u)_{11}=M_E(u,\,t)_8 \\
      M_E(t,\,u)_{12}=M_E(t,\,u)_{13}
 \end{array}
$\\
\end{tabular}
\end{center}
\end{figure}

$$D_E(t)= \left( t-1 \right) ^{3} \left( {t}^{2}+t+1
\right) ^{3}
 \left( {t}^{2}+1 \right)  \left( {t}^
{4}-{t}^{2}+1 \right)  \left( t+1 \right) ^{2} \left( {t}^{2}-t+1
 \right) ^{2}$$
$$\widehat{P}_E(t)={\frac
{-{t}^{18}+{t}^{15}-{t}^{12}-{t}^{6}+{t}^{3}-1}{ \left( t-1
 \right) ^{3} \left( {t}^{2}+t+1 \right) ^{3} \left( {t}^{2}+1
 \right)  \left( {t}^{4}-{t}^{2}+1 \right)  \left( t+1 \right) ^{2}
 \left( {t}^{2}-t+1 \right) ^{2}}}.$$

\subsection{\textsf{Type $F$}}

\noindent The group of type $F$ is the group $\langle S,\ T,\ V,\ P\rangle$, with $S,\ T,\ V$ as for the type $E$, and
$$P:=\frac{1}{\sqrt{-3}}\left(
             \begin{array}{ccc}
               1 & 1 & \zeta_3^2 \\
               1 & \zeta_3 & \zeta_3 \\
               \zeta_3 & 1 & \zeta_3 \\
             \end{array}
           \right).
$$
Here $l+1=16$, $\textrm{rank}(A^{(1)})=15$, $\Theta=( 3,\
-\zeta_3,\ -\zeta_3^2,\ 0,\ -1,\ \zeta_3,\ 1,\ \zeta_3,\ 1,\ \zeta_3,\ 1,\
\zeta_3^2,\ \zeta_3^2, \zeta_3^2,\ 3\,\zeta_3,\ 3\,\zeta_3^2 )$,
$p=3$, and  $\tau_0:=s_0s_1s_2s_3s_{4}s_{15},\ \tau_1:=s_5s_7s_{9}s_{11}s_{13},\ \tau_2:=s_6s_8s_{10}s_{12}s_{14}$.

\begin{figure}[h]
\begin{center}
\begin{tabular}{ccc}
 \includegraphics[width=4cm,height=4cm]{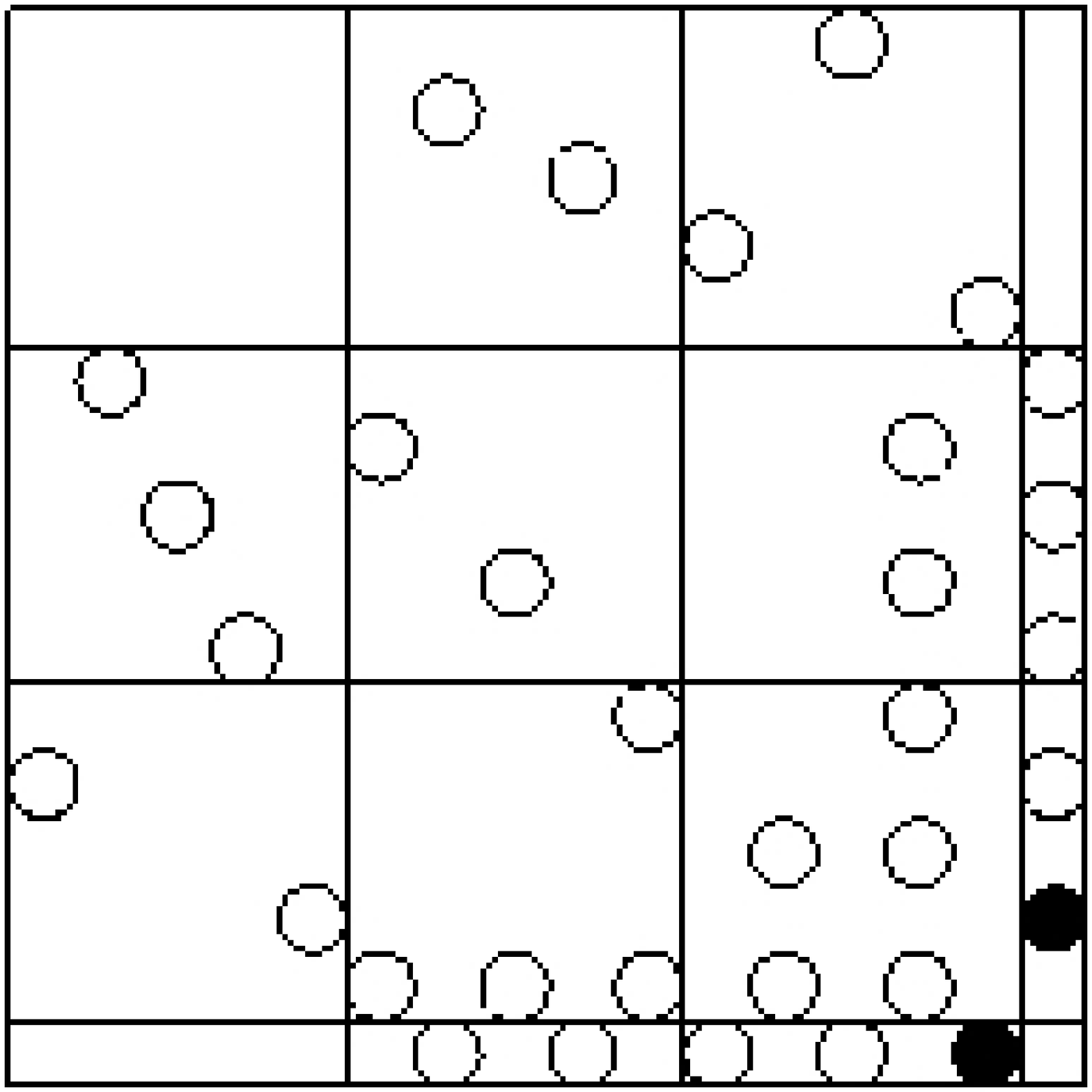} \\[-3.4cm]&
\begin{tiny}
\def\arrow(#1,#2){\ncline{-}{#1}{#2}}
$\begin{array}{c@{\hskip .1cm}c@{\hskip .1cm}c@{\hskip .1cm}c@{\hskip .1cm}c@{\hskip .1cm}c@{\hskip .1cm}c@{\hskip .1cm}
c@{\hskip .1cm}c@{\hskip .1cm}c@{\hskip .1cm}c@{\hskip .1cm}c@{\hskip .1cm}c@{\hskip .1cm}c@{\hskip .1cm}c@{\hskip .1cm}c}
      & & & & & & & \rnode{13}{\Tcircle{13}} &   &  & \\[.2cm]
    &  & \rnode{12}{\Tcircle{12}} &  & & \rnode{6}{\Tcircle{6}} & & & & \rnode{8}{\Tcircle{8}} & & & \rnode{10}{\Tcircle{10}} \\[.2cm]
  \rnode{inv}{ }  & \rnode{0}{\Tcircle{0}} & & & \rnode{1}{\Tcircle{1}} & & & \rnode{15}{\Tcircle{15}} & & & \rnode{2}{\Tcircle{2}} & & & \rnode{3}{\Tcircle{3}} & & \rnode{4}{\Tcircle{4}}\\[.2cm]
     & & \rnode{11}{\Tcircle{11}} & & & \rnode{5}{\Tcircle{5}} & & & & \rnode{7}{\Tcircle{7}} & & & \rnode{9}{\Tcircle{9}} \\[.2cm]
      & & & & & & & \rnode{14}{\Tcircle{14}} &   &  & \\[.2cm]
    \end{array}
    \everypsbox{\scriptstyle}
    \psset{nodesep=0pt,arm=.6,linearc=.4,angleA=0,angleB=90}
     \ncarc[arcangleA=-30, arcangleB=-30]{-}{4}{13}
         \ncarc[arcangleA=30, arcangleB=30]{-}{4}{14}
             \ncarc[arcangleA=-90, arcangleB=-70]{-}{13}{inv}
             \ncarc[arcangleA=90, arcangleB=70]{-}{14}{inv}
     \arrow(0,11) \arrow(1,5) \arrow(2,7) \arrow(3,9)
     \arrow(5,6) \arrow(5,14) \arrow(6,1) \arrow(6,15) \arrow(7,8)
     \arrow(7,14) \arrow(8,2) \arrow(8,15) \arrow(9,10) \arrow(9,14)
     \arrow(10,3) \arrow(10,15) \arrow(11,12) \arrow(11,14) \arrow(12,0)
     \arrow(12,15) \arrow(13,6) \arrow(13,8) \arrow(13,10) \arrow(13,12)
      \arrow(14,15) \arrow(15,5) \arrow(15,7)
     \arrow(15,9) \arrow(15,11)\arrow(15,13)
\ncline[doubleline=true]{-}{13}{15} \ncline[doubleline=true]{-}{14}{15}$
\end{tiny}
 & $\begin{array}{l}
     M_F(t,\,u)_1=M_F(t,\,u)_2=M_F(t,\,u)_3 \\
      M_F(t,\,u)_5=M_F(t,\,u)_7=M_F(t,\,u)_9  \\
    M_F(t,\,u)_{12}=M_F(u,\,t)_{11} \\
     M_F(t,\,u)_{14}=M_F(u,\,t)_{13} \\
      M_F(t,\,u)_6=M_F(t,\,u)_8=M_F(t,\,u)_{10}\\
      =M_F(u,\,t)_5
    \end{array}
 $\\
\end{tabular}
\end{center}
\end{figure}

$$D_F(t)=\left( t-1 \right) ^{3} \left( {t}^{2}+t+1
\right) ^{3}  \left( {t}^{2}+1 \right)  \left( {t}^ {4}-{t}^{2}+1
\right)  \left( t+1 \right) ^{2} \left( {t}^{2}-t+1
 \right) ^{2}$$
$$\widehat{P}_F(t)={\frac {-{t}^{18}+{t}^{15}-{t}^{9}+{t}^{3}-1}{ \left( t-1 \right) ^{3}
 \left( {t}^{2}+t+1 \right) ^{3} \left( {t}^{2}+1 \right)  \left( {t}^
{4}-{t}^{2}+1 \right)  \left( t+1 \right) ^{2} \left( {t}^{2}-t+1
 \right) ^{2}}}.$$

\subsection{\textsf{Type $G$}}

\noindent The group of type $G$ is the group $\langle S,\ T,\ V,\ U\rangle$, with $S,\ T,\ V$ as for the type $E$, and
$$U:=\left(
             \begin{array}{ccc}
               \zeta_9^2 & 0 & 0 \\
               0 & \zeta_9^2 & 0 \\
               0 & 0 & \zeta_9^2\zeta_3 \\
             \end{array}
           \right).
$$
Here $l+1=24$, $\textrm{rank}(A^{(1)})=21$,
$$\Theta=( 3,\
-\zeta_3,\ -\zeta_3^2,\ 0,\ 0,\ -\zeta_9^7,\ -\zeta_9^4,\ \zeta_9^4+\zeta_9^7,\ 0,\
-\zeta_9^5,\ -\zeta_9^2, \zeta_9^2+\zeta_9^5,\ -1,\ \zeta_3,\ 1,\ \zeta_3^2,\
3\,\zeta_3,\ 3\,\zeta_3^2,$$
$$-2\,\zeta_9^4-\zeta_9^7,\ \zeta_9^4-\zeta_9^7,\
-\zeta_9^2-2\,\zeta_9^5,\ -\zeta_9^2+\zeta_9^5,\ 2\,\zeta_9^2+\zeta_9^5,\ \zeta_9^4+2\,\zeta_9^7 ),$$
$p=3$, and $\tau_0:=s_0s_1s_2s_3s_{4}s_{5}s_{6}s_{19}s_{20}s_{21},\ \tau_1:=s_7s_8s_{9}s_{13}s_{14}s_{15}s_{22},\ \tau_2:=s_{10}s_{11}s_{12}s_{16}s_{17}s_{18}s_{23}.$

\vspace{1cm}

\begin{figure}[h]
\begin{center}
\begin{tabular}{ccc}
 \includegraphics[width=4cm,height=4cm]{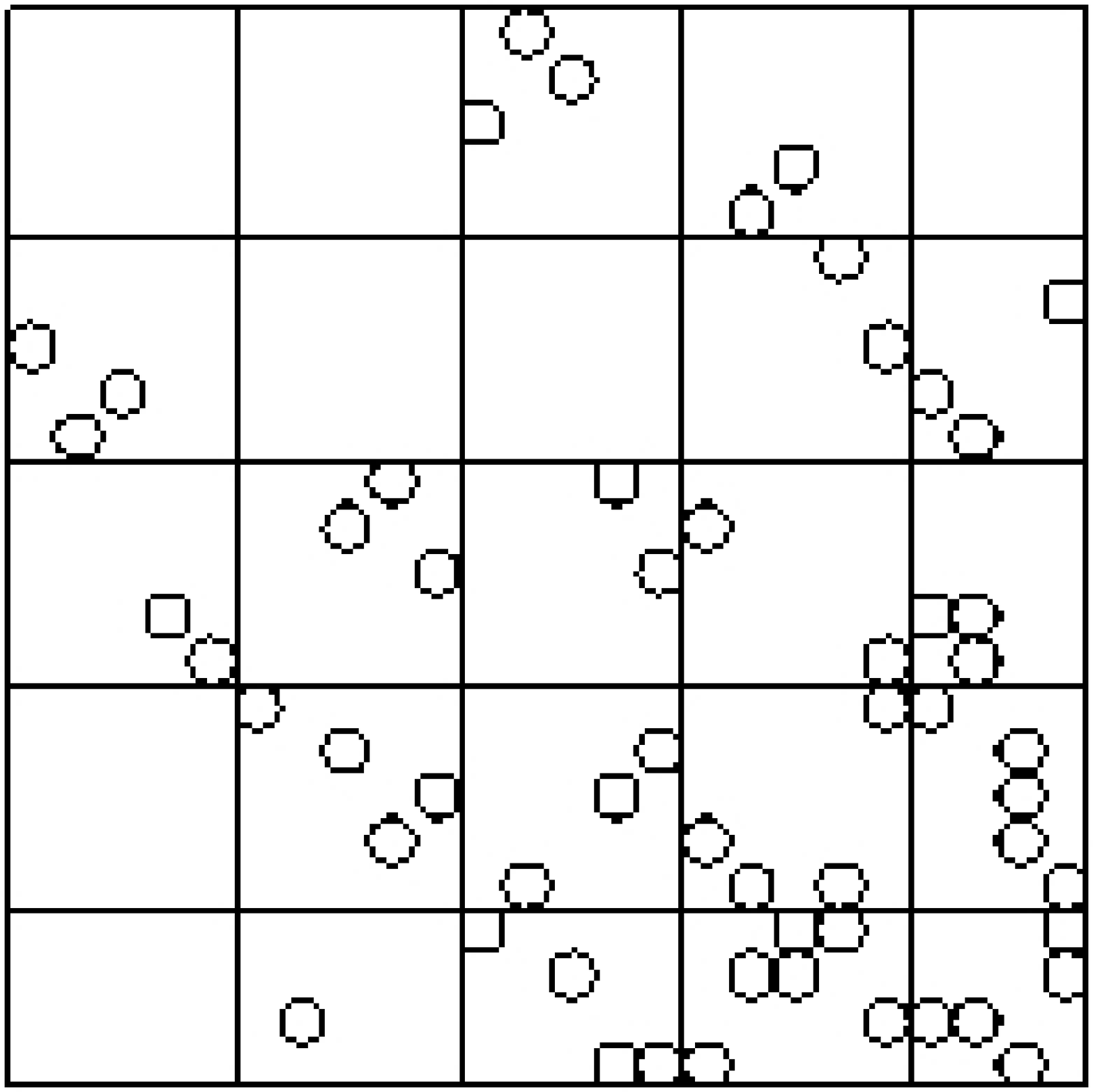} \\[-5cm]&
\begin{tiny}
\def\arrow(#1,#2){\ncline{-}{#1}{#2}}
$\begin{array}{c@{\hskip .1cm}c@{\hskip .1cm}c@{\hskip .1cm}c@{\hskip .1cm}c@{\hskip .1cm}c@{\hskip .1cm}c@{\hskip .1cm}c@{\hskip .1cm}
c@{\hskip .1cm}c@{\hskip .1cm}c}
    & & & & & \rnode{5}{\Tcircle{5}}  \\[.1cm]
     & & & & \rnode{18}{\Tcircle{18}} &   &  \rnode{15}{\Tcircle{15}} \\[.1cm]
     & \rnode{8}{\Tcircle{8}} & & & & & & & &  \rnode{11}{\Tcircle{11}} \\[.1cm]
     \rnode{2}{\Tcircle{2}} & & & \rnode{20}{\Tcircle{20}} & & & & \rnode{19}{\Tcircle{19}} & & &\rnode{0}{\Tcircle{0}}  \\[.1cm]
     & \rnode{10}{\Tcircle{10}} & & & & & & & &  \rnode{7}{\Tcircle{7}} \\[.1cm]
     & & & & \rnode{22}{\Tcircle{22}} &   &  \rnode{23}{\Tcircle{23}} \\[.1cm]
    & & & & & \\[.1cm]
     & \rnode{13}{\Tcircle{13}} & & & & \rnode{6}{\Tcircle{6}} & & & &  \rnode{16}{\Tcircle{16}} \\[.1cm]
\rnode{3}{\Tcircle{3}} &  & \rnode{17}{\Tcircle{17}} & & & &  & &  \rnode{14}{\Tcircle{14}} & & \rnode{4}{\Tcircle{4}}  \\[.1cm]
       & & & & & \rnode{21}{\Tcircle{21}}  \\[.1cm]
     & & & & \rnode{9}{\Tcircle{9}} &   &  \rnode{12}{\Tcircle{12}} \\[.1cm]
    & & & & & \rnode{1}{\Tcircle{1}}  \\[.1cm]
    \end{array}
    \everypsbox{\scriptstyle}
    \psset{nodesep=0pt,arm=.6,linearc=.4,angleA=0,angleB=90}
         \arrow(0,7) \arrow(1,9)   \arrow(2,8) \arrow(3,13) \arrow(4,14) \arrow(5,15) \arrow(6,22) \arrow(7,11)   \arrow(7,16) \arrow(8,10)
     \arrow(8,18) \arrow(9,12)   \arrow(9,17) \arrow(10,2) \arrow(10,20) \arrow(11,0) \arrow(11,19)   \arrow(12,1) \arrow(12,21) \arrow(13,10)
     \arrow(13,17) \arrow(13,23)   \arrow(14,12) \arrow(14,16) \arrow(14,23)  \arrow(15,11)    \arrow(15,18) \arrow(15,23)
     \arrow(16,4) \arrow(16,19)   \arrow(16,21) \arrow(17,3) \arrow(17,20) \arrow(17,21) \arrow(18,5) \arrow(18,19)
      \arrow(18,20) \arrow(19,7) \arrow(19,14)  \arrow(19,15) \arrow(19,22) \arrow(20,8) \arrow(20,13) \arrow(20,15) \arrow(20,22)
      \arrow(21,9) \arrow(21,13) \arrow(21,14) \arrow(21,22) \arrow(22,16) \arrow(22,17) \arrow(22,18) \arrow(22,23)
       \arrow(23,6) \arrow(23,19) \arrow(23,20) \arrow(23,21)$
\end{tiny} &
$\begin{array}{l}
       M_G(t,\,u)_{2}=M_G(u,\,t)_{1} \\
        M_G(t,\,u)_{5}=M_G(u,\,t)_{4} \\
         M_G(t,\,u)_{10}=M_G(u,\,t)_{9}\\
        M_G(t,\,u)_{11}=M_G(u,\,t)_{7} \\
         M_G(t,\,u)_{12}=M_G(u,\,t)_{8} \\
          M_G(t,\,u)_{16}=M_G(u,\,t)_{15}\\
    M_G(t,\,u)_{18}=M_G(u,\,t)_{14}  \\
     M_G(t,\,u)_{21}=M_G(u,\,t)_{20} \\
      M_G(t,\,u)_{23}=M_G(u,\,t)_{22}
   \end{array}
$\\
\end{tabular}
\end{center}
\end{figure}

\noindent $D_G(t)=\left( t-1 \right) ^{3} \left( {t}^{6}-{t}^{3}+1 \right)
 \left( {t}^{2}+t+1 \right) ^{3} \left( t+1 \right) ^{2}
  \left( {t}^{2}+1 \right)  \left( {t}^{4}-{t}^{2}+1 \right)  \left(
{t}^{2}-t+1 \right) ^{2}  \left( {t}^{6}+{t}^{3}+1 \right) ^{2} $\\
$$\widehat{P}_G(t)=\frac {-{t}^{18}-{t}^{36}-1}{ \left( t-1 \right) ^{3} \left( {t}^{6}-
{t}^{3}+1 \right)  \left( {t}^{2}+t+1 \right) ^{3} \left( t+1 \right)
^{2} \left( {t}^{2}+1 \right)  \left( {t}^{4}-{t}^{2}+1 \right)
 \left( {t}^{2}-t+1 \right) ^{2} \left( {t}^{6}+{t}^{3}+1 \right) ^{2}
}$$

\subsection{\textsf{Type $H$}}

\noindent The group of type $H$, isomorphic to the alternating group $\go{A}_5$, is the group $\langle S,\ U,\ T\rangle$, with
$$S:=\left(
       \begin{array}{ccc}
         1 & 0 & 0 \\
         0 & \zeta_{5}^4 & 0 \\
         0 & 0 & \zeta_5 \\
       \end{array}
     \right),\
     U:=\left(
          \begin{array}{ccc}
            -1 & 0 & 0 \\
            0 & 0 & -1 \\
            0 & -1 & 0 \\
          \end{array}
        \right),\
        T:=\frac{1}{\sqrt{5}}\left(
             \begin{array}{ccc}
               1 & 1 & 1 \\
               2 & \zeta_5^2+\zeta_5^3 & \zeta_5+\zeta_5^4 \\
               2 & \zeta_5+\zeta_5^4 & \zeta_5^2+\zeta_5^3 \\
             \end{array}
           \right).
$$
Here $l+1=5$, $\textrm{rank}(A^{(1)})=4$, $A^{(1)}$ is symmetric, $\Theta=( 3,\ -1,\ -\zeta_5^2-\zeta_5^3,\ -\zeta_5-\zeta_5^4,\ 0 )$,\\
$p=3$, and $\tau_0:=s_0s_3,\ \tau_1:=s_1s_2,\ \tau_2:=s_4.$

\begin{figure}[h]
\begin{center}
\begin{tabular}{cccc}
 \begin{scriptsize}$A^{(1)}=\left( \begin {array}{ccccc}
0&0&1&0&0\\\noalign{\medskip}0&0&0&1&1\\\noalign{\medskip}1&0&1&0&1\\\noalign{\medskip}0&1&0&1&1
\\\noalign{\medskip}0&1&1&1&1\end {array} \right)$
\end{scriptsize} &
\begin{scriptsize}
\def\arrow(#1,#2){\ncline{-}{#1}{#2}}
$\begin{array}{c@{\hskip .15cm}c@{\hskip .15cm}c@{\hskip .15cm}c}
      \rnode{1}{\Tcircle{1}} &   &  & \\[.2cm]
      & \rnode{4}{\Tcircle{4}} & \rnode{2}{\Tcircle{2}}  & \rnode{0}{\Tcircle{0}} \\[.2cm]
      \rnode{3}{\Tcircle{3}} & & &  \\[.2cm]
    \end{array}
    \everypsbox{\scriptstyle}
    \psset{nodesep=0pt,arm=.6,linearc=.4,angleA=0,angleB=90}
     \arrow(1,3) \arrow(1,4)   \arrow(4,3) \arrow(2,4) \arrow(2,0)$
\end{scriptsize} &
$\begin{array}{c}
\begin{tiny}
D_H(t)=\left( t-1 \right) ^{3}  \left( {t}^{2}+t+1
\right)  \left( {t}^{4}+{t}^{3}+{t}^{
2}+t+1 \right)  \left( t+1 \right) ^{2}
\end{tiny} \\ \\
\widehat{P}_H(t)=\frac {-{t}^{8}-{t}^{7}+{t}^{5}+{t}^{4}+{t}^{3}-t-1}{ \left( t-1
 \right) ^{3} \left( {t}^{2}+t+1 \right)  \left( {t}^{4}+{t}^{3}+{t}^{
2}+t+1 \right)  \left( t+1 \right) ^{2}}
 \end{array}$
\\
\end{tabular}
\end{center}
\end{figure}

\subsection{\textsf{Type $I$}}

\noindent The group of type $I$ is the group $\langle S,\ T,\ R\rangle$, with
$$S:=\left(
       \begin{array}{ccc}
         \zeta_7 & 0 & 0 \\
         0 & \zeta_{7}^2 & 0 \\
         0 & 0 & \zeta_7^4 \\
       \end{array}
     \right),\
     T:=\left(
          \begin{array}{ccc}
            0 & 1 & 0 \\
            0 & 0 & 1 \\
            1 & 0 & 0 \\
          \end{array}
        \right),\
        R:=\frac{i}{\sqrt{7}}\left(
             \begin{array}{ccc}
               \zeta_7^4-\zeta_7^3 & \zeta_7^2-\zeta_7^5 & \zeta_7-\zeta_7^6 \\
               \zeta_7^2-\zeta_7^5 & \zeta_7-\zeta_7^6 & \zeta_7^4-\zeta_7^3 \\
               \zeta_7-\zeta_7^6 & \zeta_7^4-\zeta_7^3 & \zeta_7^2-\zeta_7^5 \\
             \end{array}
           \right).
$$
Here $l+1=6$, $\textrm{rank}(A^{(1)})=5$, $\Theta=( 3,\ 0,\ 1,\ \zeta_7+\zeta_7^2+\zeta_7^4,\ \zeta_7^3+\zeta_7^5+\zeta_7^6,\ -1 )$,\\
$p=4$, and $\tau_0:=s_5s_0,\ \tau_1:=s_1s_4,\ \tau_2:=s_2,\ \tau_3:=s_3$.

\begin{figure}[h]
\begin{center}
\begin{tabular}{ccc}
 \begin{scriptsize}$A^{(1)}=\left( \begin {array}{cccccc}
0&1&0&0&0&0\\\noalign{\medskip}0&0&1&1&0&0\\\noalign{\medskip}1&0&0&0&0&1\\\noalign{\medskip}0&0&1&0&1&1
\\\noalign{\medskip}0&0&0&1&1&1\\\noalign{\medskip}0&1&0&1&1&1
\end {array} \right)$
\end{scriptsize} &
\begin{scriptsize}
\def\arrow(#1,#2){\ncline{-}{#1}{#2}}
$\begin{array}{c@{\hskip .25cm}c@{\hskip .25cm}c@{\hskip .25cm}c}
       & \rnode{1}{\Tcircle{1}}  & \rnode{3}{\Tcircle{3}} & \\[.2cm]
    \rnode{0}{\Tcircle{0}}  &  & & \rnode{4}{\Tcircle{4}} \\[.2cm]
       & \rnode{2}{\Tcircle{2}} & \rnode{5}{\Tcircle{5}} &  \\[.2cm]
    \end{array}
    \everypsbox{\scriptstyle}
    \psset{nodesep=0pt,arm=.6,linearc=.4,angleA=0,angleB=90}
          \arrow(0,1) \arrow(0,2)   \arrow(1,2) \arrow(1,3) \arrow(1,5)
               \arrow(2,3) \arrow(2,5)   \arrow(3,5) \arrow(4,3) \arrow(5,4)$
\end{scriptsize}
$\begin{array}{l}
  D_I(t)=\left( t-1 \right) ^{3}\left( {t}^{2}+t+1
 \right) \left( {t}^{2}+1 \right)  \left( t+1 \right) ^{2}\\
  \ \ \ \ \ \ \ \ \left( {t}^{6}+{t}^{5}+{t}^{4}+{t}^{3}+{t}^{2}+t+1 \right)\\ \\
  M_I(t,\,u)_{2}=M_I(u,\,t)_{1}
\end{array}$
\\
\end{tabular}
\end{center}
\end{figure}
$$\widehat{P}_I(t)=\frac {-{t}^{12}-{t}^{11}+{t}^{9}+{t}^{8}-{t}^{6}+{t}^{4}+{t}^{3}-t-1
}{ \left( t-1 \right) ^{3} \left( {t}^{2}+t+1 \right)  \left(
{t}^{6}+ {t}^{5}+{t}^{4}+{t}^{3}+{t}^{2}+t+1 \right)  \left(
{t}^{2}+1 \right) \left( t+1 \right) ^{2}}$$

\subsection{\textsf{Type $J$}}

\noindent The group of type $J$ is the group $\langle S,\ U,\ T,\ W\rangle$, with $S,\ U,\ T$ as for the type $H$, and
$W:=Diag(j,\,j,\,j).$ It is the direct product of the group of type $H$ and the center of $\mathbf{SL}_3\mathbb{C}$. Here $l+1=15$, $\textrm{rank}(A^{(1)})=12$,\\
$\Theta=( 3,\
-1,\ -\zeta_5^2-\zeta_5^3,\ -\zeta_5-\zeta_5^4,\ -\zeta_3^2,\ -\zeta_3,\ 3\zeta_3,\
 -\zeta_{15}^{11}-\zeta_{15}^{14},\ -\zeta_{15}^2-\zeta_{15}^8,\ 3\zeta_{3}^2,\
 -\zeta_{15}^7-\zeta_{15}^{13},\\
 -\zeta_{15}-\zeta_{15}^4,\ 0,\ 0,\ 0 ),$
 $p=3$, and $\tau_0:=s_{2}s_{5}s_{7}s_{10}s_{13},\ \tau_1:=s_1s_6s_{8}s_{11}s_{14},\ \tau_2:=s_0s_{3}s_{4}s_{9}s_{12}$.\vspace{.4cm}

\begin{figure}[h]
\begin{center}
\begin{tabular}{ccc}
 \includegraphics[width=4cm,height=4cm]{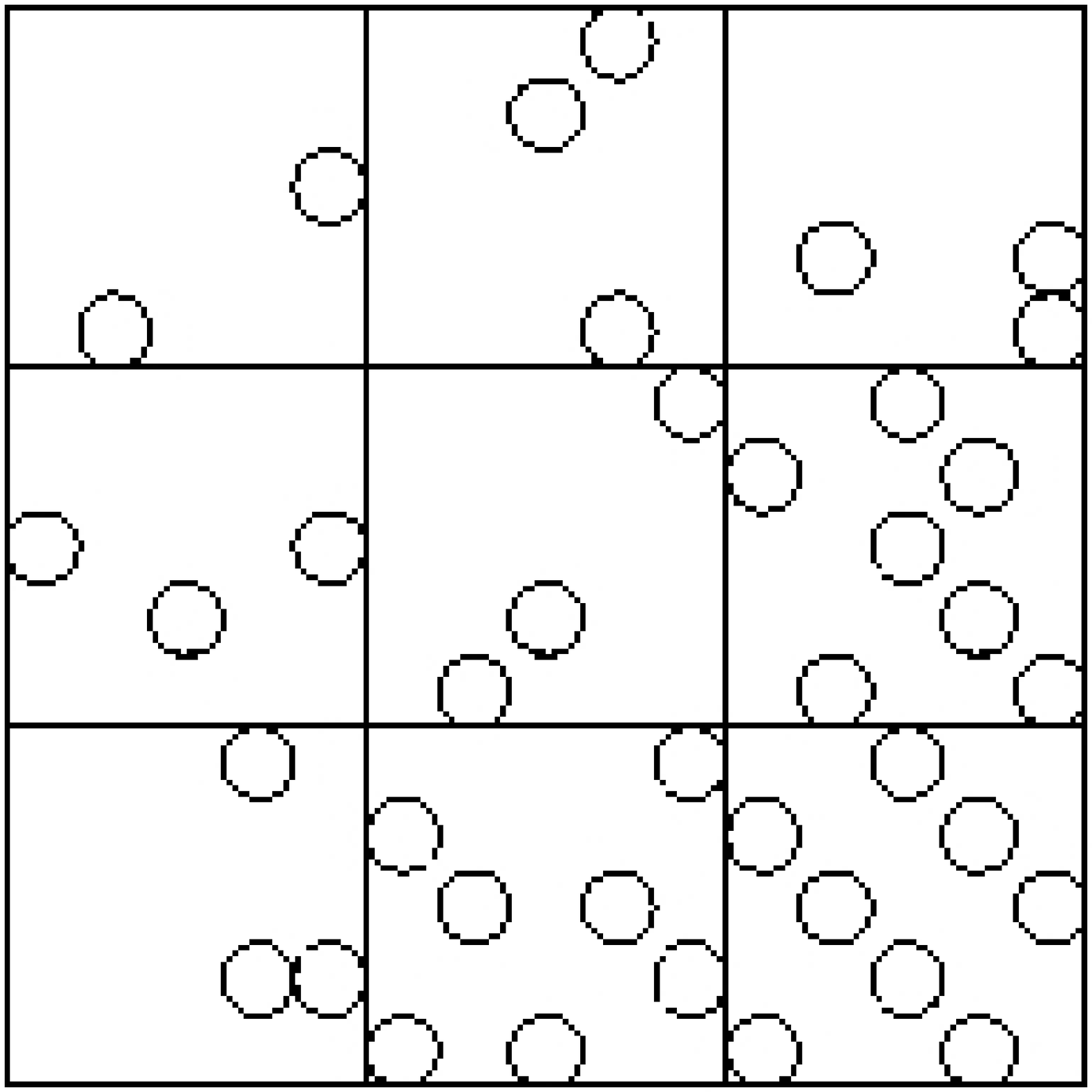} \\[-4.5cm]&
\begin{tiny}
\def\arrow(#1,#2){\ncline{-}{#1}{#2}}
$\begin{array}{c@{\hskip .15cm}c@{\hskip .15cm}c@{\hskip .15cm}c@{\hskip .15cm}c@{\hskip .15cm}c@{\hskip .15cm}c@{\hskip .15cm}c@{\hskip .15cm}c}
     & & & & \rnode{10}{\Tcircle{10}}  \\[.3cm]
      & \rnode{3}{\Tcircle{3}} &   &  \rnode{11}{\Tcircle{11}} & & \rnode{9}{\Tcircle{9}} & & \rnode{6}{\Tcircle{6}} \\[.3cm]
       & & & \rnode{inv2}{} & \rnode{5}{\Tcircle{5}} & \rnode{inv3}{}  \\[.2cm]
      \\[.2cm]
      & &  \rnode{14}{\Tcircle{14}} &   &  \rnode{13}{\Tcircle{13}} & & \rnode{12}{\Tcircle{12}} \\[.15cm]
       & & & & \rnode{inv1}{}  \\[.15cm]
        & & \rnode{4}{\Tcircle{4}} & & & & \rnode{8}{\Tcircle{8}} \\[.15cm]
      & & & & \rnode{7}{\Tcircle{7}} \\[.15cm]
            \\[.1cm]
      & &   \rnode{1}{\Tcircle{1}} & & \rnode{2}{\Tcircle{2}} & & \rnode{0}{\Tcircle{0}}
    \end{array}
    \everypsbox{\scriptstyle}
    \psset{nodesep=0pt,arm=.6,linearc=.4,angleA=0,angleB=90}
     \ncarc[arcangleA=-30, arcangleB=-25]{-}{11}{inv2}      \ncarc[arcangleA=-25, arcangleB=-30]{-}{inv2}{12}
     \ncarc[arcangleA=-30, arcangleB=-25]{-}{14}{inv3}      \ncarc[arcangleA=-25, arcangleB=-30]{-}{inv3}{9}
     \ncarc[arcangleA=-40, arcangleB=-40]{14}{12}
     \arrow(0,7) \arrow(1,4)   \arrow(2,8) \arrow(3,10) \arrow(3,13) \arrow(4,2) \arrow(4,7) \arrow(4,13)   \arrow(5,11) \arrow(5,14)
     \arrow(6,9) \arrow(6,12)   \arrow(7,1) \arrow(7,8) \arrow(7,14) \arrow(8,0) \arrow(8,4)   \arrow(8,12) \arrow(9,5) \arrow(9,10)
     \arrow(9,13) \arrow(10,6)   \arrow(10,11) \arrow(10,14) \arrow(11,3)  \arrow(11,9)    \arrow(12,5) \arrow(12,7)
     \arrow(12,10) \arrow(12,13)   \arrow(13,6) \arrow(13,8) \arrow(13,11) \arrow(13,14) \arrow(14,3) \arrow(14,4)$
\end{tiny} &
$\begin{array}{l}
  M_J(t,\,u)_{2}=M_J(u,\,t)_{1} \\
  M_J(t,\,u)_{6}=M_J(u,\,t)_{5} \\
  M_J(t,\,u)_{8}=M_J(u,\,t)_{7} \\
  M_J(t,\,u)_{11}=M_J(u,\,t)_{10}\\
  M_J(t,\,u)_{14}=M_J(u,\,t)_{13}
\end{array}$
\\
\end{tabular}
\end{center}
\end{figure}
$$D_J(t)=\left( {t}^{4}+{t}^{3}+{t}^{2}+t+1 \right)
\left( {t}^{ 8}-{t}^{7}+{t}^{5}-{t}^{4}+{t}^{3}-t+1 \right)
\left( t+1 \right) ^{2 } \left( {t}^{2}-t+1 \right) ^{2} \left(
t-1 \right) ^{3} \left( {t}^{ 2}+t+1 \right) ^{3}$$
$$\widehat{P}_J(t)={\frac {-{t}^{24}-{t}^{12}-1}{ \left( {t}^{4}+{t}^{3}+{t}^{2}+t+1
 \right)  \left( {t}^{8}-{t}^{7}+{t}^{5}-{t}^{4}+{t}^{3}-t+1 \right)
 \left( t+1 \right) ^{2} \left( {t}^{2}-t+1 \right) ^{2} \left( t-1
 \right) ^{3} \left( {t}^{2}+t+1 \right) ^{3}}}$$

\subsection{\textsf{Type $K$}}

\noindent The group of type $K$ is the group $\langle S,\ T,\ R,\ W\rangle$, with $S,\ T,\ R$ as for the type $I$, and
$W:=Diag(j,\,j,\,j)$. It is the direct product of the group of type $I$ and the center of $\mathbf{SL}_3\mathbb{C}$. Here $l+1=18$, $\textrm{rank}(A^{(1)})=15$,
$$\begin{array}{c}
   \Theta=( 3,\ 0,\
0,\ 0,\ 3\zeta_{3},\ 3\zeta_{3}^2,\ 1,\ \zeta_{7}+\zeta_{7}^2+\zeta_{7}^4,\ \zeta_{7}^3+\zeta_{7}^5+\zeta_{7}^6,\
  -1,\ \zeta_{3}^2,\ \zeta_{21}^2+\zeta_{21}^8+\zeta_{21}^{11},\ \zeta_{3}, \\
    \zeta_{21}+\zeta_{21}^4+\zeta_{21}^{16},\
  \zeta_{21}^5+\zeta_{21}^{17}+\zeta_{21}^{20},\ \zeta_{21}^{10}+\zeta_{21}^{13}+\zeta_{21}^{19},\ -\zeta_{3},\ -\zeta_{3}^2 ),
  \end{array}$$
$p=3$, and $\tau_0:=s_1s_5s_6s_{11}s_{14}s_{17},\ \tau_1:=s_0s_3s_{4}s_{9}s_{12}s_{15},\ \tau_2:=s_2s_7s_8s_{10}s_{13}s_{16}$.\vspace{1cm}

\begin{figure}[h]
\begin{center}
\begin{tabular}{ccc}
 \includegraphics[width=4cm,height=4cm]{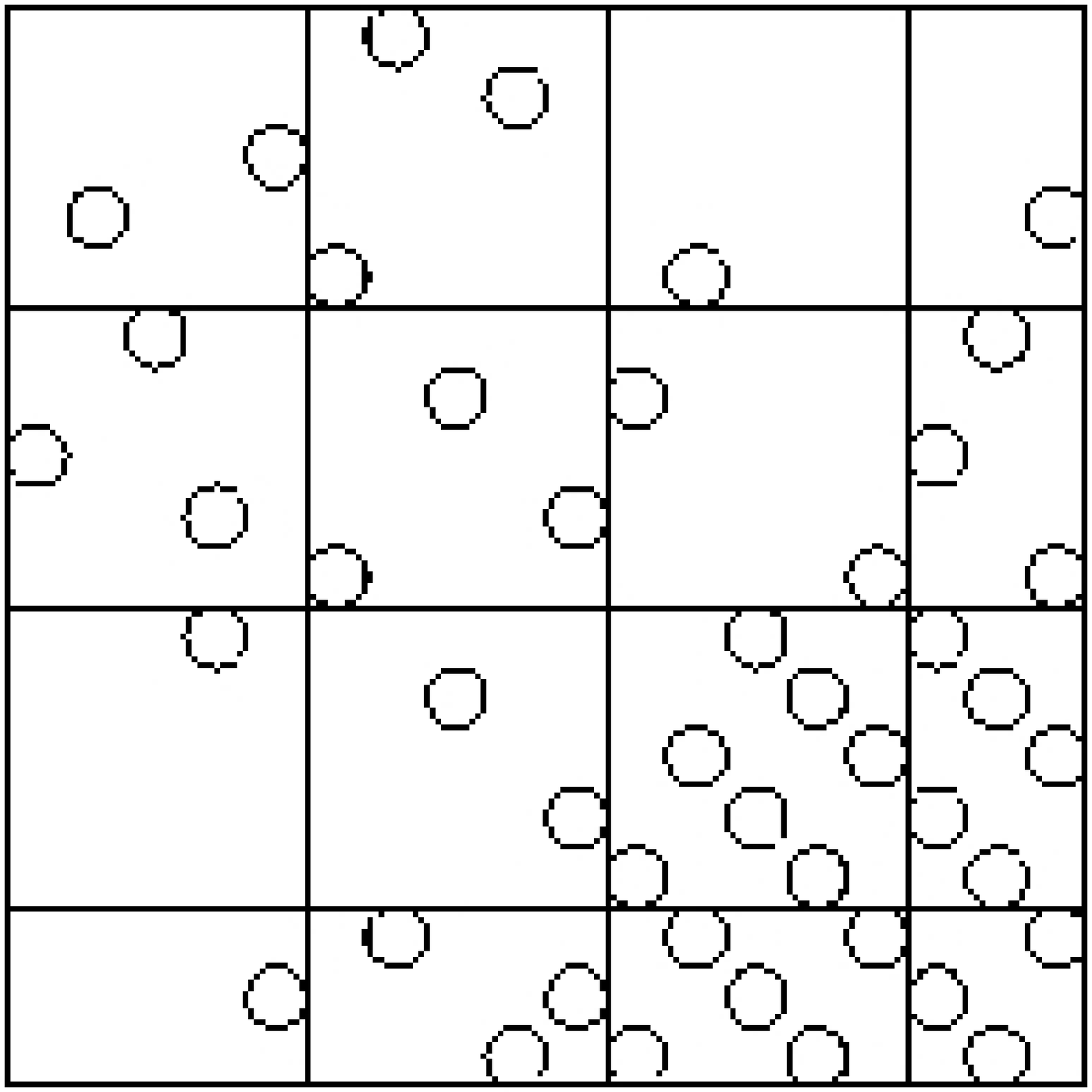} \\[-4.8cm]&
\begin{tiny}
\def\arrow(#1,#2){\ncline{-}{#1}{#2}}
$\begin{array}{c@{\hskip .1cm}c@{\hskip .1cm}c@{\hskip .1cm}c@{\hskip .1cm}c@{\hskip .1cm}c@{\hskip .1cm}c@{\hskip .1cm}c@{\hskip .1cm}c}
    & & & & \rnode{2}{\Tcircle{2}}  \\[.1cm]
     & & & \rnode{5}{\Tcircle{5}} &   &  \rnode{4}{\Tcircle{4}} \\[.1cm]
     \\[.0cm]
       & & & & \rnode{16}{\Tcircle{16}}  \\[.1cm]
       &  & \rnode{9}{\Tcircle{9}} & & & & \rnode{11}{\Tcircle{11}}  \\[.1cm]
       \\[.0cm]
              & & & & \rnode{13}{\Tcircle{13}}  \\[.1cm]
     & & & \rnode{14}{\Tcircle{14}} &   &  \rnode{12}{\Tcircle{12}} \\[.1cm]
     \\[.1cm]
      &  & \rnode{17}{\Tcircle{17}} & & & &  \rnode{15}{\Tcircle{15}}  \\[.1cm]
       & & & & \rnode{10}{\Tcircle{10}}  \\[.1cm]
       & \rnode{8}{\Tcircle{8}} & & & & & & \rnode{7}{\Tcircle{7}}  \\[.1cm]
        \rnode{1}{\Tcircle{1}} & & \rnode{3}{\Tcircle{3}} & & & & \rnode{6}{\Tcircle{6}} & & \rnode{0}{\Tcircle{0}}
    \end{array}
    \everypsbox{\scriptstyle}
    \psset{nodesep=0pt,arm=.6,linearc=.4,angleA=0,angleB=90}
     \ncarc[arcangleA=-30, arcangleB=-10]{-}{16}{17}      \ncarc[arcangleA=30, arcangleB=10]{-}{16}{15}
\ncarc[arcangleA=30, arcangleB=10]{-}{13}{17}      \ncarc[arcangleA=-30, arcangleB=-10]{-}{13}{15}
     \arrow(0,7) \arrow(1,3)   \arrow(2,5) \arrow(3,8) \arrow(3,10) \arrow(4,2) \arrow(4,16) \arrow(5,4)   \arrow(5,9) \arrow(6,0)
     \arrow(6,15) \arrow(7,6)   \arrow(7,11) \arrow(8,1) \arrow(8,17) \arrow(9,8) \arrow(9,13)   \arrow(9,16) \arrow(10,6) \arrow(10,14)
     \arrow(10,17) \arrow(11,4)   \arrow(11,12) \arrow(11,15) \arrow(12,10)  \arrow(12,13)    \arrow(12,16) \arrow(13,11)
     \arrow(13,14)   \arrow(14,9) \arrow(14,12) \arrow(14,15) \arrow(15,7) \arrow(15,10)
      \arrow(16,5) \arrow(16,11) \arrow(16,14)  \arrow(17,3) \arrow(17,9) \arrow(17,12) \arrow(17,15)$
\end{tiny} &
$\begin{array}{l}
 M_K(t,\,u)_{2}=M_K(u,\,t)_{1}\\
 M_K(t,\,u)_{4}=M_K(u,\,t)_{3}\\
 M_K(t,\,u)_{7}=M_K(u,\,t)_{6} \\
 M_K(t,\,u)_{8}=M_K(u,\,t)_{5}\\
 M_K(t,\,u)_{11}=M_K(u,\,t)_{10}\\
 M_K(t,\,u)_{14}=M_K(u,\,t)_{13} \\
 M_K(t,\,u)_{17}=M_K(u,\,t)_{16}
\end{array}$
\\
\end{tabular}
\end{center}
\end{figure}

$$\begin{array}{c}
   D_K(t)=\left( t-1
 \right) ^{3} \left( {t}^{2}+t+1 \right) ^{3}  \left( {t}^{2}+1 \right)  \left( {t}^
{4}-{t}^{2}+1 \right)  \left(
{t}^{6}+{t}^{5}+{t}^{4}+{t}^{3}+{t}^{2}+ t+1 \right)\\
      \left(
{t}^{12}-{t}^{11}+{t}^{9}-{t}^{8}+{t}^{6}-{t}^{4}+ {t}^{3}-t+1
\right)  \left( t+1 \right) ^{2} \left( {t}^{2}-t+1
 \right) ^{2}
  \end{array}$$
$\widehat{P}_K(t)= \frac {-{t}^{36}-{t}^{18}-1}{ \left(
t-1 \right) ^{3} \left( {t}^{2}+ t+1 \right) ^{3} \left( {t}^{2}+1
\right)  \left( {t}^{4}-{t}^{2}+1
 \right)  \left( {t}^{6}+{t}^{5}+{t}^{4}+{t}^{3}+{t}^{2}+t+1 \right)
 \left( {t}^{12}-{t}^{11}+{t}^{9}-{t}^{8}+{t}^{6}-{t}^{4}+{t}^{3}-t+1
 \right)  \left( t+1 \right) ^{2} \left( {t}^{2}-t+1 \right) ^{2}}$

\subsection{\textsf{Type $L$}}

\noindent The group of type $L$ is the group $\langle S,\ U,\ T,\ V\rangle$, with $S,\ U,\ T$ as for the type $H$, and
$$S:=\left(
       \begin{array}{ccc}
         1 & 0 & 0 \\
         0 & \zeta_{5}^4 & 0 \\
         0 & 0 & \zeta_5 \\
       \end{array}
     \right),\
     U:=\left(
       \begin{array}{ccc}
         -1 & 0 & 0 \\
         0 & 0 & -1 \\
         0 & -1 & 0 \\
       \end{array}
     \right),\
     T:=\frac{1}{\sqrt{5}}\left(
                            \begin{array}{ccc}
                              1 & 1 & 1 \\
                              2 & s & t \\
                              2 & t & s \\
                            \end{array}
                          \right),\
                          V:=\frac{1}{\sqrt{5}}\left(
                               \begin{array}{ccc}
                                 1 & \lambda_1 & \lambda_1 \\
                                 2\lambda_2 & s & t \\
                                 2\lambda_2 & t & s \\
                               \end{array}
                             \right),
$$
where $s:=\zeta_5^2+\zeta_5^3$, $t:=\zeta_5+\zeta_5^4$, $\lambda_1:=\frac{-1+i\sqrt{15}}{4}$, and $\lambda_2:=\frac{-1-i\sqrt{15}}{4}$. Here $l+1=17$, $\textrm{rank}(A^{(1)})=15$,\\ $p=3$, and $\tau_0:=s_0s_5s_6s_{9}s_{10}s_{11}s_{14},\ \tau_1:=s_1s_2s_{7}s_{12}s_{15},\ \tau_2:=s_3s_4s_8s_{13}s_{16}.$\\
\begin{small}
$\Theta=( 3,\
-\zeta_{3},\ -\zeta_{3}^2,\ 3\zeta_{3}^2,\ 3\zeta_{3},\ -1,\ \zeta_{3}^2,\ \zeta_{3},\ 1,\
-\zeta_{15}^2-\zeta_{15}^8,\
  -\zeta_{15}-\zeta_{15}^4,
  -\zeta_{15}^{11}-\zeta_{15}^{14},\ -\zeta_{15}^7-\zeta_{15}^{13},\
  -\zeta_{5}^2-\zeta_{5}^3,\\
  -\zeta_{5}-\zeta_{5}^4,\ 0,\ 0 ).$
\end{small}
\\

\begin{figure}[h]
\begin{center}
\begin{tabular}{ccc}
 \includegraphics[width=4cm,height=4cm]{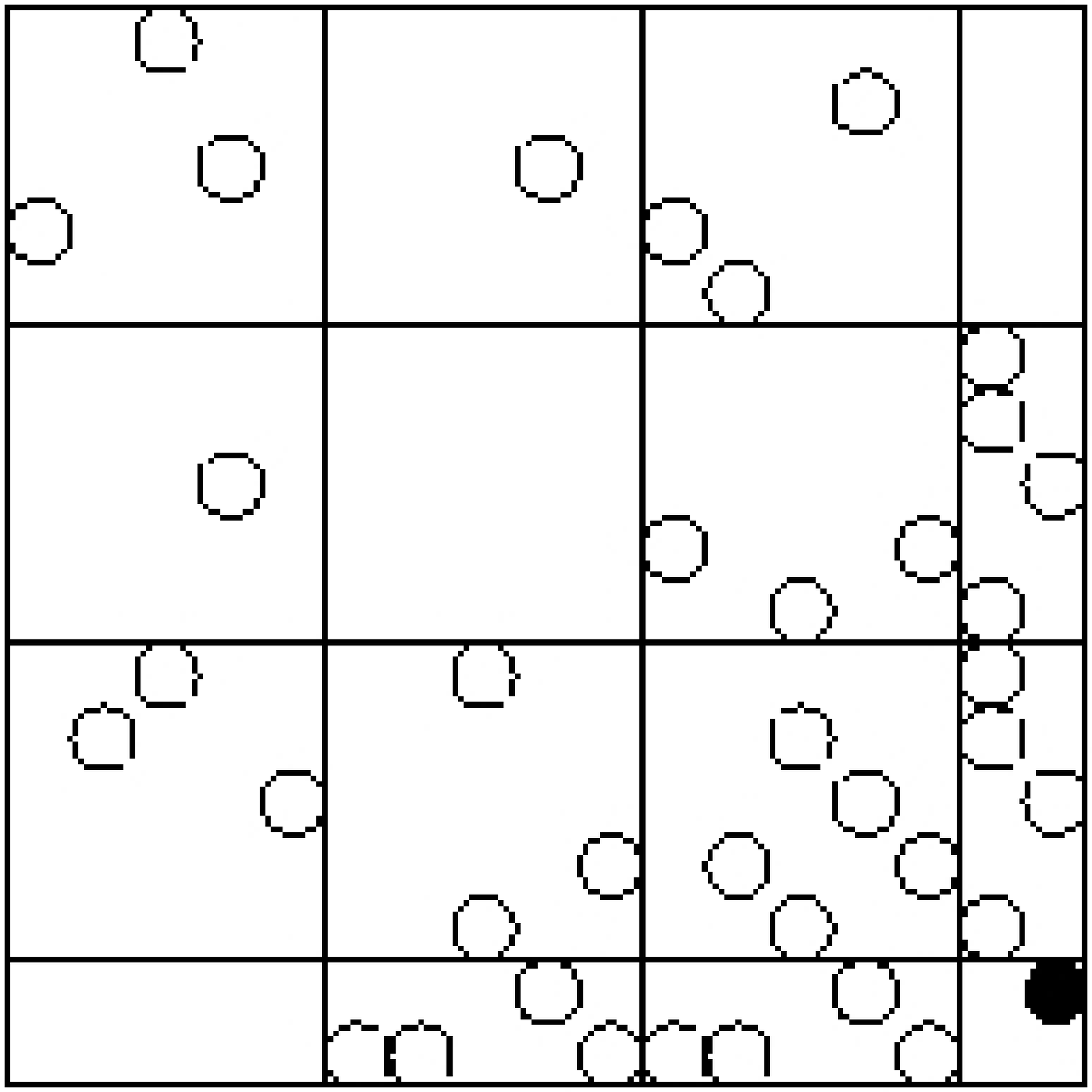} \\[-5cm]& &
\begin{tiny}
\def\arrow(#1,#2){\ncline{-}{#1}{#2}}
$\begin{array}{c@{\hskip .1cm}c@{\hskip .1cm}c@{\hskip .1cm}c@{\hskip .1cm}c@{\hskip .1cm}c@{\hskip .1cm}c}
     & & & \rnode{0}{\Tcircle{0}}  \\[.15cm]
      & & \rnode{2}{\Tcircle{2}} &   &  \rnode{3}{\Tcircle{3}} \\[.15cm]
       & \rnode{8}{\Tcircle{8}} & &  \rnode{10}{\Tcircle{10}} & & \rnode{7}{\Tcircle{7}}  \\[.35cm]
      & &  & \rnode{14}{\Tcircle{14}} \\[.15cm]
           \\[.15cm]
       & \rnode{16}{\Tcircle{16}} & & & & \rnode{15}{\Tcircle{15}}  \\[.15cm]
        \\[.15cm]
        \rnode{6}{\Tcircle{6}} & & & & & & \rnode{5}{\Tcircle{5}} \\[.15cm]
      & & & \rnode{9}{\Tcircle{9}} \\[.15cm]
      & &   \rnode{12}{\Tcircle{12}} & & \rnode{13}{\Tcircle{13}} \\[.15cm]
            & & & \rnode{11}{\Tcircle{11}} \\[.15cm]
       & & \rnode{4}{\Tcircle{4}} & & \rnode{1}{\Tcircle{1}}\\[.15cm]
    \end{array}
    \everypsbox{\scriptstyle}
    \psset{nodesep=0pt,arm=.6,linearc=.4,angleA=0,angleB=90}
    \ncarc[arcangleA=30, arcangleB=10]{-}{7}{16}      \ncarc[arcangleA=-30, arcangleB=-10]{-}{8}{15}
     \arrow(0,3) \arrow(1,11)   \arrow(2,0) \arrow(2,10) \arrow(3,2) \arrow(3,7) \arrow(4,12) \arrow(5,16)   \arrow(6,16) \arrow(7,10)
     \arrow(7,14) \arrow(8,2)    \arrow(9,13) \arrow(9,16) \arrow(10,3) \arrow(10,8)   \arrow(10,16) \arrow(11,4) \arrow(11,13)
     \arrow(11,16) \arrow(12,9)   \arrow(12,11) \arrow(12,14) \arrow(13,1)  \arrow(13,12)    \arrow(13,15) \arrow(14,8)
     \arrow(14,13) \arrow(14,16)   \arrow(15,5) \arrow(15,6) \arrow(15,9) \arrow(15,10) \arrow(15,11) \arrow(15,14)
     \arrow(16,12) \arrow(16,15)
     \ncline[doubleline=true]{-}{15}{16}$
\end{tiny}
$\begin{array}{l}
   M_L(t,\,u)_{3}=M_L(u,\,t)_{2}\\
   M_L(t,\,u)_{5}=M_L(t,\,u)_{6}\\
   M_L(t,\,u)_{8}=M_L(u,\,t)_{7} \\
    M_L(t,\,u)_{13}=M_L(u,\,t)_{12}\\
    M_L(t,\,u)_{16}=M_L(u,\,t)_{15}
 \end{array}
$\\
\end{tabular}
\end{center}
\end{figure}

$$\begin{array}{c}
D_L(t)= \left( t-1 \right) ^{3} \left( {t}^{2}+t+1
\right) ^{3}
 \left( {t}^{2}+1 \right)  \left( {t}^{4}-{t}^{2}+1 \right)  \left( {t
}^{4}+{t}^{3}+{t}^{2}+t+1 \right) \\
  \left(
{t}^{8}-{t}^{7}+{t}^{5}-{t}^ {4}+{t}^{3}-t+1 \right) \left( t+1
\right) ^{2} \left( {t}^{2}-t+1
 \right) ^{2}
  \end{array}
$$
\begin{center}
$\widehat{P}_L(t)=
  \frac {-{t}^{30}+{t}^{15}-1}{ \left( t-1
\right) ^{3} \left( {t}^{2}+ t+1 \right) ^{3} \left( {t}^{2}+1
\right)  \left( {t}^{4}-{t}^{2}+1
 \right)  \left( {t}^{4}+{t}^{3}+{t}^{2}+t+1 \right)  \left( {t}^{8}-{
t}^{7}+{t}^{5}-{t}^{4}+{t}^{3}-t+1 \right)  \left( t+1 \right)
^{2}
 \left( {t}^{2}-t+1 \right) ^{2}}$
 \end{center}

\begin{small}

\end{small}


\begin{thebibliography}{99}

\bibitem[BKR01]{BKR01}
Bridgeland T.,  King A., Reid M., \textit{ The McKay correspondence as an equivalence of derived categories }
\textit{J. Amer. Math. Soc.} \textbf{14} (2001), 535--554.

\bibitem[Bl17]{Bl17}
Blichfeldt H. F., \textit{Finite collineation groups},
 The Univ. Chicago Press, Chicago, 1917.

\bibitem[DHZ05]{DHZ05}
Dais D. I., Henk M., Ziegler G. M., On the Existence of Crepant
Resolutions of Gorenstein Abelian Quotient Singularities in
Dimensions $4$, arXiv:math/0512619v2 [math.AG], 2006.

\bibitem[GNS04]{GNS04}
Gomi Y., Nakamura I., Shinoda K., Coinvariant Algebras of Finite
Subgroups of $\mathbf{SL}_3\mathbb{C}$,
 {\it Canad. J. Math.} \textbf{56} (3), 495--528, 2004.

\bibitem[GSV83]{GSV83}
Gonzalez-Sprinberg G., Verdier J.-L., Construction g\'{e}om\'{e}trique de la correspondance de McKay,
 {\it Annales scientifiques de l'E. N. S., 4\`{e}me s\'{e}rie} \textbf{16} (3), 409--449, 1983.

\bibitem[HH98]{HH98}
Hanany A., He Y.-H., Non-Abelian Gauge Theories,
 arXiv:hep-th/9811183v3, 1998.

\bibitem[I94]{I94}
Ito Y., Crepant Resolution of Trihedral Singularities,
arXiv:alg-geom/9404008v1.

\bibitem[Ka]{Ka}
Kac V. G., Infinite Dimensional Lie Algebras, Bikha\"{a}ser.

\bibitem[Kos85]{Kos85}
Kostant B., The McKay Correspondence, the Coxeter Element and
Representation Theory, {\it SMF, Ast\'erisque, hors s\'erie}, 209--255, 1985.

\bibitem[Kos06]{Kos06}
Kostant B., The Coxeter element and the branching law for the finite subgroups of $\mathbf{SU}(2)$,
The Coxeter legacy, 63--70, Amer. Math. Soc., Providence, RI, 2006, and arXiv:math/0411142v1 [math.RT], 2004.

\bibitem[McK99]{McK99}
McKay, Semi-Affine Coxeter-Dynkin Graphs and $G\subseteq\mathbf{SU}_2\mathbb{\mathbb{C}}$,
 {\it Canad. J. Math.} \textbf{51} (6), 1226--1229, 1999.

\bibitem[Sp77]{Sp77}
Springer T.A., Invariant theory,
 {\it Lecture Notes in Math.} \textbf{585}, Springer-Verlag, 1977.

\bibitem[St08]{St08}
Stekolshchik R., Notes on Coxeter Transformations and the McKay Correspondence,
Springer-Verlag, 2008.

\bibitem[YY93]{YY93}
Yau S. S.-T., Yu Y., Gorenstein Quotient Singularities in
Dimension Three,
 {\it Mem. AMS} \textbf{105}, (505), 1993.

\end{thebibliography}
\end{document}